\documentclass[]{amsart}

\usepackage[pdfencoding=auto, psdextra]{hyperref}
\hypersetup{
    colorlinks=true,
    linkcolor=black,
    filecolor=black, 
    urlcolor=black 
    pdftitle={Overleaf Example},
    pdfpagemode=FullScreen,
    urlcolor = black,
    }
\usepackage{bookmark}
\usepackage[english]{babel} 
\usepackage[T1]{fontenc}

\usepackage{amsthm}
\usepackage{amsmath}
\usepackage[a4paper]{geometry}
\usepackage{diagbox}

\usepackage{verbatim} 
\usepackage[colorinlistoftodos]{todonotes} 
\usepackage{mathtools}
\usepackage{calrsfs}
\usepackage[all]{xy}
\usepackage{amssymb,mathrsfs,amscd,array,multirow,enumerate,amsfonts,euscript} 
\usepackage{comment}
\usepackage{upgreek} 

\xyoption{poly} 
\usepackage{url}

\allowdisplaybreaks
\usepackage{color}
\usepackage{xcolor}

\usepackage{tikz-cd}
\usepackage{tikz}
\usetikzlibrary{matrix,arrows}

\theoremstyle{plain}
\newtheorem{thm}{Theorem}[section]
\newtheorem{lem}[thm]{Lemma}
\newtheorem{prop}[thm]{Proposition}
\newtheorem{cor}[thm]{Corollary}

\theoremstyle{definition}
\newtheorem{defn}[thm]{Definition}

\theoremstyle{remark}
\newtheorem{rmk}[thm]{Remark}
\newtheorem{ex}[thm]{Example}
\numberwithin{equation}{section}

\newtheorem{theor}{\rm{\textbf{Theorem}}}
\newtheorem*{theorem*}{\rm{\textbf{Theorem}}}

\def\parref#1{\autoref{#1}}
\def\thmref#1{Theorem\parref{#1}}
\def\propref#1{Proposition\parref{#1}}
\def\corref#1{Corollary\parref{#1}}
\def\remref#1{Remark\parref{#1}}
\def\lemref#1{Lemma\parref{#1}}
\def\exref#1{Example\parref{#1}}
\def\secref#1{Section~\ref{#1}}
\def\defref#1{Definition\parref{#1}}

\def\makeop#1{\expandafter\def\csname#1\endcsname{\mathop{\rm#1}\nolimits}\ignorespaces}
\makeop{Gal}    \makeop{Br}

\def\makebb#1{\expandfter\def\csname bb#1\endcsname{\mathbb{\rm{#1}}}\ignorespaces}
\def\makebf#1{\expandfter\def\csname bf#1\endcsname{\bf{#1}}\ignorespaces}
\def\mathscr#1{\expandfter\def\csname scr#1\endcsname{{\EuScript{#1}}}\ignorespaces}
\def\makecal#1{\expandfter\def\csname cal#1\endcsname{\mathcal{#1}}\ignorespaces}

\def\qed{\qedmark\medbreak}
\def\qedmark{{\enspace\vrule height 6pt width 5 pt depth 1.5pt}}

\newcommand{\FF}{\mathbb{F}} 
\newcommand{\ZZ}{\mathbb{Z}}
\newcommand{\End}{\mathrm{End}}

\newcommand{\QQ}{\mathbb{Q}}

\newcommand{\Hom}{\mathrm{Hom}}
\newcommand{\Aut}{\text{Aut}}
\makeatletter
\newcommand*{\rom}[1]{\expandafter\@slowromancap\romannumeral #1@}

\def\embed{\hookrightarrow}

\def\to{\rightarrow}

\usepackage{fancyhdr}
\let\@mkboth\@gobbletwo
\let\@oddhead\@empty
\let\@evenhead\@empty
\makeatother

\usepackage{marginnote}
\usepackage{float}
\title{CM-liftability of simple superspecial abelian surfaces over prime fields}
\author{Hsin-Yi Yang}
\address{University of Amsterdam}
\email{h.y.yang@uva.nl}

\date{\today}
\begin{document}
\makeatletter
\@namedef{subjclassname@2020}{\textup{2020} Mathematics Subject Classification}
\begin{abstract} 
For any prime~$p>0$, we prove that simple superspecial abelian surfaces over~$\FF_{p}$ admit CM liftings, 
by using the residual reflex condition (RRC) and Lie types.
   The CM-liftability of ordinary simple abelian surfaces is proved by Serre-Tate,  and the CM-liftability of almost ordinary simple abelian surfaces is proved by Oswal-Shankar and Bergstr\"om-Karemaker-Marseglia, respectively.
    As there can only be ordinary, almost ordinary, or supersingular simple abelian surfaces over~$\FF_{p}$, our work is another step to complete the CM-liftability of simple abelian surfaces over~$\FF_{p}$.

  
\end{abstract}

\thanks{\text{Keywords.} abelian surfaces over prime fields, complex multiplications, CM-liftability.}
\maketitle


\section{Introduction.}
Let~$A$ be a simple abelian variety over a field~$K$ of dimension~$g>0$.
Let~$L/\QQ$ be a CM field with~$[L:\QQ]=2g$. 
It is an important property that~$A/K$ admits \emph{sufficiently many complex multiplications (smCM) by~$L$}.
That is, there exists a~$\QQ$-algebra homomorphism\begin{equation} \label{eq:introL}
	\iota: L\embed \End^{0}_{K}(A):=\End_{K}(A)\otimes_{\ZZ}\QQ,
\end{equation}which is a~$\QQ$-algebra injection from the CM field~$L$ into the endomorphism algebra~$\End^{0}_{K}(A)$ of~$A$.
We also call~$\iota$ a \emph{CM structure} of~$A/K$.
There are many consequences of~$A$ admitting smCM when the field of definition~$K$ of~$A$ is of characteristic~$0$.
We call such an abelian variety \emph{of CM type}.
For example, if~$A$ is defined over~$K:=\mathbb{C}$ and admits smCM by~$L$, then~$A$ can be defined over a number field~$E$, see \cite[Theorem~1]{O73} or \cite[Chapter III, Proposition~26]{Sh06}.
That is, there exists a simple abelian variety~$B/E$, which may not be unique, admitting smCM by~$L$ with CM structure~$\iota'$, and a~$K$-isomorphism such that\[(A,\iota)\cong (B,\iota')\otimes_{E}K.\] 
The~$K$-isomorphism is called \emph{$L$-linear with respect to the CM structures~$\iota,\iota'$}.

By a reduction of~$A/\overline{\QQ}$, we mean a reduction of some~$B/E$ as above.
In particular, as the abelian variety~$B/E$ admits smCM, it has potentially good reduction by Serre and Tate, see \cite[Theorem~6]{ST68}.
Let~$O_{L}$ be the ring of integers in~$L$.
Then, the Frobenius endomorphism of a reduction of~$A$ is induced by an element~$\pi$ in~$O_{L}$ by \cite[Theorem A]{H68}.
Let~$p>0$ be a prime, and let~$q$ be a power of~$p$.
A \emph{Weil-$q$ number} (of weight~$1$) is an algebraic integer having absolute value~$q^{1/2}$ via any embedding to~$\mathbb{C}$. 
An ideal factorization of~$\pi O_{L}$ as in \textit{loc. cit.} and abelian varieties admitting smCM by~$L$ are used to prove the
so-called Honda-Tate theorem, which we recall in the following.
\begin{theorem*}(Honda-Tate) \cite{T66}, \cite{H68}, \cite{T06}
There exists a bijection, given by assigning a simple abelian variety over~$\FF_{q}$ to its Frobenius endomorphism, from the set of isogeny classes of simple abelian varieties over a finite field~$\FF_{q}$ to the set of~$\Gal(\overline{\QQ}/\QQ)$-conjugacy classes of Weil-$q$ numbers.
\end{theorem*}



From now on, we focus on the case when~$A$ is defined over a finite field~$K:=\FF_{q}$.
Tate showed that it has smCM by a CM field~$L$ in \cite{T66}.
Let~$(R,\mathfrak{m})$ be a pair, where~$R$ is a local domain with characteristic~$0$, and~$\mathfrak{m}$ is the maximal ideal of~$R$, such that~$R/\mathfrak{m}$ is isomorphic to~$\FF_{q}$ (e.g.~$R$ is a discrete valuation ring).
Oort asked whether~$A/\FF_{q}$ has a \emph{CM lifting} (CML) in \cite{F92}, that is, if there exist such a pair~$(R,\mathfrak{m})$ and an abelian scheme~$\mathfrak{A}$ over~$R$ of relative dimension~$g>0$ equipped with a~$\QQ$-algebra injection~$L \hookrightarrow \End^{0}_{R}(\mathfrak{A})$, 
such that the special fiber of~$\mathfrak{A}$ at~$\mathfrak{m}$ is~$\FF_{q}$-isomorphic to~$A/\FF_{q}$.

Let~$\overline{\FF}_{p}$ be an algebraically closed field containing~$\FF_{q}$.
In \cite[Theorem~A]{F92}, Oort showed that~$A\otimes\overline{\FF}_{p}$ admits a CM lifting if~the~$p$-rank of~$A\otimes\overline{\FF}_{p}$ is~$g$ or~$g-1$, where~$g>0$ is the dimension of~$A$. 
(An integer~$f$ is called the \emph{$p$-rank} of~$A\otimes\overline{\FF}_{p}$ if the~$p$-torsion group~$A(\overline{\FF}_{p})[p]$ is isomorphic to~$(\ZZ/p\ZZ)^{f}$.)
However, in general, suppose that~$A/\FF_{q}$ admits a \emph{formal deformation} to a formal abelian scheme~$\mathfrak{A}$ over~$W(\FF_{q})$, where~$W(\FF_{q})$ is the Witt ring of~$\FF_{q}$-vectors. 
If we demand a CM structure of~$A/\FF_{q}$ to deform, 
then~$\mathfrak{A}$ may not be algebraizable, that is, it may not be the formal completion of a proper scheme over~$W(\FF_{q})$, see \cite[Example~1.4.5.4]{CCO14}. 
Therefore, not all abelian varieties over finite fields admit CM liftings, even if they admit deformations to formal abelian schemes over some Witt rings.
In fact, Oort showed that there exist infinitely many abelian varieties over~$\overline{\FF}_{p}$, which are of dimension~$g\geq 3$ and~$p$-rank~$f$ such that~$0\leq f \leq g-2$, that do not admit CM liftings in \cite[Theorem~B]{F92}. 

The CM-liftability problem is solved by Serre and Tate when the abelian varieties are ordinary, see \cite[Chapter V, Theorem~(3.3)]{Messing72}.
Moreover, as is shown in \cite[Proposition~2.3]{OS20} by using Grothendieck-Messing theory, any one-dimensional supersingular~$p$-divisible group over~$\overline{\FF}_{p}$, whose endomorphism ring contains an integrally closed rank two~$\ZZ_{p}$-algebra, and its endomorphism ring lift to~$W(\overline{\FF}_{p})$.
Therefore, the CM-liftability of almost ordinary abelian varieties over~$\overline{\FF}_{p}$ is clear after applying the Serre-Tate lifting equivalence, see \cite[Chapter V, Theorem~(2.3)]{Messing72}.
Additionally, Bergstr\"om-Karemaker-Marseglia proved that an almost ordinary abelian variety over~$\FF_{q}$ and its endomorphism ring lift if its endomorphism ring is abelian in \cite[Proposition~2.10]{BKM}.
On the other hand, it is proven in \cite[Corollary~3.7]{K21} that there are only finitely many supersingular principally polarized abelian varieties over~$\overline{\FF}_{p}$  with a given dimension admitting CM liftings.

In the following, we introduce our result.
Let~$\pi$ be a Weil-$p$ number of our concern, that is,
\begin{equation*}
\pi=\begin{cases}
	\sqrt{p}\zeta_{3},\sqrt{p}\zeta_{8},\sqrt{p}\zeta_{12}, \,\textrm{or}\,\pm\sqrt{5}\zeta_{5} \,\,\textrm{when}~p\geq 5,\\
	\sqrt{3}\zeta_{3},\,\textrm{or}\,\sqrt{3}\zeta_{8} \qquad\quad\quad\quad\quad\,\,\,\,\,\,\,\,\textrm{when}~p=3,\\
	\sqrt{2}\zeta_{3},\sqrt{2}\zeta_{12},\,\textrm{or}\pm\sqrt{2}\zeta_{24} \qquad\,\,\,\,\,\textrm{when}~p=2.
\end{cases}
\end{equation*}
Then,~$\pi$ corresponds to a simple \emph{supersingular} abelian surface~$A_{\pi}$ over~$\FF_{p}$ by the Honda-Tate theorem.
Moreover, these are all Weil-$p$ numbers corresponding to simple supersingular abelian surfaces over~$\FF_{p}$, see \cite[Section~5.2]{XYY16}.
An abelian variety~$A$ over a finite field~$\FF_{q}$ is said to be \emph{superspecial} if~$A\otimes_{\FF_{q}}\overline{\FF}_{p}$ is~$\overline{\FF}_{p}$-isomorphic to a product of supersingular elliptic curves over~$\overline{\FF}_{p}$.
Our result states that all simple superspecial abelian surfaces~$A_{\pi}$ over prime fields, in the isogeny class determined by~$\pi$, have CM liftings.
 
\begin{theor}(\thmref{cor:surfacesCML}, \propref{rem:p=3})
	Let~$\pi$ be a Weil-$p$ number of our concern. 
	Let~$L:=\QQ(\pi)$.
	Let~$A_{\pi}/\FF_{p}$ be any simple superspecial abelian surface in the isogeny class determined by~$\pi$, which is obtained by the Honda-Tate theorem.
	Then,~$A_{\pi}/\FF_{p}$ has a CML, admitting smCM by~$L$.
\end{theor}
We treat the case when~$\pi=\sqrt{3}\zeta_{3}$ separately in \propref{rem:p=3}.
This case is more complicated because the index~$[O_{L}:\End_{\FF_{3}}(A_{\pi})]$ can be~$3$, which means that~$\End_{\FF_{3}}(A_{\pi})$ is not locally isomorphic to~$O_{L}$ at~$3$.

In addition, we have an analogous result when the Weil-$p$ number is~$\sqrt{p}$.
Any simple abelian surface~$A_{\pi}$ over~$\FF_{p}$, obtained by the Honda-Tate theorem, in the isogeny class determined by~$\pi:=\sqrt{p}$ is superspecial, see \cite[p.~528]{W69}.
It is different from the cases considered above, because the endomorphism algebra~$\End_{\FF_{p}}^{0}(A_{\pi})$ is a quaternion algebra over~$\QQ(\sqrt{p})$.
This means that we can embed more than one CM field into it.
However, we find a CM field and prove the same result:
\begin{theor}(\thmref{thm:liftingreal})
Let~$L:=\QQ(\sqrt{p}\zeta_{3})$ be a CM field over~$\QQ$.
In the isogeny class determined by~$\pi=\sqrt{p}$, every simple superspecial abelian surface over~$\FF_{p}$ has a CML, admitting smCM by~$L$.
\end{theor}


In the remainder of the introduction, we mention the tools used in the proofs of our results. 
We combine the residual reflex condition (RRC), Lie types, \cite[Proposition~2.22]{BKM}, \cite[Proposition~4.2.9]{CCO14}, Grothendieck-Messing theory (see \cite[Chapter V, Theorem~(1.6) or~ Theorem~(1.10)]{Messing72}), and Serre-Tate deformation theory (see \cite[Theorem~1.2.1]{K81} or \cite[Chapter V, Theorem~(2.3)]{Messing72}) in our proofs. 
Our strategy of the proof of the theorems is roughly as followings. 
Let~$A_{\pi}$ be a simple superspecial abelian surface over~$\FF_{p}$ in the isogeny class determined by a Weil-$p$ number~$\pi$ of our concern. 
We check whether~$A_{\pi}$ satisfies RRC over~$\FF_{p}$.   
By \cite[Theorem~2.5.3]{CCO14}, this implies that~$A_{\pi}$ has a CM lifting, up to~$L$-linear isogeny.
However, it does not imply that~$A_{\pi}$ has a CM lifting without an~$L$-isogeny. 
To show the CM liftability of~$A_{\pi}/\FF_{p}$, we want to apply \cite[Proposition~2.22]{BKM}, which says the following: Let~$A_{0},A_{0}'$ be abelian varieties over~$\FF_{q}$ with CM by a CM field~$L$. 
Assume that~$A_{0}$ admits a CM lifting to a local domain~$R$ of characteristic~$0$ with residue field~$\FF_{q}$. 
Assume that there exists an~$\FF_{q}$-separable isogeny~$A_{0}\to A_{0}'$. 
Then,~$A_{0}'$ also admits a CM lifting to~$R$. 
So, it is sufficient to find a \emph{separable} isogeny between~$A$ and the~$L$-isogenous abelian variety of \cite[Theorem~2.5.3]{CCO14} if~$A$ satisfies RRC.

We use Lie types to find such an isogeny. 
Composing the CM structure~$\iota$ in \eqref{eq:introL} and the canonical embedding of endomorphism rings~$\End_{\FF_{p}}(A_{\pi})\embed\End_{\overline{\FF}_{p}}(A_{\pi}\otimes_{\FF_{p}}\overline{\FF}_{p})$, we have the following embedding
\begin{equation}\label{eq:embedintro}
\iota: O_{L}\otimes_{\ZZ}\ZZ_{p}\embed \End_{\overline{\FF}_{p}}(A_{\pi}\otimes\overline{\FF}_{p})\otimes_{\ZZ}\ZZ_{p},
\end{equation}which we denote by~$\iota$ again with slight ambiguity.
In particular, \eqref{eq:embedintro} induces an embedding from~$O_{L}\otimes\ZZ_{p}$ into~$\End_{\overline{\FF}_{p}}((A_{\pi}\otimes\overline{\FF}_{p})[p^{\infty}])$. 
Such a~$p$-divisible group is said to be \emph{$O_{L}$-linear}. 
Note from \eqref{eq:embedintro} that one may do base change for computing Lie types.  
Over~$\overline{\FF}_{p}$, Yu used \emph{Lie types}, see \defref{def:glt}
, to classify~$O_{L}$-linear CM~$p$-divisible groups of height~$[L:\QQ]$ (up to~$L$-linear~$\overline{\FF}_{p}$-isomorphism) in \cite{Yu04}.
We do not apply Yu's result because we want the CM liftability of~$A_{\pi}$ without base changing and because~$O_{L}$ is assumed to be embedded into the endomorphism ring of an abelian variety in \textit{loc. cit}. 

In \propref{lem:sameLietype}, we compute the Lie type of~$A_{\pi}$, in the isogeny class determined by~$\pi$. 
Then, it follows from \cite[Proposition~4.2.9]{CCO14}, which provides an~$O_{L}\otimes\ZZ_{p}$-linear isomorphism~$\xi$ between two isomorphism classes over~$\FF_{p}$ in the isogeny class, and from \cite[Main theorem]{T66} that there exists a separable isogeny between~$A_{\pi}$ and its~$L$-isogenous abelian surface. 
Since~$\xi$ is an isomorphism of Dieudonn\'e modules, there exists a sequence of homoomorphisms of corresponding abelian surfaces over~$\FF_{p}$ appoaching~$\xi$ as explained in the proof of \cite[Section~1, Theorem~7]{MW71}.
In particular, as the element in the sequence is closed enough, it is an isogeny, and it induces an isomorphism between Dieudonn\'e modules.
Although we sometimes obtain the CM liftability of~$A_{\pi}$ after base changing, we use \thmref{lem:changeofcoe} to show that~$A_{\pi}$ over~$\FF_{p}$ has a CM lifting.
In the proof of \thmref{lem:changeofcoe}, the deformation theory is used. 

The following is the outline of the paper.
In \secref{chap:CMandG}, we set up notations and recall the reflex residual condition (RRC).
In \secref{sec:deformation}, we recall deformations. 
In \secref{sec:Lie type}, we recall definitions and constructions of Lie types, and show that any two simple superspecical abelian surfaces over~$\FF_{p}$, in an isogeny class determined by a Weil-$p$ number of our concern, have the same Lie type. 
In \secref{ex:1}, we work out an example of the CM liftability explicitly with the Weil-$7$ number~$\sqrt{7}\zeta_{3}$.
In \secref{chap:ss}, we prove in \thmref{cor:surfacesCML}, \propref{rem:p=3} and \thmref{thm:liftingreal} that a simple superspecial abelian surface over~$\FF_{p}$, in the isogeny class determined by~$\pi$, has a CM lifting for any Weil-$p$ number~$\pi$ of our concern or~$\sqrt{p}$ and for any prime~$p>0$.

\section{CM types.} \label{chap:CMandG}
In this section, we recall definitions and set notations, which will be used throughout the paper.
When we write an algebraically closed field~$\overline{\FF}_{p}$ or~$\overline{\QQ}$, we always mean the one containing the field of definition of any abelian variety or the one containing the CM field mentioned in the corresponding context.

\subsection{CM types (algebraic).}

We recall some definitions related to CM types of a CM field.

\begin{defn}\label{def:CMtype}Let~$L/\QQ$ be a CM field with maximal totally real subfield~$L_{0}$. 
Let~$c\in\Aut(L/\QQ)$ be the non-trivial element of order~$2$, which is induced by the complex conjugation for any embedding~$L\embed \mathbb{C}$.
Let~$\Phi$ be a subset in~$\Hom(L,\overline{\QQ})$ and~$\overline{\Phi}:= \{\overline{\varphi}:=\varphi\circ c\,|\,\varphi\in\Phi\}$. 
The subset~$\Phi$ is a \emph{CM type of~$L$} if~$\Phi \cap \overline{\Phi} =\emptyset$ and~$\Hom(L,\overline{\QQ})=\Phi\sqcup\overline{\Phi}$.
We also call~$(L,\Phi)$ a \emph{CM type of~$L$} if~$\Phi$ is.\end{defn} 
Denote~$\Sigma_{(\cdot),p}$ the set of all places in a number field~$(\cdot)$ lying above~$p$.
\begin{defn}\label{def:emb}
Fix an embedding~$\overline{\QQ}_{}\embed\overline{\QQ}_{p}$. 
Then, we may identify~$\Hom(L,\overline{\QQ}_{p})$ and~$\Hom(L,\overline{\QQ})$. 
Let~$L$ be a CM field, and~$\Phi$ be one of its CM types in the sense of \defref{def:CMtype}. 
\begin{enumerate}
	\item We specifically call~$\Phi$ a \emph{$\overline{\QQ}_{p}$-valued CM type of~$L$} when we regard it as a subset of~$\Hom(L,\overline{\QQ}_{p})$. 
  \item For any~$w\in\Sigma_{L,p}$, the subset~$\Phi_{w}$ of the~$\overline{\QQ}_{p}$-valued CM type~$\Phi$ of~$L$ consists of elements in~$\Phi$, which induce the place~$w$. 
    That is, any element~$\varphi\in\Phi_{w}$ factors through the completion~$L_{w}$ of~$L$ at~$w$ with respect to the fixed embedding~$\overline{\QQ}\embed\overline{\QQ}_{p}$.
    \item The \emph{slope of a CM type~$\Phi$} is the set~$\{|\Phi_{w}|/[L_{w}:\QQ_{p}]\}_{w\in\Sigma_{L,p}}$.
\end{enumerate}
\end{defn}

By \cite[Chapter II, Proposition 28]{Sh06} or \cite[Chapter 1, Section 5]{L12}, we have the following definition.

\begin{defn}\label{def:reflexfield}
Let~$L/\QQ$ be a CM field of CM type~$(L,\Phi)$.
The \emph{reflex field} of~$(L,\Phi)$ is the subfield~$L':=\QQ(\sum_{\sigma\in \Phi} x^{\sigma}|\,x\in L)$ of~$\overline{\QQ}$.
\end{defn}

\subsection{CM types (abelian varieties)}
In this subsection, we connect CM types of a CM field and abelian varieties.

\begin{defn}\cite[Definition~1.3.1.2]{CCO14}, \cite[Chapter~3, Section~1]{L12} \label{def:avCM}
\begin{enumerate}
	\item 
    Let~$A/K$ be a simple abelian variety of dimension~$g>0$ over a field~$K$.
     If there exists a CM field~$L$ with~$[L:\QQ]=2g$, embedded in~$\End^{0}_{K}(A)$ by an injection~$\iota$ of~$\QQ$-algebras, then we say~$A/K$ admits \emph{sufficiently many complex multiplications (smCM) by~$L$ over~$K$}. 
    We call the embedding~$\iota$ a \emph{CM structure} of~$A/K$. 
Moreover, we say that~\emph{$(A,\iota)$ is defined over~$K$} as~$A$ is defined over~$K$, and the CM structure~$\iota$ of~$A$ embeds~$L$ into~$\End_{K}^{0}(A)$. 
    \item Let~$(A,\iota_{A}),(B,\iota_{B})$ be simple abelian varieties of dimension~$g>0$ defined over~$K$, which admit smCM by CM field~$L$ with respective CM structures~$\iota_{A}, \iota_{B}$ of~$A,B$. 
    They are\emph{~$L$-linearly~$K$-isogenous} if there exists a~$K$-isogeny~$\phi: A\to B$, such that the following diagram \begin{equation}
    	\begin{tikzcd}
    		A\arrow[r,"\phi"]\arrow[d,"\iota_{A}(\ell)"'] & B \arrow[d,"\iota_{B}(\ell)"]\\
    		A \arrow[r,"\phi"'] & B
    	\end{tikzcd}
    \end{equation}commutes for any~$\ell\in L$.
    \item 
    Let~$L$ be a CM field with~$[L:\QQ]=2g>0$, and let~$\Phi\subset \Hom(L,\overline{\QQ})$ be a CM type of~$L$.
    \begin{enumerate} \item 
    Let~$(A,\iota)$ be a simple abelian variety of dimension~$g$ over~$\overline{\QQ}$, admitting smCM by CM field~$L$.
    The abelian variety~$(A,\iota)$ is \emph{of CM type}~$(L,\Phi)$ if the induced action by~$\iota(\ell)$ on the Lie algebra of~$A$ can be represented by the diagonal matrix \begin{equation}\label{rem:related}
\text{diag}(\ell^{\varphi_{1}},\dots,\ell^{\varphi_{g}}),
\end{equation}
where~$\Phi=\{\varphi_{1},\dots, \varphi_{g}\}$, for any~$\ell\in L$.
Moreover, we say that~$\iota$ \emph{is induced by}~$\Phi$.
\item Let~$(A,\iota)$ be a simple abelian variety of dimension~$g$ over a number field~$K$, admitting smCM by CM field~$L$.
The abelian variety~$(A,\iota)$ is \emph{of CM type}~$(L,\Phi)$ if~$(A\otimes_{K}\overline{\QQ},\iota_{\textrm{com}})$ is of CM type~$(L,\Phi)$, where~$\iota_{\textrm{com}}$ is the composition of the CM structure~$\iota$ of~$A/K$ and the canonical embedding~$\End_{K}(A)\embed\End_{\overline{\QQ}}(A\otimes_{K}\overline{\QQ})$.
That is, \[\iota_{\textrm{com}}: L\xhookrightarrow{\iota}\End_{K}^{0}(A)\xhookrightarrow{\textrm{can}}\End_{\overline{\QQ}}^{0}(A\otimes_{K}\overline{\QQ}),\] where~$\End_{K}^{0}(A)\xhookrightarrow{\textrm{can}}\End_{\overline{\QQ}}^{0}(A\otimes_{K}\overline{\QQ})$ is the canonical embedding.
The embedding~$\iota_{\textrm{com}}$ is a CM structure of~$A\otimes_{K}\overline{\QQ}$.
We say such~$\iota$ \emph{is induced by~$\Phi$} if the composition~$\iota_{\textrm{com}}$ is induced by~$\Phi$.
With slight ambiguity, we will denote~$\iota_{\textrm{com}}$ by~$\iota$.
\end{enumerate}
   \end{enumerate}
    \end{defn}

\subsection{Residual reflex condition.}\label{sec:RRC}
Let~$A_{0}/\FF_{q}$ be a simple abelian variety of dimension~$g>0$, admitting smCM by CM field~$L/\QQ$.
In this subsection, we recall the residual reflex condition, which helps to detect if~$A_{0}/\FF_{q}$ has a CM lifting, up to~$L$-isogeny.
To be more precise, if~$A_{0}/\FF_{q}$ satisfies the residual reflex condition, then it has a CM lifting up to~$L$-linear~$\FF_{q}$-isogeny by \cite[Theorem~2.5.3]{CCO14}.
Conversely, if~$A_{0}/\FF_{q}$ admits a CM lifting, then it satisfies the residual reflex condition, see \cite[Section~2.1.5]{CCO14}.

We recall the definition of CM liftings.
\begin{defn}\cite[1.8.5]{CCO14}\label{def:CML}
Let~$\FF_{q}$ be a finite field of characteristic~$p>0$, where~$q$ is a power of~$p$.
Let~$A_{0}/\FF_{q}$ be a simple abelian variety of dimension~$g>0$.
	The abelian variety~$A_{0}/\FF_{q}$ has a \emph{CM lifting (CML)} if there exist a CM field~$L/\QQ$ with degree~$[L:\QQ]=2g$ and a local domain $(R,\mathfrak{m})$ such that~$\textrm{char}(R)=0$ and~$\textrm{char}(R/\mathfrak{m})=p$ (e.g.~$R$ is a discrete valuation ring), and an abelian scheme~$A$ over~$R$ of relative dimension~$g$ equipped with a~$\QQ$-algebra injection~$L \hookrightarrow \End^{0}_{R}(A)$, 
such that the special fiber of~$A$ at~$\mathfrak{m}$ is~$\FF_{q}$-isomorphic to~$A_{0}/\FF_{q}$.
\end{defn}

We recall the Shimura-Taniyama formula.
\begin{defn}\cite[2.1.4.1]{CCO14}\label{def:ST}
Fix an embedding~$\overline{\QQ}\embed\overline{\QQ}_{p}$.
Let~$A_{0}/\FF_{q}$ be a simple abelian variety over a finite field~$\FF_{q}$, admitting smCM by CM field~$L$.
    Let~$\textrm{Frob}_{A_{0},q}$ be the Frobenius endomorphism of~$A_{0}/\FF_{q}$.
    The \emph{Shimura-Taniyama formula} says there exists a~$\overline{\QQ}_{p}$-valued CM type~$\Phi$ of~$L$ such that \begin{equation}\label{ST}
    \frac{\text{ord}_{w}({\rm{Frob}}_{A_{0},q})}{\text{ord}_{w}(q)} = \frac{|\{\phi\in\Phi\,|\,\,\phi\,\, \text{induces}\,\,w\}|}{[L_{w}:\QQ_{p}]}
\end{equation} for all~$w\in \Sigma_{L,p}$.
We call the right-hand side of \eqref{ST}
the \emph{slope} of~$A_{0}/\FF_{q}$ and the left-hand side of \eqref{ST} the \emph{slope} of~$\Phi$.
We also say that a CM type~$\Phi$ of~$L$ \emph{is realized by}~$A_{0}/\FF_{q}$ if the Shimura-Taniyama formula \eqref{ST} holds with~$\Phi$ for all~$w\in\Sigma_{L,p}$. 
\end{defn}
\begin{rmk}
	After fixing the embedding~$\overline{\QQ}\embed\overline{\QQ}_{p}$ as in \defref{def:ST}, we regard a CM type~$\Phi$ of~$L$ as a~$\overline{\QQ}_{p}$-valued CM type of~$L$.
Let $\Phi_{w}:= \{\phi\in\Phi\,|\,\,\phi\,\, \text{induces}\,\,w\}$ as in \defref{def:emb} for any~$w\in\Sigma_{L,p}$. 
Then, we can rewrite the Shimura-Taniyama formula \eqref{ST} to
\begin{equation*}
    \frac{\text{ord}_{w}({\rm{Frob}}_{A_{0},q})}{\text{ord}_{w}(q)} = \frac{|\Phi_{w}|}{[L_{w}:\QQ_{p}]}
\end{equation*}for all~$w\in\Sigma_{L,p}$.
That is, the slope of the CM type~$\Phi$ of~$L$ is the same as the slope of~$A_{0}/\FF_{q}$ if the Shimura-Taniyama formula \eqref{ST} holds with~$\Phi$.
\end{rmk}

\begin{rmk}
Let~$A_{0}/\FF_{q}$ be a simple abelian variety.
Suppose that~$A_{0}/\FF_{q}$ admits smCM by CM field~$L$ with CM structure~$\tilde{\iota}:L\embed\End_{\FF_{q}}^{0}(A_{0})$.
We may define the slope in a slightly general way as in \cite[Section~2.1.4.1]{CCO14}, which we recall as follows.
Let~$F(T)\in\ZZ[T]$ be the characteristic polynomial of the Frobenius endomorphism~${\rm{Frob}}_{A_{0},q}$ of~$A_{0}/\FF_{q}$. 
Fix an embedding~$\overline{\QQ}\hookrightarrow \overline{\QQ}_{p}$.
Let~$\{\lambda_{i}\}_{i=1,\dots,2g}$ be the set of roots (counted with multiplicity) of~$F(T)$ in~$\overline{\QQ}_{p}$.    
Then, the set~$\{\text{ord}_{p}(\lambda_{i})/\text{ord}_{p}(q)\}_{i=1,\dots,2g}$ (counted with multiplicity) is the slope of~$A_{0}/\FF_{q}$, 
where~$\text{ord}_{p}$ is a valuation at~$p$ such that~$\text{ord}_{p}(p):=1$.
Moreover,\[F(T)=\prod_{w\in\Sigma_{L,p}}\prod_{\varphi\in\Hom(L_{w},\overline{\QQ}_{p})}(T-\textrm{Frob}_{A_{0},q}^{\varphi}),\] since~$L\otimes\QQ_{p}=\prod_{w\in\Sigma_{L,p}}L_{w}$ and~$\textrm{Frob}_{A_{0},q}$ is an algebraic interger.
Therefore,~$\{\lambda_{i}\}_{i=1,\dots,2g}=\{\textrm{Frob}_{A_{0},q}^{\varphi}\}_{\varphi\in \Hom(L,\overline{\QQ}_{p})}$, and hence\[\{\text{ord}_{p}(\lambda_{i})/\text{ord}_{p}(q)\}_{i=1,\dots,2g}=\{\text{ord}_{p}(w)\textrm{ord}_{w}(\textrm{Frob}_{A_{0},q})/\text{ord}_{p}(q)\}_{w\in\Sigma_{L,p}},\] where the right-hand side is counted with multiplicity.
We point out that the number of distinct places~$|\Sigma_{L,p}|$ might not equal the number of distinct slopes of~$A_{0}/\FF_{q}$, see \secref{ex:1} for example.
\end{rmk}


We recall the definition of the residual reflex condition.
\begin{defn}\cite[2.1.5]{CCO14}\label{def:RRC}
Let~$A_{0}/\FF_{q}$ be a simple abelian variety, admitting smCM by CM field~$L$.
    Let~$\Phi\subset \Hom(L,\overline{\QQ})$ be a CM type of~$L$.
    We say that~$(L,\Phi)$ satisfies the \emph{residual reflex condition (RRC)} if there exists an embedding~$i:\overline{\QQ}\embed\overline{\QQ}_{p}$, such that~$\Phi$ is~$\overline{\QQ}_{p}$-valued CM type of~$L$, and the following hold.
    \begin{enumerate}
        \item The slope of~$A_{0}$ satisfies the Shimura-Taniyama formula \eqref{ST}.
        \item Let~$\phi:L'\embed\overline{\QQ}$ be the canonical embedding of the reflex field~$L'$, which is a subfield of~$\overline{\QQ}$.
        Let~$v$ be the place of~$L'$, induced by~$i\circ \phi$.
        That is, we have a commutative diagram \begin{equation*}
        	\begin{tikzcd}
        		L' \arrow[rr,"i\circ \phi"]\arrow[rd, hook', "\textrm{incl}"']& &\overline{\QQ}_{p}\\
        		&L'_{v},\arrow[ru,hook']&
        	\end{tikzcd}
        \end{equation*}where the map~$\textrm{incl}$ is the canonical embedding.
        The residue field of~$L'$ at $v$ can be realized as a subfield of~$\FF_{q}$. 
    \end{enumerate}
    We say~$A_{0}/\FF_{q}$ satisfies RRC with~$(L,\Phi)$ if~$(L,\Phi)$ satisfies RRC.
\end{defn}
\begin{rmk}
Keep the notation as in \defref{def:RRC}.
From the definition, if one wants to show that a simple abelian variety~$A_{0}$ over a finite field~$\FF_{q}$ does not satisfy RRC, one may consider every residue field of~$L'$ (as a subfield in~$\overline{\QQ}$) at every place lying above~$p$, and show that none of them can be realized as a subfield of~$\FF_{q}$.
\end{rmk}

\section{Deformations.}\label{sec:deformation}
Let~$A_{0}$ be a simple abelian variety over a finite field~$\FF_{q}$, which admits smCM by a CM field~$L$.
In \thmref{lem:changeofcoe}, we reduce the problem of finding a CM lifting of~$A_{0}/\FF_{p}$ to the problem of finding a CM lifting of~$A_{0}\otimes\overline{\FF}_{p}$, by considering a deformation, under certain situations.

We start with recalling the following definition.
\begin{defn}\cite[Definition~1.4.4.4]{CCO14}\label{def:deformation}
Let~$R$ be a complete local noetherian ring with residue field~$\FF_{q}$.
	Let~$A_{0}/\FF_{q}$ be an abelian variety or a~$p$-divisible group.
	Let~$\alpha_{0}:O\embed\End_{\FF_{q}}(A_{0})$ be an injective homomorphism from a~$\ZZ$-finite associative ring~$O$. 
 A \emph{(resp. formal) deformation of~$(A_{0},\alpha_{0})$ to an (resp.~a formal) abelian scheme} over~$R$ is a pair~$(A,i)$, where~$A$ is an (resp. formal) abelian scheme or a~$p$-divisible group over~$R$ with an injection~$O\embed\End_{R}(A)$, and~$i:A\otimes\FF_{q}\cong A_{0}$ is an~$O$-linear isomorphism over~$\FF_{q}$.
\end{defn}

We compare the definition of a CM lifting and of a deformation in the following remark.
\begin{rmk}\label{rmk:defandCM}
A CM lifting is different from a deformation.
The isomorphism~$i$ in \defref{def:deformation} is~$O$-linear, while only an~$\FF_{q}$-isomorphism is required for a CM lifting, see \defref{def:CML}.
Therefore, the condition of having a deformation is stronger than of having a CM lifting, see \defref{def:CML}.

\end{rmk}

From now on, let~$(A_{0},\alpha_{0})$ be as in \defref{def:deformation}.
We recall the deformation functor with respect to the pair~$(A_{0},\alpha_{0})$ now.
\begin{defn}\cite[Section~2.1]{CCO14}
Let~$\FF_{q}$ be a finite field of characteristic~$p>0$.
Let~$\Lambda:=W(\FF_{q})$ be the ring of Witt vectors with residue field~$\FF_{q}$.
Let~$C_{\Lambda}$ be the category of artinian local~$\Lambda$-algebras~$R$ with local structure map~$\Lambda\to R$ and residue field~$\FF_{q}$.
\begin{enumerate}
\item	Define the functor
\[\textrm{Def}_{\Lambda}(A_{0},\alpha_{0}):C_{\Lambda}\to \textrm{Set},\] sending every~$R\in C_{\Lambda}$ to the set of isomorphism classes of deformations of~$(A_{0},\alpha_{0})$ over~$R$.
\item The functor~$\textrm{Def}_{\Lambda}(A_{0},\alpha_{0})$ is \emph{pro-representable} if there is a complete local noetherian~$\Lambda$-algebra~$\mathcal{R}$ with local structure map~$\Lambda\to \mathcal{R}$ and residue field~$\FF_{q}$ such that~$\textrm{Def}_{\Lambda}(A_{0},\alpha_{0}) \cong \Hom_{\Lambda}(\mathcal{R},\cdot)$ (using local~$\Lambda$-algebra homomorphisms).
	Moreover, such~$\mathcal{R}$ is called a \emph{formal deformation ring}.
\end{enumerate}
\end{defn}



\begin{thm}\cite[Theorem~1.4.4.5, Theorem~1.4.5.5]{CCO14}\label{thm:pro}
Let~$A_{0}$ be an abelian variety (resp.~$p$-divisible group) of dimension~$g$ (resp. height~$h$).
Let~$\Lambda[[x_{1},\cdots, x_{g^2}]]$ be a formal power series ring. 
Then, the deformation functors~$\rm{Def}_{\Lambda}(A_{0},\alpha_{0})$ and~$\rm{Def}_{\Lambda}(A_{0},\alpha_{0},\lambda_{0})$ are pro-representable by quotients of~$\Lambda[[x_{1},\cdots, x_{g^2}]]$ (resp.~$\Lambda[[x_{1},\dots, x_{g(h-g)}]]$).
\end{thm}

\begin{rmk}
Let~$\widehat{C}_{\Lambda}$ be the category of complete noetherian local~$\Lambda$-algebras~$R$ with maximal ideal~$\mathfrak{m}$ for which~$R/\mathfrak{m}^n\in C_{\Lambda}$ for all~$n\in \ZZ_{>0}$.
As in \cite[Section~2]{S68}, we can extend~$\textrm{Def}_{\Lambda}(A_{0},\alpha_{0})$ to a functor~$\widehat{\textrm{Def}}_{\Lambda}(A_{0},\alpha_{0})$ on~$\widehat{C}_{\Lambda}$ by the following\[ R \mapsto \projlim_{n}\textrm{Def}_{\Lambda}(A_{0},\alpha_{0})(R/\mathfrak{m}^n)=:\widehat{\textrm{Def}}_{\Lambda}(A_{0},\alpha_{0})(R),\]where~$R\in \widehat{C}_{\Lambda}$ with maximal ideal~$\mathfrak{m}$.
We denote the extension by~$\widehat{\textrm{Def}}_{\Lambda}(A_{0},\alpha_{0})$.
Note that the functor~$\textrm{Def}_{\Lambda}(A_{0},\alpha_{0})$ is pro-representable by a deformation ring~$\mathcal{R}$ by \thmref{thm:pro}.
So, we have isomorphisms~$\widehat{\textrm{Def}}_{\Lambda}(A_{0},\alpha_{0})(R)\cong \projlim_{n}\Hom_{\Lambda}(\mathcal{R},R/\mathfrak{m}^n)\cong \Hom_{\Lambda}(\mathcal{R},\projlim_{n}R/\mathfrak{m}^n)$.
Since~$R$ is complete, we see\[\widehat{\textrm{Def}}_{\Lambda}(A_{0},\alpha_{0})(R)\cong\Hom_{\Lambda}(\mathcal{R},R).\]
In particular, by definition of~$\widehat{\textrm{Def}}_{\Lambda}(A_{0},\alpha_{0})$, we see that any element in~$\widehat{\textrm{Def}}_{\Lambda}(A_{0},\alpha_{0})(R)$ is a \emph{formal} deformation of~$(A_{0},\alpha_{0})$ over~$R$ for any~$R\in \widehat{C}_{\Lambda}$.

\end{rmk}

\begin{rmk}\cite[pp.~46-47]{CCO14}\label{rmk:uniformal}
Let~$\mathcal{R}$ be the quotient, pro-representing~$\textrm{Def}_{\Lambda}(A_{0},\alpha_{0})$, as in \thmref{thm:pro}.
There is a universal \emph{formal} abelian scheme~$\mathfrak{A}$ over~$\textrm{Spf}(\mathcal{R})$ with structure morphism~$\rm{s}:\mathfrak{A}\to \textrm{Spf}(\mathcal{R})$.
That is,~$\mathfrak{A}$ represents~$\widehat{\textrm{Def}}_{\Lambda}(A_{0},\alpha_{0})$ on the category~$\widehat{C}_{\Lambda}$.
In particular, let~$R\in \widehat{C}_{\Lambda}$ be a complete noetherian local~$\Lambda$-algebra.
Given a local~$\Lambda$-algebra homomorphism~$\textrm{f}:\mathcal{R}\to R$, we have a Cartesian diagram \begin{equation*}
	\begin{tikzcd}
		\mathfrak{A}\widehat{\otimes} R \arrow[r]\arrow[d]\drar[phantom, "\square"] & \mathfrak{A}\arrow[d,"\rm{s}"] \\
		\textrm{Spf}(R)\arrow[r,"\textrm{Spf}(\textrm{f})"'] & \textrm{Spf}(\mathcal{R}),
	\end{tikzcd}
\end{equation*} where~$\textrm{Spf}(\textrm{f})$ is the induced morphism of formal schemes by~$\textrm{f}$,~$\mathfrak{A}\widehat{\otimes}R$ is the fiber product in the category of formal schemes, and~$\textrm{Spf}(R)$ and~$\textrm{Spf}(\mathcal{R})$ are formal affine schemes over the formal affine scheme~$\textrm{Spf}(\Lambda)$.
Furthermore,~$\mathfrak{A}\widehat{\otimes} R $ is a formal deformation of~$(A_{0},\alpha_{0})$ over~$R$, which corresponds to~$\textrm{f}\in\Hom_{\Lambda}(\mathcal{R},R)$.
Conversely, every formal deformation of~$(A_{0},\alpha_{0})$ over~$R$ is a pullback of~$\mathfrak{A}$ over~$\textrm{Spf}(\mathcal{R})$ to~$\textrm{Spf}(R)$. 

\end{rmk}

We recall the Serre-Tate deformation theory in the following.
\begin{thm}(Serre-Tate deformation theory) \cite[Theorem~1.4.5.3]{CCO14} or \cite[Theorem~1.2.1]{K81}\label{thm:SerreTate}
	Let~$R$ be a ring in which a prime~$p$ is nilpotent, and let~$I\subset R$ be a nilpotent ideal.
	Let~$R_{0}:=R/I$.
	Let~$\mathfrak{G}(R)$ be the category of abelian scheme over~$R$.
	Let~${\rm{Def}}(R,R_{0})$ be the category of triples~$(A_{0},G,\varepsilon)$, where~$A_{0}$ is an abelian scheme over~$R_{0}$,~$G$ is a~$p$-divisible group over~$R$, and~$\varepsilon: G_{0}\cong A_{0}[p^{\infty}]$ is an isomorphism of~$p$-divisible groups over~$R_{0}$.
	Then the functor\[\mathfrak{G}(R) \to {\rm{Def}}(R,R_{0}),\] given by \[A\mapsto (A\otimes (R/I) ,A[p^{\infty}]), {\rm{can}}:A[p^{\infty}] \otimes R/I \cong (A\otimes (R/I))[p^{\infty}]),\] is an equivalence of categories, where~${\rm{can}}:A[p^{\infty}] \otimes (R/I) \cong (A\otimes (R/I))[p^{\infty}])$ is a canonical isomorphism.
\end{thm}
\begin{rmk}\cite[Example~1.4.5.4]{CCO14}\label{rem:STformal}
Let~$R$ be a complete local noetherian ring with mixed characteristic. 
Let~$A_{0}$ be an abelian variety over the residue field of~$R$. 
Then, a deformation of the~$p$-divisible group~$A_{0}[p^{\infty}]$ of~$A_{0}$ to a~$p$-divisible group over~$R$ corresponds to a deformation of~$A_{0}$ to a \emph{formal} abelian scheme~$\mathfrak{A}$ over~$R$ by \thmref{thm:SerreTate}.
\end{rmk}

Note that~$\mathfrak{A}$ in \remref{rem:STformal} is not necessarily \emph{algebraizable}.
That is,~$\mathfrak{A}$ is not necessarily a formal completion of some abelian scheme. 
However, under certain situation, we may apply the following theorem, such that~$\mathfrak{A}$ is algebraizable, which is one of the steps in the proof of \thmref{lem:changeofcoe}. 
\begin{thm}\cite[Theorem~2.2.3]{CCO14}\label{thm:algebrai}
	Let~$R$ be a~$\ZZ_{p}$-flat~$1$-dimensional complete local noetherian ring with residue characteristic~$p$, and let~$\mathfrak{U}$ be a formal abelian scheme over~$R$ of relative dimension~$g>0$.
	Assume that the finite dimensional~$\QQ$-algebra~$\QQ\otimes_{\ZZ}\End_{R}(\mathfrak{U})$ contains a CM subalgebra~$L$ with~$[L:\QQ]= 2g$.
	Then~$\mathfrak{U}$ is an algebraizable if and only if the~$R[1/p]$-linear~$L$-action on~${\rm{Lie}}(\mathfrak{U})[1/p]$ is given by a CM type over each geometric point of~${\rm{Spec}}(R[1/p])$ in the sense of \defref{def:avCM}~$(3)$. 
\end{thm}
\begin{cor}\label{cor:onetoeach}
Keep notations as in \thmref{thm:algebrai}. 
	Assume that~${\rm{Spec}}(R[1/p])$ is connected. 
	If the~$R[1/p]$-linear~$L$-action on~${\rm{Lie}}(\mathfrak{U})[1/p]$ is given by a CM type over a geometric point, then the same holds for every geometric point of~${\rm{Spec}}(R[1/p])$.
\end{cor}
\begin{pf}
Let~$K$ be a geometric point of~$\textrm{Spec}(R[1/p])$,
such that the~$R[1/p]$-linear~$L$-action~$L\embed \End_{K}(\textrm{Lie}(\mathfrak{U})[1/p]\otimes K)$ is given by a CM type~$\Phi$ of~$L$.
	Then, we have decomposition~$\textrm{Lie}(\mathfrak{U})[1/p]\otimes K= \prod_{\varphi\in\Phi}\left(\textrm{Lie}(\mathfrak{U})[1/p]\otimes K\right)_{\varphi}$, where~$\left(\textrm{Lie}(\mathfrak{U})[1/p]\otimes K\right)_{\varphi}$ is the eigenspace of the~$L$-action with respect to~$\varphi$ for each~$\varphi\in\Phi$.
	Moreover,~$\dim_{K} \left(\textrm{Lie}(\mathfrak{U})[1/p]\otimes K\right)_{\varphi} = 1$ for each~$\varphi\in\Phi$ by \cite[Lemma~1.5.2]{CCO14}. 
	
We consider a function\[{{\rm{rk}}_{\varphi}:\rm{Spec}}(R[1/p]) \to \ZZ,\]
which is given by~$x\mapsto \dim_{\kappa(x)}\textrm{Lie}(\mathfrak{U})[1/p]_{\varphi}\otimes\kappa(x)$, where~$x\in\textrm{Spec}(R[1/p])$ corresponds to any prime ideal of~$R[1/p]$,~$\kappa(x)$ is the residue field of~$R[1/p]$ at~$x$, and~$\textrm{Lie}(\mathfrak{U})[1/p]_{\varphi}$ is the eigenspace of~$\varphi$ for each~$\varphi\in\Phi$.
The function~$\textrm{rk}_{\varphi}$ is locally constant for each~$\varphi\in\Phi$ because~$\textrm{Lie}(\mathfrak{U})[1/p]$ (and hence~$\textrm{Lie}(\mathfrak{U})[1/p]_{\varphi}$ for any~$\varphi\in\Phi$) is a locally free~$\mathcal{O}_{\textrm{Spec}(R[1/p])}$-sheaf of finite rank, where~$\mathcal{O}_{\textrm{Spec}(R[1/p])}$ is the structure sheaf of~$\textrm{Spec}(R[1/p])$.
Since~$\textrm{Spec}(R[1/p])$ is connected by assumption, the function~$\textrm{rk}_{\varphi}$ is constant for each~$\varphi\in\Phi$.
Then,\[\dim_{K'} (\textrm{Lie}(\mathfrak{U})[1/p]\otimes K')_{\varphi} = \dim_{K} (\textrm{Lie}(\mathfrak{U})[1/p]\otimes K)_{\varphi},\]which implies that the~$R[1/p]$-linear~$L$-action on~$\textrm{Lie}(\mathfrak{U})[1/p]\otimes K'$ is given by the CM type~$\Phi$ by \textit{loc. cit} for any geometric point~$K'$ of~$\textrm{Spec}(R[1/p])$.
This completes the proof.\qed
\end{pf}

\begin{rmk}
	There exists a counterexample that the~$L$-action on~$\mathfrak{U}$ is not given by a CM type of~$L$, see \cite[Section~4.1.2]{CCO14}.
\end{rmk}

Therefore, let~$R$ be a complete local noetherian ring with residue field~$\FF_{q}$, a finite field containing~$\FF_{p}$. 
Let~$A_{0}$ be an abelian variety over~$\FF_{q}$, admitting smCM by CM field~$L$.
Assume~$O_{L}\otimes\ZZ_{p}\embed\End_{\FF_{q}}(A_{0})\otimes\ZZ_{p}$, where~$O_{L}$ is the ring of integers of~$L$.  
A deformation of the~$p$-divisible group~$A_{0}[p^{\infty}]$ with~$L\otimes\QQ_{p}\embed\End_{\FF_{q}}^{0}(A_{0}[p^{\infty}])$ over~$R$ corresponds to a CM lifting of~$A_{0}$.

\section{Lie types.}\label{sec:Lie type}
In this section, we show a base-change result in \thmref{lem:changeofcoe} and show that the Lie types of simple superspecial abelian surfaces over~$\FF_{p}$ in the isogeny class are the same, see \propref{lem:sameLietype}. 
\subsection{Definitions.}
Let~$L/\QQ$ be a CM field with~$[L:\QQ]:=2g>0$, and let~$O_{L}$ be the ring of integers of~$L$.
Let~$A_{0}/\FF_{q}$ be a simple abelian variety of dimension~$g$ over a finite field~$\FF_{q}$ that admits smCM by~$L$ with CM structure~${\iota}:L\embed\End^{0}_{\FF_{q}}(A_{0})$. 
Furthermore,~${\iota}$ induces a CM structure of~$A_{0}\otimes_{\FF_{q}}\overline{\FF}_{p}$ by composing~${\iota}$ with the canonical inclusion~$\End_{\FF_{q}}(A_{0})\embed\End_{\overline{\FF}_{p}}(A_{0}\otimes_{\FF_{q}}\overline{\FF}_{p})$.
With slight ambiguity, we denote by~${\iota}$ the induced CM structure of~$A_{0}\otimes_{\FF_{q}}\overline{\FF}_{p}$. 

We recall a definition of the Lie type of~$(A_{0}\otimes\overline{\FF}_{p},{\iota})$, which is a datum obtained from the Dieudonn\'e module of~$A_{0}\otimes\overline{\FF}_{p}$.

We start with recalling~$O_{L}$-linear CM~$p$-divisible groups.
\begin{defn}\cite[Section~3.7.1.2]{CCO14} or \cite[Section~4.2.2]{CCO14}\label{def:CMp} 
Let~$L/\QQ$ be a CM field with ring of integers~$O_{L}$.
For any~$w\in\Sigma_{L,p}$, let~$L_{w}$ be the completion of~$L$ at~$w$ with valuation ring~$O_{L_{w}}$.
\begin{enumerate}
	\item An~\emph{$O_{L_{w}}$-linear (resp.~$L_{w}$-linear) CM~$p$-divisible group}~$(X,{\iota})$ is a pair consisting of a~$p$-divisible~$X$ group over a field~$K$ of characteristic~$p>0$ with height~$h>0$, and~${\iota}$ is a ring injection~$O_{L_{w}} \embed \End_{K}(X)$ (resp.~$\QQ_{p}$-algebra injection~$L_{w}\embed \End^{0}_{K}(X):=\End_{K}(X)\otimes_{\ZZ_{p}}\QQ_{p}$) such that~$[L_{w}:\QQ_{p}]=h$ for any~$w\in\Sigma_{L,p}$.
	We say that~$X$ \emph{has CM by~$L_{w}$}.
	\item An~\emph{$O_{L}$-linear (resp.~$L$-linear) CM~$p$-divisible group}~$(X,{\iota})$ is a pair consisting of a~$p$-divisible~$X$ group over a field~$K$ of characteristic~$p>0$ with height~$h>0$, and~${\iota}$ is a ring injection~$O_{L}\otimes_{\ZZ}\ZZ_{p}\embed \End_{K}(X)$ (resp.~$\QQ_{p}$-algebra injection~$L\otimes_{\QQ}\QQ_{p}\embed \End^{0}_{K}(X):=\End_{K}(X)\otimes_{\ZZ_{p}}\QQ_{p}$) such that~$[L\otimes\QQ_{p}:\QQ_{p}]=h$.
	We say that~$X$ \emph{has CM by~$L\otimes\QQ_{p}$}.
	\item A morphism~$\varphi$ of~$O_{L}$-linear CM~$p$-divisible groups~$(X_{1},{\iota}_{1})$,~$(X_{2},{\iota}_{2})$ is~$O_{L}$-linear (resp.~$L$-linear) if it is~$O_{L}\otimes\ZZ_{p}$-linear (resp.~$L\otimes\QQ_{p}$-linear), that is,~$\varphi\circ {\iota}_{1}(\ell)={\iota}_{2}(\ell)\circ \varphi$ for all~$\ell\in O_{L}\otimes\ZZ_{p}$ (resp.~$\ell\in L\otimes\QQ_{p}$).
\end{enumerate}
\end{defn}

Similar to \defref{def:CMtype}, we have following definitions. 
\begin{defn} \cite[Definition~3.7.1.4]{CCO14}\label{def:CMtypep}
	Fix an algebraic closure~$\overline{\QQ}_{p}$ of~$\QQ_{p}$.
	\begin{enumerate}
		\item For a finite-dimensional commutative semisimple~$\QQ_{p}$-algebra~$F$, a\emph{~$\overline{\QQ}_{p}$-valued~$p$-adic CM type} for~$F$ is a subset~$\Phi\subset \Hom_{\QQ_{p}}(F,\overline{\QQ}_{p})$.
		\item Let~$F$ be a finite-dimensional non-zero commutative semisimple~$\QQ_{p}$-algebra such that~$F=\prod_{j}F_{j}$, where~$F_{j}$ is a finite field extension over~$\QQ_{p}$ for any~$j$.
		The \emph{reflex field}~$E_{j}:=E(F_{j},\Phi_{j})\subset\overline{\QQ}_{p}$ over~$\QQ_{p}$ is the subfield corresponding to the open subgroup of elements~$\sigma\in\Gal(\overline{\QQ}_{p}/\QQ_{p})$ such that~$\sigma(\Phi_{j})=\Phi_{j}$ inside~$\Hom_{\QQ_{p}}(F_{j},\overline{\QQ}_{p})$ for every~$j$.
	The \emph{reflex field}~$E:=E(F,\Phi)\subset \overline{\QQ}_{p}$ is the composite field over~$\QQ_{p}$ of~$E_{j}$ in~$\overline{\QQ}_{p}$ for all~$j$.
	\end{enumerate}
\end{defn}
\begin{ex}\cite[p.~171]{CCO14}\label{ex:rrp}
Let~$A$ be an abelian scheme over a complete discrete valuation ring~$R$, whose residue field~$\kappa$ is of characteristic~$p>0$ and fraction field~$K$ is of characteristic~$0$.
Suppose that~$A$ admits smCM by a CM field~$L$.
Fix an algebraic closure~$\overline{K}$ of~$K$.
Let~$\overline{\QQ}_{p}$ be an algebraic closure of~$\QQ_{p}$ inside~$\overline{K}$.
Let~$\overline{\QQ}\subset \overline{\QQ}_{p}$ be an algebraic closure of~$\QQ$.

Let~$A[p^{\infty}]$ be the associated~$p$-divisible group of~$A$.
Let~$(L,\Phi)$ be an associated CM type of~$A$, where\begin{equation}\label{eq:padicCM}
	\Phi\subset \Hom_{\textrm{ring}}(L,\overline{K}) = \bigsqcup_{w|p}\Hom_{\QQ_{p}}(L_{w},\overline{\QQ}_{p}) = \Hom_{\QQ_{p}}(L\otimes_{\QQ}\QQ_{p},\overline{\QQ}_{p})
\end{equation} is a subset of~$\Hom_{\QQ_{p}}(L\otimes_{\QQ}\QQ_{p},\overline{\QQ}_{p})$, and the index runs through all places in~$L$ lying above~$p$.

Recall the notation~$\Phi_{w}$ from \defref{def:emb}~$(2)$.
Then\begin{equation}\label{eq:Phi}
	\Phi=\bigsqcup_{w|p}\Phi_{w},
\end{equation} where the index runs through places~$w$ in~$L$ lying above~$p$, with the chosen embedding~$\overline{\QQ}\embed\overline{\QQ}_{p}$. 
As we may identify~$\textrm{Lie}(A\otimes_{R}K)$ and~$\textrm{Lie}(A[p^{\infty}]\otimes_{R}K)$, we see that~$\Phi$ is a~$p$-adic CM type of~$A[p^{\infty}]$.
That is,~$\dim A[p^{\infty}]_{w} = |\Phi_{w}|$ for any place~$w$ in~$L$ lying above~$p$, where~$A[p^{\infty}]_{w}$ is an~${L_{w}}$-linear~$p$-divisible group and~$A[p^{\infty}]=\prod_{w|p}A[p^{\infty}]_{w}$ following the decomposition of~$L\otimes\QQ_{p}$.

Moreover, let~$E(L,\Phi)$ be the reflex field as in \defref{def:reflexfield}.
Let~$\nu$ be the place in~$E(L,\Phi)$, which is induced by the chosen embedding~$\overline{\QQ}\embed\overline{\QQ}_{p}$.
The completion of~$E(L,\Phi)$ at~$\nu$ is the reflex field~$E(L\otimes_{\QQ}\QQ_{p},\Phi)$, as in \defref{def:CMtypep}.
Indeed, the reflex field~$E(L,\Phi)=\QQ(\sum_{\sigma\in\Phi}\sigma x|x\in L)$ is a subfield of~$\overline{\QQ}$.
So, its completion~$E(L,\Phi)_{\nu}$ is the fixed field of~$\langle\sigma\in\Gal(\overline{\QQ}_{p}/\QQ_{p})\,|\,\sigma\Phi = \Phi\rangle$ in~$\overline{\QQ}_{p}$.
By \eqref{eq:Phi}, we have~$\langle\sigma\in\Gal(\overline{\QQ}_{p}/\QQ_{p})\,|\,\sigma\Phi = \Phi\rangle = \bigcap_{w|p} \langle\sigma\in\Gal(\overline{\QQ}_{p}/\QQ_{p})\,|\,\sigma\Phi_{w} = \Phi_{w}\rangle$, where the index runs through places~$w$ in~$L$ lying above~$p$, and~$\Phi_{w}:=\{\varphi\in\Phi\,|\,\varphi\,\,\textrm{induces}~w\}$ is a~$p$-adic CM type of~$L_{w}$ for any place~$w$ in~$L$ lying above~$p$.
This implies\[E(L,\Phi)_{\nu}=E(L\otimes_{\QQ}\QQ_{p},\Phi)\] is the reflex field in~$\overline{\QQ}_{p}$.
\end{ex}

In the following, we define a Lie type.
\begin{defn}\label{def:glt}\cite[Section~4]{Yu04}
	Let~$\overline{\FF}_{p}$ be an algebraically closed field of characteristic~$p>0$.
	Let~$L$ be a CM field with the maximal totally real subfield~$L_{0}$.
	Let~$Y$ be an~$O_{L}\otimes_{\ZZ}\ZZ_{p}$-linear~$p$-divisible group over~$\overline{\FF}_{p}$.
	That is, we fix an embedding\[\iota: O_{L}\otimes_{\ZZ}\ZZ_{p}\embed \End_{\overline{\FF}_{p}}(Y).\] 
	Let~$\textrm{V}_{Y}$ be the Verschiebung morphism of the (covariant) Dieudonn\'e module~$M_{Y}$ of~$Y$ over the Witt ring~$W:=W(\overline{\FF}_{p})$ of~$\overline{\FF}_{p}$.
	The morphism~$\textrm{V}_{Y}$ is obtained functorially by the relative Frobenius morphism of~$Y$. 
	
	Note that~$M_{Y}$ can be decomposed into~$\bigoplus_{w|p}M_{Y,w}$ as an~$O_{L}\otimes_{\ZZ}\ZZ_{p}$-module, where the index runs through all places~$w$ in~$L$ lying above~$p$, and~$M_{Y,w}$ is an~$O_{L_{w}}$-module for any place~$w|p$ in~$L$.
	Let~$v$ be a place in~$L_{0}$ lying above~$p$, with inertia degree~$f(v|p)$, and let~$\Hom(O_{L_{0,v}^{\textrm{ur}}}, W):=\{\sigma_{i}\}_{i\in\ZZ/f(v|p)\ZZ}$ be the inertia group of~$v$ over~$p$, which is isomorphic to~$\ZZ/f(v|p)\ZZ$.
	Then,\begin{equation}\label{eq:MYw}
		M_{Y,w}=\begin{cases}
			\bigoplus_{i\in\ZZ/f(v|p)\ZZ}(M_{Y,w}^{i} \oplus M_{Y,\overline{w}}^{i}) \qquad\,\,\,\,\,\textrm{if~$v=w\overline{w}$ splits in~$L$},\\
			\bigoplus_{i\in\ZZ/f(v|p)\ZZ} (M_{Y,w}^{i}\oplus M_{Y,w}^{i+f(v|p)}) \quad\textrm{if~$v=w$ is inert in~$L$},\\
			\bigoplus_{i\in\ZZ/f(v|p)\ZZ} M_{Y,w}^{i} \qquad\qquad\quad\quad\,\,\,\,\,\,\textrm{if~$v=w$ is ramified in~$L$},
		\end{cases}
	\end{equation} where~$M_{Y,w}^{i}\oplus M_{Y,\overline{w}}^{i}$ (resp.~$M_{Y,w}^{i}\oplus M_{Y,w}^{i+f(v|p)}$,~$M_{Y,w}^{i}$) is~$W^{i} \oplus W^{i}$-module (resp.~$W^{i}\oplus W^{\overline{i}}$,~$W^{i}$) of rank~$1$ if~$v$ splits (resp.~$v$ is inert, is ramified) in~$L$, where~$W^{i}:=W\otimes_{O_{L_{w}^{\textrm{ur}}},\sigma_{i}}O_{L_{w}^{\textrm{ur}}}$ and~$W^{\overline{i}}:=W\otimes_{O_{L_{w}^{\textrm{ur}}},\sigma_{i+f(v|p)}}O_{L_{w}^{\textrm{ur}}}$ (if~$v$ is inert, is ramified) for any~$i\in\ZZ/f(v|p)\ZZ$.
	
	There exists a~$W^{i}$-generator~$x_{w}^{i}$ of~$M_{Y,w}^{i}$ such that~$\textrm{V}_{Y}x_{w}^{i+1}=\iota(\pi_{w})^{e_{w}^{i}}x_{w}^{i}$, where~$\pi_{w}$ is a uniformizer of~$O_{L_{w}}$, for any~$i\in\ZZ/f(v|p)\ZZ$, see \cite[Section~4]{Yu04}.
	The following tuples are the\emph{~$w$-component of the Lie type of~$Y$(with respect to~$\iota$)}:
	\begin{align*}
		& (e_{w}^{i},e_{\overline{w}}^{i})_{i\in\ZZ/f(v|p)\ZZ} \qquad\quad\,\,\,\textrm{if~$v=w\overline{w}$ splits in~$L$},\\
		& (e_{w}^{i}, e_{w}^{i+f(v|p)})_{i\in\ZZ/f(v|p)\ZZ}  \quad\textrm{if~$v=w$ is inert in~$L$}, \\
		& (e_{w}^{i})_{i\in\ZZ/f(v|p)\ZZ} \qquad\quad\quad\,\,\,\,\,\,\textrm{if~$v=w$ is ramified in~$L$}.
	\end{align*}
A \emph{Lie type of~$Y$ (with respect to~$\iota$)} is the tuple consisting of the~$w$-component for all places~$w$ in~$L$ lying above~$p$.
\end{defn}


\begin{rmk}\label{rmk:lietypeuse}
	Let~$L$ be a CM field with ring of integers~$O_{L}$.
	Let~$Y$ be an~$O_{L}\otimes_{\ZZ}\ZZ_{p}$-linear CM~$p$-divisible group over a finite field~$\FF_{q}\supseteq\FF_{p}$.
	Let~$w$ be a place in~$L$ lying above~$p$, such that the inertia field~$\FF_{w}$ of~$L$ at~$w$ is a subfield of~$\FF_{q}$ (up to an isomorphism.)
	We can compute the~$w$-component of the Lie type of~$Y$, which we will explain in the following.
	In particular, the Lie type of~$Y$ is the same as the one of~$Y\otimes_{\FF_{q}}\overline{\FF}_{p}$ by \cite[Lemma~4.2.4]{CCO14}.
	
	Let~$M_{Y}$ be the (covariant) Dieudonn\'e module of~$Y$ over~$W(\FF_{q})$.
	From \eqref{eq:MYw}, we see that the construction of a~$w$-component of a Lie type depends on the group~$\Hom(O_{L_{w}^{\textrm{ur}}}, W(\overline{\FF}_{p}))$, where~$L_{w}^{\textrm{ur}}$ is the maximal unramified subextension of~$L_{w}$ over~$\QQ_{p}$, whose ring of integers is~$O_{L_{w}^{\textrm{ur}}}$.
	Since~$\FF_{w}$ is a subfield of~$\FF_{q}$, we have the equality\begin{equation}\label{eq:wq}
		\Hom(O_{L_{w}^{\textrm{ur}}}, W(\overline{\FF}_{p}))= \Hom(O_{L_{w}^{\textrm{ur}}}, W(\FF_{q})).
	\end{equation}
	Indeed,~$\Hom(O_{L_{w}^{\textrm{ur}}}, W(\overline{\FF}_{p}))$ is bijective to~$\Hom(\FF_{w},\overline{\FF}_{p})$.
	Fix an embedding~$\FF_{q}\embed\overline{\FF}_{p}$.
	Then, we may identify~$\Hom(\FF_{w},\FF_{q})$ and~$\Hom(\FF_{w},\overline{\FF}_{p})$ and hence obtain the Equation \eqref{eq:wq}.

Therefore, we may compute the Lie type of~$Y$ over~$\FF_{q}$ if~$\FF_{w}$ is a subfield of~$\FF_{q}$ for any place~$w$ in~$L$ lying above~$p$.
In the remainder of the paper, we will sometimes omit the field of definition of a~$p$-divisible group when we mention the Lie type of it due to this reason.
\end{rmk}

In the following example, we will see that one may use a CM type of a CM field to compute a Lie type under certain conditions.
\begin{ex}\cite[Section~(5.3)]{Yu04}\label{ex:glt}
Let~$L/\QQ$ be a CM field with maximal totally real subfield~$L_{0}$.
Let~$A_{0}/{\FF}_{q}$ be a simple abelian variety, admitting smCM by CM field~$L$.
Let~${\iota}$ be a CM structure of~$A_{0}$, such that it induces an embedding~$O_{L}\otimes\ZZ_{p}\embed\End_{\overline{\FF}_{p}}(A_{0}\otimes\overline{\FF}_{p})\otimes\ZZ_{p}$, which we denote by~${\iota}$ again.
Let~$A_{0}[p^{\infty}]$ be the~$p$-divisible group of~$A_{0}$ over~$\FF_{q}$.
Then,~$(X:=A_{0}[p^{\infty}]\otimes\overline{\FF}_{p},{\iota})$ is an~$O_{L}$-linear CM~$p$-divisible group over~$\overline{\FF}_{p}$. 

Fix an embedding~$\overline{\QQ}\embed\overline{\QQ}_{p}$.
Suppose that~$(A_{0}\otimes\overline{\FF}_{p},{\iota})$ has a CML.
Let~$w$ be any place in~$L$ lying above~$p$. 
Let~$\Phi$ be a~$\overline{\QQ}_{p}$-valued CM type of~$L$, such that~${\iota}$ is induced by~$\Phi$.
That is, there exists an abelian variety over~$\overline{\QQ}$ of CM type~$(L,\Phi)$, such that its CM structure~$\iota_{\Phi}$, which is induced by~$\Phi$, reduces to~$\tilde{\iota}$.
Such a CM type of~$L$ exists because we assume that~$(A_{0}\otimes\overline{\FF}_{p},{\iota})$ has a CML.

There are three cases of the place~$v$, which is below~$w$, in~$L_{0}$. 
Let~$f_{v}$ be the inertia degree of~$v$ over~$p$.
Define the subset~$\Phi_{w}^{i}$ (resp.~$\Phi_{w}^{i+f_{v}}$ if~$v$ is inert in~$L$) 
to be the intersection of the preimage of~$r^{-1}(i)$ (resp.~$r^{-1}(i+f_{v})$) and~$\Phi_{w}$, where 
\begin{equation*} r:
    \begin{cases}
        \Hom(L_{w},\overline{\QQ}_{p}) \to \Hom(L_{w}^{\textrm{ur}},\overline{\QQ}_{p})=\ZZ/f_{v}\ZZ,\quad\,\,\,\text{if $v=w\overline{w}$ splits completely in~$L$;}\\
        \Hom(L_{w},\overline{\QQ}_{p}) \to \Hom(L_{w}^{\textrm{ur}},\overline{\QQ}_{p})=\ZZ/f_{v}\ZZ, \quad\,\,\,\text{if~$v$ is ramified in~$L$;}\\
        \Hom(L_{w},\overline{\QQ}_{p}) \to \Hom(L_{w}^{\textrm{ur}},\overline{\QQ}_{p})=\ZZ/2f_{v}\ZZ,\quad\text{if~$v$ is inert in~$L$,}\\
    \end{cases}
\end{equation*}is the restriction map, and~$L_{w}^{\textrm{ur}}$ is the maximal unramified subextension of~$L_{w}$ for all $i=1,\dots,f_{v}$. 
We have the datum
\begin{equation}\label{eq:CMLie}
        \begin{cases}
        \left(|\Phi_{w}^{i}|,|\Phi_{\overline{w}}^{i}|\right)_{i\in \ZZ/f_{v}\ZZ}, \quad\quad\, \text{if $v=w\overline{w}$ splits completely in~$L$;}\\
        \left(|\Phi_{w}^{i}|\right)_{i\in \ZZ/f_{v}\ZZ},\qquad\qquad\,\,\, \text{if $v$ is ramified in~$L$;}\\
        \left(|\Phi_{w}^{i}|,|\Phi_{w}^{i+f_{v}}|\right)_{i\in \ZZ/f_{v}\ZZ},\quad \text{if $v$ is inert in~$L$,}
    \end{cases}\end{equation} which is the Lie type of~$X$ (with respect to~$\iota$).
    Since~$X$ is the~$p$-divisible group arising from the abelian variety~$A_{0}$, we also call \eqref{eq:CMLie} the \emph{Lie type of~$A_{0}$ (with respect to~$\iota$.)} 
In particular, \begin{equation*}
        \begin{cases}
       \dim \textrm{Lie}(X_{w})^{i}=|\Phi_{w}^{i}|,\,\,\dim \textrm{Lie}(X_{\overline{w}})^{i}=|\Phi_{\overline{w}}^{i}|, \quad\quad\, \text{if $v=w\overline{w}$ splits completely in~$L$;}\\
        \dim \textrm{Lie}(X_{w})^{i}=|\Phi_{w}^{i}|,\qquad\qquad\qquad\qquad\qquad\qquad \text{if $v$ is ramified in~$L$;}\\
        \dim \textrm{Lie}(X_{w})^{i}=|\Phi_{w}^{i}|,\,\, \dim \textrm{Lie}(X_{w})^{\overline{i}}=|\Phi_{w}^{\overline{i}}|,\,\qquad \text{if $v$ is inert in~$L$,}
    \end{cases}\end{equation*}where~$\overline{i}:=i+f_{v}$ for any~$i=1,\dots,f_{v}$.
In addition, let~$e_{v}$ be the ramification index of~$w$ at~$p$.
We have
\begin{equation}\label{eq:gltchar}
        \begin{cases}
        |\Phi_{w}^{i}|+|\Phi_{\overline{w}}^{i}| =e_{v}, \quad\quad\, \text{if $v=w\overline{w}$ splits completely in~$L$;}\\
        |\Phi_{w}^{i}|=e_{v},\,\,\,\,\,\,\,\,\,\,\,\quad\qquad\,\,\, \text{if $v$ is ramified in~$L$;}\\
        |\Phi_{w}^{i}|+|\Phi_{w}^{i+f_{v}}|=e_{v},\quad \text{if $v$ is inert in~$L$,}
    \end{cases}\end{equation}
    for any~$i=1,\dots, f_{v}$.
    Such a Lie type is called \emph{good Lie type}, see \cite[Section~(5.3)]{Yu04}. 
\end{ex}

\subsection{Base change}
We are now in the position to prove the following theorem,
which will be used in the proof of the main theorem in \secref{chap:ss}.
We remark that the last step of the proof is similar to \cite[Section~4.6.4.3]{CCO14}.

\begin{thm}\label{lem:changeofcoe}
Let~$A_{0}/\FF_{q}$ be an abelian variety, admitting smCM by a CM field~$L$.
Let~$O_{L}$ be the ring of integers of~$L$. 
Let~$R\supseteq \ZZ_{p}$ be a discrete valuation ring with residue field~$\overline{\FF}_{p}$.
Suppose that~$R$ is a finitely generated~$W(\overline{\FF}_{p})$-module.
Assume that~$A_{0}\otimes\overline{\FF}_{p}$ has a CM lifting~$A$ over~$R$. 
Then,~$A_{0}/\FF_{q}$ has a CM lifting.
\end{thm}
\begin{pf}
Let~$\alpha_{0}:L\embed \End_{\FF_{p}}^{0}(A_{0})$ be a CM structure of~$A_{0}$ over~$\FF_{q}$.
In particular, we have induced embeddings~$\alpha_{0}:O\embed\End_{\FF_{q}}(A_{0})$ integrally and~$\alpha_{0}[p^{\infty}]: O_{}\otimes\ZZ_{p}\embed\End_{\FF_{q}}(A_{0}[p^{\infty}])$ locally at~$p$, where~$O:=\alpha_{0}^{-1}(\End_{\FF_{q}}(A_{0}))\cap L$ is a CM order of~$L$, and~$A_{0}[p^{\infty}]$ is the~$p$-divisible group of~$A_{0}$ over~$\FF_{q}$.
With slight ambiguity of notations, we denote the embedding~$O\embed\End_{\overline{\FF}_{p}}(A_{0}\otimes\overline{\FF}_{p})$, which is obtained by the composition of~$\alpha_{0}$ and the canonical inclusion, by~$\alpha_{0}$ again.

In the following, we outline the proof. 
We first want to show that~$\widehat{\textrm{Def}}_{\Lambda'}(A_{0}\otimes\overline{\FF}_{p},\alpha_{0})(\mathcal{O})$ is non-empty, where~$\Lambda':=W(\overline{\FF}_{p})$ is a complete local noetherian ring with residue field~$\overline{\FF}_{p}$, and~$\mathcal{O}$ is a discrete valuation ring over~$\Lambda$ with residue field~$\overline{\FF}_{p}$.
And, we wil apply \thmref{thm:algebrai} to show that an element in~$\widehat{\textrm{Def}}_{\Lambda'}(A_{0},\alpha_{0})(\mathcal{O})$ is algebraizable.
Then, we show that~$\widehat{\textrm{Def}}_{\Lambda}(A_{0},\alpha_{0})(D)$ is non-empty by using the element, where~$\Lambda:=W(\FF_{q})$ is a complete local noetherian ring with residue field~$\FF_{q}$, and~$D$ is a domain over~$\Lambda$ with residue field~$\FF_{q}$. 
Together with \thmref{thm:algebrai} again, this will imply that~$A_{0}/\FF_{q}$ has a CM lifting by definition of deformations, see \remref{rmk:defandCM}.

Let~$O_{L}$ be the ring of integers of~$L$.
By assumption, we have an embedding~$O_{L}\otimes \ZZ_{p}\embed \End_{\overline{\FF}_{p}}(A_{0}[p^{\infty}]\otimes\overline{\FF}_{p})$ after composing the canonical embedding~$\End_{\FF_{q}}(A_{0}[p^{\infty}])\embed\End_{\overline{\FF}_{p}}(A_{0}[p^{\infty}]\otimes\overline{\FF}_{p})$.
Therefore, we can consider the Lie type of~$A_{0}\otimes\overline{\FF}_{p}$.
Moreover, since~$A_{0}\otimes\overline{\FF}_{p}$ has a CM lifting by assumption, the Lie type of~$A_{0}\otimes\overline{\FF}_{p}$ is good by \cite[Lemma~4.4.4.(i)]{CCO14}.
Let~$\Phi$ be a CM type of~$L$, such that the Lie type of~$A_{0}\otimes\overline{\FF}_{p}$ is obtained by~$\Phi$ as in \exref{ex:glt}.
More precisely, recall the notation~$\Phi_{w}$ from \defref{def:emb} for any place~$w$ in~$L$ lying above~$p$.
In particular,~$\Phi_{w}$ realizes the Lie type of~$A_{0}\otimes\overline{\FF}_{p}$ as in \eqref{eq:CMLie} when the index runs through every place~$w$ in~$L$ lying above~$p$.

Set~$K_{0}:= W(\overline{\FF}_{p})[1/p]$ with algebraic closure~$\overline{K}_{0}$.
It follows from \cite[Proposition~4.4.6(ii)]{CCO14} that for any finite subextension~$K_{w}$ in~$\overline{K}_{0}/K_{0}$, containing the reflex field of~$(L_{w},\Phi_{w})$, there exists an~$O_{L_{w}}$-linear~$p$-divisible group~$\mathcal{Y}_{w}$ over~$O_{K_{w}}$ such that\[\mathcal{Y}_{w}\otimes\overline{\FF}_{p}\cong A_{0}[p^{\infty}]_{w}\otimes\overline{\FF}_{p}\]
$O_{L_{w}}$-linearly, where~$O_{K_{w}}$ is the ring of integers of~$K_{w}$.
Therefore, let~$K$ be the composite field of~$K_{w}$, for any place~$w$ in~$L$ lying above~$p$, and let~$\mathcal{Y}:=\prod_{w}\mathcal{Y}_{w}$ be the~$p$-divisible group over~$O_{K}=:\mathcal{O}$.
Then,~$\mathcal{Y}\otimes\overline{\FF}_{p}$ is~\emph{$O_{L}$-linearly}~$\overline{\FF}_{p}$-isomorphic to~$A_{0}[p^{\infty}]\otimes\overline{\FF}_{p}$. 
In particular, this means that a deformation of the CM~$p$-divisible group~$A_{0}[p^{\infty}]\otimes\overline{\FF}_{p}$ exists.
By \thmref{thm:SerreTate} and \thmref{thm:algebrai}, we see that~$(A_{0}\otimes\overline{\FF}_{p}, \alpha_{0})$ has a deformation~$(\mathcal{A},\alpha)$ over~$O_{K}$, where~$\mathcal{A}$ is an abelian scheme over~$O_{K}$, and~$\alpha:O\embed\End_{O_{K}}(\mathcal{A})$ is an embedding of rings, whose induced~$L$-action on~$\textrm{Lie}(\mathcal{A})[1/p]$ is given by a CM type.

Lastly, we want to find a domain~$D$ over~$\Lambda$ with residue field~$\FF_{q}$, such that~$\widehat{\textrm{Def}}_{\Lambda}(A_{0},\alpha_{0})(D)$ is non-empty and has an algebraizable element.
Here, we regard~$\alpha_{0}$ as an embedding~$O\embed\End_{\FF_{q}}(A_{0})$ of rings.
We start with considering the deformation functor of~$(A_{0}[p^{\infty}],\alpha_{0}[p^{\infty}])$. 
The functor has a formal deformation ring~$\mathcal{R}$ by \thmref{thm:pro}.
And the universal formal deformation~$\mathfrak{X}$ of the functor over~$\mathcal{R}$ is also a universal deformation over~$\mathcal{R}$ since~$p$-divisible groups are built from torsions that are \emph{finite flat} over the base, see \cite[p.~49]{CCO14}. 
Let~$(\mathcal{A}[p^{\infty}],\alpha[p^{\infty}])$ be the~$p$-divisible group of~$\mathcal{A}$, where~$\alpha[p^{\infty}]$ is the embedding induced by~$\alpha$.
So,~$((A_{0}\otimes\overline{\FF}_{p})[p^{\infty}],\alpha_{0}[p^{\infty}])$ has a deformation~$(\mathcal{A}[p^\infty],\alpha[p^{\infty}])$ to~$O_{K}$ as we have seen before. 
There is the canonical local map (and hence non-zero) of rings\[{\rm{g}}:\mathcal{R}\to O_{K}.\]
We consider the domain~$D:={\rm{g}}(\mathcal{R})\subseteq O_{K}$ of characteristic~$0$.
Since~$D$ is isomorphic to a quotient~$\mathcal{R}/\ker \textrm{g}$ of~$\mathcal{R}$, its residue field is~$\FF_{q}$.

In the remainder of the proof, we first show that~$D$ is an order in a finite extension of~$\Lambda'[1/p]$.
So,~$A_{0}[p^{\infty}]$ admits an~$O_{L}$-linear lifting over~$D$, which descends~$(\mathcal{A}[p^{\infty}],\alpha[p^{\infty}])$, by pulling back the universal deformation~$\mathfrak{X}$ over~$\mathcal{R}$ along~$\rm{g}$. 
Then, it follows from \thmref{thm:SerreTate} and \thmref{thm:algebrai} that~$A_{0}$ over~$\FF_{q}$ has a CM lifting over~$D$.  

Note that we have the following commutative diagram:
\begin{equation}\label{eq:gRDOK}
	\begin{tikzcd}
		& \textrm{Spec}(D)\arrow[dl,"\iota"'] & \\
	\textrm{Spec}(\mathcal{R})	& & \textrm{Spec}(O_{K}), \arrow[ll,"\textrm{Spec}(\mathrm{g})"]\arrow[ul,"h"']
	\end{tikzcd}
\end{equation}where~$\iota$ is a closed immersion, and~$h$ is induced by the inclusion~$D \subseteq O_{K}$.
Since~$D$ is a domain, it has the generic point~$\eta$, such that its Zariski closure~$\overline{\{\iota(\eta)\}}^{\textrm{Zar}}$ equals~$\iota(\textrm{Spec}(D))$.
Moreover, the image~$\iota(\eta)$ in~$\textrm{Spec}(\mathcal{R})$ corresponds to the prime ideal~$\ker\textrm{g}$. 
In particular,~$p\notin \ker\textrm{g}$ because~$D$ has characteristic~$0$.
Let~$\mathfrak{q}\subset \mathcal{R}$ be a minimal prime ideal such that\begin{equation}\label{eq:inclusionqg}
	\mathfrak{q}\subseteq \ker\mathrm{g}
\end{equation}
and hence~$p\notin \mathfrak{q}$.
Let~$\eta_{\mathfrak{q}}$ be the generic point in~$\textrm{Spec}(\mathcal{R})$ corresponding to~$\mathfrak{q}\subset \mathcal{R}$.
(We say that the generic point~$\eta_{\mathfrak{q}}$ has \emph{characteristic~$0$}.)
Let~$D':= \mathcal{R}/\mathfrak{q}$ be the quotient ring of~$\mathcal{R}$, such that~$\textrm{Spec}(D')=\overline{\{\eta_{\mathfrak{q}}\}}^{\textrm{Zar}}$ is the Zariski closure of~$\{\eta_{\mathfrak{q}}\}$ in~$\textrm{Spec}(\mathcal{R})$.
Then, we see \begin{equation}\label{eq:generic}
	\textrm{Spec}(D') \supset \{\iota(\eta)\}
\end{equation} from \eqref{eq:inclusionqg}. 
Taking the Zariski closure on both sides of \eqref{eq:generic}, we have~$\textrm{Spec}(D') \supset \iota(\textrm{Spec}(D))$ in~$\textrm{Spec}(\mathcal{R})$.
In particular,~$\textrm{Spec}(D)$ is a closed subscheme of~$\textrm{Spec}(D')$ via~$\iota$.
Furthermore, the diagram \eqref{eq:gRDOK} factors through
\begin{equation*}
	\begin{tikzcd}
		& \textrm{Spec}(D) \arrow[dl,"\iota"']&\\
		\textrm{Spec}(D')& & \textrm{Spec}(O_{K}).\arrow[ll,"\textrm{Spec}(\mathrm{g})"]\arrow[ul,"h"']
	\end{tikzcd}
\end{equation*}
By \cite[Lemma~2.1.1.2]{CCO14}, we see that~$\textrm{Spec}(D')$ is~$\Lambda$-finite.
This implies that~$\textrm{Spec}(D)$ is~$\Lambda$-finite becuase a closed immersion is always finite by \cite[Proposition~12.11]{GW10}.
Also,~$D$ is~$\Lambda$-flat because~$D$ is torsion-free and~$\Lambda=W(\FF_{q})$ is a principal ideal domain by \cite[Chapter XVI, Proposition~3.2]{Lang12}. 
To conclude,~$D$ is a~$\Lambda$-finite flat domain with local inclusion into~$O_{K}$.
So,~$D$ is an order in a finite extension of~$\Lambda[1/p]$.
This completes the proof.\qed



\end{pf}
\begin{rmk}\label{rmk:Yu}
Let~$L$ be a CM field with~$[L:\QQ]=2g>0$ with maximal totally real subfield~$L_{0}$.
Let~$A_{0}$ be a simple abelian variety of dimension~$g$ over a finite field~$\FF_{q}$, admitting smCM by~$L$.
In Yu's work, he assumed that the ring of integers~$O_{L}$ is embedded into~$\End_{\FF_{q}}(A_{0})$, that is,\begin{equation}\label{eq:OL}
	O_{L}\embed\End_{{\FF}_{q}}(A_{0}).
\end{equation}
Then, a good Lie type of~$A_{0}\otimes\overline{\FF}_{p}$ can detect the CM liftability of~$A_{0}\otimes\overline{\FF}_{p}$, see \cite[Proposition~(5.5)]{Yu04}, which states that there is a bijection between two sets of abelian varieties, admitting embeddings from~$O_{L}$ to their endomorphism rings, the set of abelian varieties over~$\overline{\QQ}$ of CM type~$(L,\Phi)$, and the set of abelian varieties over~$\overline{\FF}_{p}$ satisfying the dimension condition, with good Lie type realized by~$(L,\Phi)$.  
The dimension condition means that the dimensions of~$\text{Lie}(A_{0})_{w}$ and of~$A_{0}[p^{\infty}]_{w}$ are the same when the place~$v\in\Sigma_{L_{0},p}$ splits~$v=w\overline{w}$ in~$L$. 

However, we do not have the embedding \eqref{eq:OL} in some of our cases as we will see in \propref{rem:p=3} later.
\end{rmk}

\subsection{Simple superspecial abelian surfaces over~$\FF_{p}$.}\label{ex:ss ab surf}
Let~$p>0$ be a prime.
An abelian variety over a finite field~$\FF_{q}\supseteq \FF_{p}$ is said to be \emph{superspecial} if it is isomorphic to a power of a supersingular elliptic curve over~$\overline{\FF}_{p}$. 
A \emph{Weil-$p$ number} is an algebraic integer~$\pi$ whose absolute value in~$\mathbb{C}$ is a square root of~$p$ for any embedding to~$\mathbb{C}$.
Note that~$\pi$ depends on~$p$.

We denote by~$W_{p}^{\textrm{ss}}(2)$ the set of Weil-$p$-numbers,
in which every element corresponds to a simple supersingular abelian surface over~$\FF_{p}$ (up to~$\FF_{p}$-isogeny) by the Honda-Tate theorem, following \cite[Section~5.2]{XYY16}.
From now on, any simple abelian variety over~$\FF_{p}$ in the isogeny class determined by a Weil-$p$ number is obtained by the Honda-Tate theorem, unless otherwise specified.
The following is the list of elements in~$W_{p}^{\textrm{ss}}(2)$ as in \cite[Section~5.2]{XYY16},
\begin{equation}\label{Weil-p}
   \quad\quad\qquad\,\, W_{2}^{\textrm{ss}}(2) = \{\sqrt{2}, \sqrt{2} \zeta_{3}, \sqrt{2} \zeta_{12}, \pm \sqrt{2} \zeta_{24}\};\end{equation}
\begin{equation} \label{Weil-p1}  W_{3}^{\textrm{ss}}(2) = \{\sqrt{3}, \sqrt{3} \zeta_{3}, \sqrt{3} \zeta_{8}\};\end{equation}
\begin{equation}\label{Weil-p2}
  \qquad\qquad\quad\quad\,\,\,  W_{5}^{\textrm{ss}}(2) = \{\sqrt{5}, \sqrt{5} \zeta_{3}, \sqrt{5} \zeta_{8}, \sqrt{5}\zeta_{12},\pm \sqrt{5}\zeta_{5}\};
	\end{equation}
\begin{equation}
\label{Weil-p3}   \qquad\qquad\quad\qquad\,\, W_{p}^{\textrm{ss}}(2) = \{\sqrt{p}, \sqrt{p} \zeta_{3}, \sqrt{p} \zeta_{8},\sqrt{p}\zeta_{12}\}\,\,\text{for $p\geq7$}.\end{equation}

\begin{defn}
	A Weil-$p$ number~$\pi$ is \textit{of our concern} if~$\pi\in\{\sqrt{p}\zeta_{3},\sqrt{p}\zeta_{8},\sqrt{p}\zeta_{12}, \pm\sqrt{5}\zeta_{5}\}$ when~$p\geq 5$, or~$\pi\in\{\sqrt{3}\zeta_{3},\sqrt{3}\zeta_{8}\}$ when~$p=3$, or~$\pi\in \{\sqrt{2}\zeta_{3},\sqrt{2}\zeta_{12},\pm\sqrt{2}\zeta_{24}\}$ when~$p=2$.
	In particular, a Weil-$p$ number of our concern is a non-real element in~$W_{p}^{\textrm{ss}}(2)$ for~$p>0$.
\end{defn}

We recall the following lemma, which appears in the proof of \cite[Proposition~5.1]{XYY16}.
\begin{lem} \label{lem:Dp}
Let~$p>0$ be a prime.
	Let~$\pi$ be a Weil-$p$ number of our concern, and let~$A_{\pi}/\FF_{p}$ be a simple supersingular abelian surface in the isogeny class determined by~$\pi$, which is obtained by the Honda-Tate theorem.
	Let~$M$ be the (covariant) Dieudonn\'e module of~$A_{\pi}$ over~$\ZZ_{p}$.
	Then,~$A_{\pi}/\FF_{p}$ is superspecial if and only if\[{\rm{V}^2}M=pM,\] where~${\rm{V}}$ is the operator of the Dieudonn\'e module~$M$, obtained from the Frobenius endomorphism~$\textrm{Frob}_{A_{\pi},p}$ of~$A_{\pi}$ functorially.
\end{lem}
\begin{pf}
	The sufficiency follows the definition of superspecial abelian varieties.
	We explain the necessity.
	Let~$\pi$ be a Weil-$p$ number of our concern.
	Let~$L:=\QQ(\pi)$ be a CM field. 
	Suppose that~$A_{\pi}/\FF_{p}$ is superspecial.
	Then, we have~$\textrm{V}^2(M\otimes_{\ZZ_{p}} W(\overline{\FF}_{p}))=p(M\otimes_{\ZZ_{p}} W(\overline{\FF}_{p}))$, and hence\[\textrm{V}^2/p\in \End_{R_{\overline{\FF}_{p}}}(M\otimes W(\overline{\FF}_{p})).\]
	Note that~$\textrm{V}$ is obtained from the Frobenius endomorphism~$\textrm{Frob}_{A_{\pi},p}$ of~$A_{\pi}/\FF_{p}$ functorially, and that the quotient~$\textrm{Frob}_{A_{\pi},p}^2/p$ is a root of unity in~$L$. 
Since~$A_{\pi}/\FF_{p}$ is simple, the endomorphism algebra~$\End_{\FF_{p}}^{0}(A_{\pi})$ equals~$L$.
In particular,~$\textrm{Frob}_{A_{\pi},p}^2/p\in \End_{\FF_{p}}^{0}(A_{\pi})$.
We have an isomorphism~$\End_{\FF_{p}}^{0}(A_{\pi})\otimes\QQ_{p}\cong\End_{R_{\FF_{p}}}(M)\otimes_{\ZZ_{p}}\QQ_{p}$ functorially by~$\mathbb{D}$, see \cite[p.~525]{W69}.
So,\[\textrm{V}^2/p\in \End_{R_{{\FF}_{p}}}(M)\otimes_{\ZZ_{p}}\QQ_{p}.\]
Therefore, we have~$\textrm{V}^2/p\in\End_{R_{\FF_{p}}}(M)$.	
This completes the proof.
\qed
\end{pf}
For a simple superspecial abelian surface over~$\FF_{p}$, which is in the isogeny class determined by a Weil-$p$ number of our concern, we have the following. 
\begin{cor}\label{def:ssp}
	Keep the notation as in \lemref{lem:Dp}.
	Suppose that~$A_{\pi}/\FF_{p}$ is superspecial.
	Let\[R_{{\rm{sp}}}:= \ZZ[{\rm{Frob}}_{A_{\pi},p},p{\rm{Frob}}_{A_{\pi},p}^{-1}, {\rm{Frob}}_{A_{\pi},p}^2/p]\] be an order in~$L$.
	Then~$R_{{\rm{sp}}}\otimes\ZZ_{p}\subseteq \End_{\FF_{p}}(A_{\pi})\otimes\ZZ_{p}$.
\end{cor}
\begin{pf}
Keep the notation in the proof of \lemref{lem:Dp}.
Then,~$\textrm{Frob}_{A_{\pi},p}^2/p\in\End_{\FF_{p}}(A_{\pi})\otimes\ZZ_{p}$ because~$\textrm{V}^2/p\in\End_{R_{\FF_{p}}}(M)$ from \textit{loc. cit}.
In addition,~$\End_{\FF_{p}}(A_{\pi})$ contains~$\textrm{Frob}_{A_{\pi},p}$ and~$p\textrm{Frob}_{A_{\pi},p}^{-1}$ by \cite[Proposition~3.5]{W69}.
This completes the proof.
\qed
\end{pf}
	



In fact, we will see that~$\End_{\FF_{p}}(A_{\pi})$ is isomorphic to~$O_{L}$ locally at~$p$, for almost all Weil-$p$ numbers~$\pi$ of our concern, in the following lemma.
\begin{lem}\label{lem:maxp}
	Let~$p>0$ be any prime, and let~$\pi$ be a Weil-$p$ number of our concern except~$\sqrt{3}\zeta_{3}$.
	Let~$L:=\QQ(\pi)$ be a CM field.
Let~$A_{\pi}/\FF_{p}$ be a simple superspecial abelian surface, in the isogeny class determined by~$\pi$, obtained by the Honda-Tate theorem. 
	The endomorphism ring~$\End_{\FF_{p}}(A_{\pi})$ is equal to~$O_{L}$ locally at~$p$,
	that is, \begin{equation}\label{eq:locmax}
    O_{L}\otimes\ZZ_{p}=\End_{\FF_{p}}(A_{\pi})\otimes\ZZ_{p}.
\end{equation} 
\end{lem}
\begin{pf}
	We recall the result in \textit{loc. cit.} 
Fix a Weil-$p$ number~$\pi$ of our concern.
Let~$L:=\QQ(\pi)$ be a field with ring of integers~$O_{L}$, and let~$A_{\pi}/\FF_{p}$ be a simple superspecial abelian surface in the isogeny class determined by~$\pi$.
Let~$\textrm{Frob}_{A_{\pi},p}$ be the Frobenius endomorphism of~$A_{\pi}/\FF_{p}$.
The order\[R_{\textrm{sp}}:=\ZZ[\textrm{Frob}_{A_{\pi},p}, p\textrm{Frob}_{A_{\pi},p}^{-1},\textrm{Frob}_{A_{\pi},p}^2/p]\] of~$L$ is a subring of~$\End_{\FF_{p}}(A_{\pi})$ locally at~$p$ because of \corref{def:ssp}.
Note that~$\Gal(L/\QQ)$ acts on~$R_{\textrm{sp}}$.
Hence,~$R_{\textrm{sp}}=\ZZ[\pi,p\pi^{-1},\pi^{2}/p]$. 
By \textit{loc. cit.}, when~$\pi\in\{\pm\sqrt{5}\zeta_{5}, \sqrt{2}\zeta_{24}, \sqrt{p}\zeta_{8} \,\,(p\neq 2) \}$, we have~$R_{\textrm{sp}}=O_{L}$. 
Therefore, the index~$[O_{L}:R_{\textrm{sp}}]=1$ in these cases.
For the remaining cases, we have the following tables, in which~$B$ is an endomorphism ring of some abelian variety over~$\FF_{p}$ by \cite[Theorem~6.1]{W69}, containing~$R_{\textrm{sp}}$. 
\begin{table}[H]
    \centering
    \begin{tabular}{|c |c|c|}\hline
       $\pi=\sqrt{p}\zeta_{3}$ & 
         $[O_{L}:R_{\textrm{sp}}]$ & $R_{\textrm{sp}} \subset B \subset O_{L}$ \\\hline
        $p=2$& $1$ & $O_{L}$ \\\hline
        $p=3$& $3$ & $R_{\textrm{sp}},O_{L}$ \\\hline
        $p \equiv 3 \,\,(\textrm{mod}\,\, 4), p \neq 3$& $1$ & $O_{L}$ \\\hline
        $p\equiv 1 \,\,(\textrm{mod}\,\, 4)$ &$4$ & $R_{\textrm{sp}},O_{L}$ \\\hline
    \end{tabular}
    \caption{Indices of~$R_{\textrm{sp}}$ in~$O_{L}$ when~$\pi=\sqrt{p}\zeta_{3}$.}\label{table:zeta3}
\end{table}

\begin{table}[H]
    \centering
    \begin{tabular}{|c |c|c|}\hline
       $\pi=\sqrt{p}\zeta_{12} \,\,(p\neq 3)$ & 
         $[O_{L}:R_{\textrm{sp}}]$ & $R_{\textrm{sp}} \subset B \subset O_{L}$ \\\hline
        $p=2$& $1$ & $O_{L}$ \\\hline
        $p \equiv 3 \,\,(\textrm{mod}\,\, 4)$ & $4$ & $R_{\textrm{sp}}, O_{L}$ \\\hline
        $p\equiv 1 \,\,(\textrm{mod}\,\, 4)$ &$1$ & $O_{L}$ \\\hline
    \end{tabular}
    \caption{Indices of~$R_{\textrm{sp}}$ in~$O_{L}$ when~$\pi=\sqrt{p}\zeta_{12}$ and~$p\neq 3$.} \label{table:zeta12}
\end{table}
This completes the proof.\qed
\end{pf}

\begin{rmk}
    Let~$C/\FF_{7}$ be a principal abelian surface in the isogeny class determined by the Weil-$7$ number~$\pi=\sqrt{7}\zeta_{3}$. 
    Although~$\End_{\FF_{7}}(C)=O_{L}$ is the maximal order in~$\End_{\FF_{7}}^{0}(C)=L$ from \autoref{table:zeta3}, it is not necessarily maximal in~$\End^{0}_{\FF_{7^i}}(C\otimes\FF_{7^i})$ for $i=2,\dots,6$ since the endomorphism algebra may change after the base change. 
    We observe that the endomorphism algebras of~$C\otimes\FF_{7^i}$ are different when~$i=1,3,6$.
    The endomorphism algebra of~$C\otimes\FF_{7^3}$ is the quaternion algebra over $\QQ(\sqrt{7})$ ramified exactly at two infinite places since~${\rm{Frob}}_{C,7}^3$ is real, \cite[Case 2, p.~528]{W69}. 
    Since~$\pi^6\in\QQ$ and hence~${\rm{Frob}}_{C,7}^{6}\in\QQ$, we see that~$\End_{\FF_{7^6}}^{0}(C\otimes~\FF_{7^6}) =~\End_{\overline{\FF}_{7}}^{0}(C\otimes\overline{\FF}_{7})$ and hence is a two-by-two matrix algebra over a quaternion algebra over~$\QQ$ ramified at~$7$ and the infinite place.
    In particular, there exist strict inclusions $\End_{\FF_{7}}(C)\subset \End_{\FF_{7^3}}(C\otimes \FF_{7^3}) \subset \End_{\FF_{7^6}}(C\otimes\FF_{7^6})$. Here, we consider the cases when~$\pm{\rm{Frob}}_{C,7}$ equals one of~$\pi$ or~$7\pi^{-1}$ since~$\pi^{3}$ is real implies~$(7\pi^{-1})^3$ is real and vice versa. 
\end{rmk}

\begin{cor}\label{cor:maxp}
Let~$\pi$ be a Weil-$p$ number of our concern except~$\sqrt{3}\zeta_{3}$.
Let~$A_{\pi}$ be a simple superspecial abelian surface over~$\FF_{p}$ in the isogeny class determined by~$\pi$.
Then, any CM structure of~$A_{\pi}/\FF_{p}$ induces an isomorphism~$O_{L}\otimes\ZZ_{p}\cong\End_{\FF_{p}}(A_{\pi})\otimes\ZZ_{p}$.	
\end{cor}
\begin{pf}
	Since the endomorphism algebra~$\End_{\FF_{p}}^{0}(A_{\pi})$ is the CM field~$L=\QQ(\pi)$, there is the unique (up to isomorphism) maximal order in~$\End_{\FF_{p}}^{0}(A_{\pi})\otimes_{}\QQ_{p}=L\otimes\QQ_{p}$.
Let~${\iota}$ be a CM structure of~$A_{\pi}/\FF_{p}$.
Note that~${\iota}:L\cong\End_{\FF_{p}}^{0}(A_{\pi})$ is an isomorphism, which induces an isomorphism\[{\iota}:L\otimes\QQ_{p}\cong\End_{\FF_{p}}^{0}(A_{\pi})\otimes\QQ_{p},\] which we denote by~${\iota}$ again with slight ambiguity.
Then, we see that~$u{\iota}(O_{L}\otimes\ZZ_{p})u^{-1}=O_{L}\otimes\ZZ_{p}$ for some~$u\in(O_{L}\otimes\ZZ_{p})^{\times}$ by \cite[Theorem~(17.3)]{R75}.
Hence,~${\iota}(O_{L}\otimes\ZZ_{p}) = O_{L}\otimes\ZZ_{p}$.
So, by \lemref{lem:maxp}, we have an isomorphism~${\iota}:O_{L}\otimes\ZZ_{p}\cong\End_{\FF_{p}}(A_{\pi})\otimes\ZZ_{p}$, which is the restriction of~${\iota}$ to~$O_{L}\otimes\ZZ_{p}$.
This completes the proof.\qed
\end{pf}

The following proposition, related to Lie types, is used in \secref{chap:ss}. 
\begin{prop}\label{lem:sameLietype}
	Any simple superspecial abelian surface over~$\FF_{p}$ in the isogeny class determined by~$\pi$, where~$\pi$ is a Weil-$p$ number of our concern, has the same Lie type (with respect to some CM structures of~$A_{\pi}\otimes\overline{\FF}_{p}$ factoring through~$\End_{\FF_{p}}^{0}(A_{\pi})$ (resp.~$\End_{\FF_{3^6}}^{0}(A_{\pi}\otimes\FF_{3^6})$) when~$\pi \neq \sqrt{3}\zeta_{3}$ (resp.~$\pi=\sqrt{3}\zeta_{3}$)).
\end{prop}
\begin{pf}
Let~$\pi$ be a Weil-$p$ number of our concern.
Let~$A_{\pi}$ over~$\FF_{p}$ be a simple superspecial abelian surface in the isogeny class determined by~$\pi$.
Let~$L:=\QQ(\pi)$ be a CM field. 

We separate the discussion into two cases.
\begin{enumerate}
	\item When~$\pi\neq \sqrt{3}\zeta_{3}$, we have~$O_{L}\otimes\ZZ_{p}\cong \End_{\FF_{p}}(A_{\pi})\otimes\ZZ_{p}$ by \corref{cor:maxp}.
	By a CM structure factoring through~$\End_{\FF_{p}}^{0}(A_{\pi})$, we mean  a CM structure~$\iota: L\embed\End_{\overline{\FF}_{p}}^{0}(A_{\pi}\otimes\overline{\FF}_{p})$, 
such that the following diagram commutes
\begin{equation*}
	\begin{tikzcd}
		L \arrow[hook,rr,"\iota"]\arrow[rd,hook']&& \End_{\overline{\FF}_{p}}^{0}(A_{\pi}\otimes\overline{\FF}_{p})\\
		&\End_{\FF_{p}}^{0}(A_{\pi}),\arrow[ru,hook'] &
	\end{tikzcd}
\end{equation*}where the inclusion on the right-hand side is canonical.
With slight ambiguity, we denote the embedding on the left-hand side by~$\iota$ again.
In particular,~$\iota$ induces an embedding\[O_{L}\otimes\ZZ_{p}\embed \End_{\overline{\FF}_{p}}(A_{\pi}\otimes\overline{\FF}_{p})\otimes\ZZ_{p}. \]
	Any isomorphism class in the isogeny class determined by~$\pi$ has the same Frobenius endomorphism up to units in~$O_{L}\otimes\ZZ_{p}$.
	Note that the semisimple algebra~$\End_{\FF_{p}}^{0}(A_{\pi})$ determines the isogeny class.
	So,~$\iota$ induces a CM structure on every isomorphism class over~$\FF_{p}$ in the isogeny class.  
	Then, every isomorphism class over~$\FF_{p}$ in the isogeny class has the same Lie type (with respect to~$\iota$), since units are automorphisms~$\left(\End_{\overline{\FF}_{p}}(A_{\pi}\otimes\overline{\FF}_{p})\otimes\ZZ_{p}\right)^{\times}$ of the (covariant) Dieudonn\'e module of~$A_{\pi}\otimes\overline{\FF}_{p}$. 
	\item When~$\pi=\sqrt{3}\zeta_{3}$ (and hence~$p=3$), 
	there are two cases from \autoref{table:zeta3}.
	If~$\End_{\FF_{p}}(A_{\pi})\cong O_{L}$, then the argument is the same as in the one in~$(1)$. 
	Suppose~$\End_{\FF_{p}}(A_{\pi}) \cong R_{\textrm{sp}}$.
	Let~$B_{\pi}$ be a simple superspecial abelian surface over~$\FF_{p}$ in the same isogeny class as~$A_{\pi}$ such that~$\End_{\FF_{p}}(B_{\pi})\cong O_{L}$. 

	Let~$q=p^6$. 
	By \cite[Theorem~2(d5)]{T66}, we know that~$\End_{\FF_{q}}^{0}(A_{\pi}\otimes\FF_{q})\cong \End_{\FF_{q}}^{0}(E\times E)$ because the center of the endomorphism algebra of~$A_{\pi}$ becomes rational after base change to~$\FF_{q}$, where~$E/\FF_{q}$ is a supersingular elliptic curve such that~$\End_{\overline{\FF}_{p}}(E\otimes\overline{\FF}_{p})=\End_{\FF_{q}}(E)$.
	In particular,~$\End_{\FF_{q}}(A_{\pi}\otimes\FF_{q}) = \End_{\overline{\FF}_{p}}(A_{\pi}\otimes\overline{\FF}_{p})$.
	Since~$A_{\pi}$ and~$B_{\pi}$ are superspecial, we have~$A_{\pi}\otimes\overline{\FF}_{p}\cong B_{\pi}\otimes\overline{\FF}_{p}$. 
	This implies that~$O_{L}$ is embedded into~$\End_{\FF_{q}}(A_{\pi}\otimes\FF_{q})$, which provides a CM structure~$\iota$ of~$A_{\pi}\otimes\overline{\FF}_{p}$.
	Note that~$\iota$ factors through~$\End_{\FF_{q}}(A_{\pi}\otimes\FF_{q})$.	
	This holds true for any isomorphism class with endomorphism ring, which is isomorphic to~$R_{\textrm{sp}}$ over~$\FF_{p}$, in the isogeny class determined by~$\pi$ since~$A_{\pi}$ is arbitrarily chosen.
	
	Because of the choice of a CM structure, there exists an element~$\ell\in L$ such that~$\iota(\ell)=\textrm{Frob}_{A_{\pi},q}$ by \cite[Theorem~5, Corollary~1]{ST68}.
	Let~$M$ be the (covariant) Dieudonn\'e module of~$A_{\pi}\otimes\FF_{q}$.
	We can decompose~$M=\bigoplus_{w} M_{w}$, where the index runs through every place~$w$ in~$L$ lying above~$p$, because it is an~$O_{L}\otimes\ZZ_{p}$-module (of rank~$1$).  
	Let~$\textrm{V}$ be the Verschiebung of~$M$.
	
	Note that~$3$ has inertia degree~$2$ in~$L$. 
	We can decompose~$M_{w}=M_{w}^1 \oplus M_{w}^2$ into eigenspaces of the~$O_{L_{w}^{\textrm{ur}}}$-action, where~$L_{w}^{\textrm{ur}}$ is the maximal unramified subextension of~$L_{w}$.
	Since the arguments for~$M_{w}^{1}$ and~$M_{w}^{2}$ are the same, we consider~$M_{w}^{1}$ in the following and omit~$M_{w}^2$. 
	Let~$m_{w}$ be a generator of~$M_{w}^{1}\otimes W(\overline{\FF}_{p})$.
	Then,~$\textrm{V}m_{w} = \iota(\pi_{w})^{e_{w}^1}m_{w}$, where~$\pi_{w}$ is a uniformizer of~$O_{L_{w}}$, and~$e_{w}^1$ is the~$w$-component of the Lie type of~$A_{\pi}\otimes\overline{\FF}_{p}$. 
	In particular,\[\textrm{V}^6 m_{w} = \iota(\pi_{w})^{6e_{w}^1}m_{w},\] which is functorially obtained by~$\iota(\ell)$.
	Note that Frobenius endomorphisms in the same isogeny class as~$A_{\pi}\otimes\FF_{q}$ are~$\pm\textrm{Frob}_{A_{\pi},q}$ because~$\textrm{Frob}_{A_{\pi},q}$ is rational. 
	Since~$\pm 1$ is a unit in~$O_{L}$, the remainder of the proof is similar to the one in~$(1)$. 
	
\end{enumerate}
The proof is hence complete. \qed
\end{pf}

\begin{rmk}
When~$\pi\neq \sqrt{3}\zeta_{3}$, the Lie type of~$A_{\pi}$ can be computed over~$\FF_{p}$ (resp.~$\FF_{p^2}$) if~$L/L_{0}$ splits (resp. is inert) because~$O_{L}\otimes\ZZ_{p}$ is embedded into the endomorphism ring~$\End_{\FF_{p}}(A_{\pi})\otimes\ZZ_{p}$, see \remref{rmk:lietypeuse}.
\end{rmk}

We compute the Lie type of~$A_{\pi}\otimes\overline{\FF}_{7}$ explicitly, where~$\pi=\sqrt{7}\zeta_{3}$ in the following example, to end the subsection.
\begin{ex}\label{ex:A7}
	Let~$\pi=\sqrt{7}\zeta_{3}$, and let~$A_{\pi}/\FF_{7}$ be a simple superspecial abelian surface over~$\FF_{7}$ in the isogeny class determined by~$\pi$.
	Let~${\rm{Frob}}_{A_{\pi},7}$ be the Frobenius endomorphism of~$A_{\pi}/\FF_{7}$.
	Then, we have~$\pm{\rm{Frob}}_{A_{\pi},7}=\sqrt{7}\zeta_{3}$ or~$\sqrt{7}\zeta_{3}^{2}$.
	Let~$L:=\QQ(\pi)$.
	Then,~$L=\QQ({\rm{Frob}}_{A_{\pi},7})=\End_{\FF_{7}}^{0}(A_{\pi})$ since~$A_{\pi}/\FF_{7}$ is simple.

	Let~$L_{0}:=\QQ(\sqrt{7})$ be the maximal totally real subfield of~$L$, and let~$v$ be the place in~$L_{0}$ lying over~$7$.
	Let~$w$ be a place in~$L$ lying over~$v$.
	We know that~$7$ is ramified in~$L_{0}=\QQ(\sqrt{7})$ and~$v=w\overline{w}$ splits completely in~$L$ since the Legendre symbol~$\left(\frac{-3}{7}\right)$ equals~$1$.
	Therefore, the completions~$L_{w}$ and~$L_{0,v}=\QQ_{7}(\sqrt{7})$ are identical, and the place~$v$ corresponds to the maximal ideal~$\sqrt{7}\ZZ_{7}[\sqrt{7}]$.
	Furthermore, we see the valuation of~${\rm{Frob}}_{A_{\pi},7}$ at~$w$ is~$1$, that is,\begin{equation}\label{eq:good}
		w({\rm{Frob}}_{A_{\pi},7})=1,
	\end{equation}since~$\zeta_{3}$ is a unit in~$\ZZ_{7}[\sqrt{7}]\subseteq \End_{\FF_{7}}(A_{\pi})\otimes\ZZ_{7}$.
	
	
	Let~$M$ be the (contravariant) Dieudonn\'e module of~$A_{\pi}\otimes\overline{\FF}_{7}$ over the Witt ring~$W(\overline{\FF}_{7})$.
	Note that~$O_{L}\otimes\ZZ_{7}\cong\End_{\FF_{7}}(A_{\pi})\otimes\ZZ_{7}$ by \corref{cor:maxp}.
	As an~$O_{L}\otimes W(\overline{\FF}_{7})$-module, we denote the~$w$-component of~$M$ by~$M_{w}$ for any~$w\in\Sigma_{L,7}$.
	Moreover, as an~$O_{L_w}\otimes_{\ZZ_{7}}\overline{\FF}_{7}$-module, we have\[\overline{M}_{w}/\text{F}\overline{M}_{w} \cong \overline{\FF}_{7}[\pi_{w}]/\pi_{w}^{e_{w}^{1}},\] where~$\overline{M}_{w}$ is~$M_{w}$ mod~$7$, and~$\pi_{w}$ is a uniformizer in~$O_{L_{w}}$.
	Similarly, we have \[\overline{M}_{\overline{w}}/\text{F}\overline{M}_{\overline{w}} \cong \overline{\FF}_{7}[\pi_{\overline{w}}]/\pi_{\overline{w}}^{e_{\overline{w}}^{1}},\] where~$\overline{M}_{\overline{w}}$ is~$M_{\overline{w}}$ mod~$7$, and~$\pi_{\overline{w}}$ is a uniformizer in~$O_{L_{\overline{w}}}$.
	Then,~$(e_{w}^{1},e_{\overline{w}}^{1})$ is the Lie type of~$A_{\pi}\otimes\overline{\FF}_{7}$.
	By Equation \eqref{eq:good}, we see~$e_{w}^{1}=1$.
	Similarly,~$e_{\overline{w}}^{1}=1$.
	In particular, we have\[e_{w}^{1}+e_{\overline{w}}^{1}=2,\] which is the condition of being good Lie type, see \exref{ex:glt}.
	Hence,~$A_{\pi}/\FF_{7}$ has good Lie type. \end{ex}
\section{Example: superspecial abelian surfaces with \texorpdfstring{$\pi= \sqrt{7}\zeta_{3}$}{}.}\label{ex:1}

Let~$\pi=\sqrt{7}\zeta_{3}$, which is a Weil-$7$ number of our concern.
Let~$L:=\QQ(\pi)$ be a CM field with maximal totally real subfield~$L_{0}:=\QQ(\sqrt{7})$.
In \thmref{thm:ApisCML}, we show that any simple superspecial abelian surface over~$\FF_{7}$ in the isogeny class determined by~$\pi$ has a CML.
We will generalize the result to any Weil-$p$ number of our concern and for any prime~$p>0$ in Section~\ref{chap:ss}. 

\subsection{Residual reflex condition }\label{sec:CML1}
Let~$A_{\pi}/\FF_{7}$ be a simple superspecial abelian surface, in the isogeny class determined by~$\pi$. 
Let~${\iota}:L\embed \End_{{\FF}_{7}}^{0}(A_{\pi})$ be a CM structure of~$A_{\pi}/\FF_{7}$.
In this subsection, we will see that~$A_{\pi}/\FF_{7}$ satisfies RRC in \propref{prop:RRC} for some (but not all) CM types of~$L$. 


We consider all CM types of~$L=\QQ(\pi)$ and check if~$A_{\pi}/\FF_{7}$ satisfies RRC. 
Fix an embedding~$i: \overline{\QQ}\hookrightarrow \overline{\QQ}_{7}$ such that \[\sqrt{7}\mapsto\sqrt{7}, \sqrt{-3}\mapsto -\sqrt{-3}.\] 
We consider all elements in~$\Hom(L,\overline{\QQ})$, which are \begin{align*}
    \varphi_{1}: \begin{cases}
        \sqrt{7} \mapsto \sqrt{7}\\
        \sqrt{-3} \mapsto -\sqrt{-3},
    \end{cases} & \varphi_{2}: \begin{cases}
        \sqrt{7} \mapsto -\sqrt{7}\\
        \sqrt{-3}\mapsto \sqrt{-3},
    \end{cases} \end{align*}\begin{align*}
    \varphi_{3}: \begin{cases}
        \sqrt{7} \mapsto -\sqrt{7}\\
        \sqrt{-3}\mapsto -\sqrt{-3}, \end{cases}& \varphi_{4}: \begin{cases}
        \sqrt{7} \mapsto \sqrt{7}\\
        \sqrt{-3}\mapsto \sqrt{-3}.
    \end{cases}
\end{align*} 

Let~$\Phi_{1}:=\{\varphi_{1},\varphi_{2}\}$,~$\Phi_{2}:=\{\varphi_{2},\varphi_{4}\}$,~$\Phi_{3}:=\{\varphi_{3},\varphi_{4}\}$, and~$\Phi_{4}:=\{\varphi_{1},\varphi_{3}\}$ be all CM types of~$L$.
We will use these notations in the remainder of the section, unless otherwise specified.
We find the reflex field associated with~$\Phi_{i}$ in \lemref{lem:reflex}, which is listed in \cite[\href{https://www.lmfdb.org/NumberField/4.0.7056.3}{Number Field~4.0.7056.3}]{lmfdb-reflex} for~$i=1,\dots,4$ as well. 
Recall the definition of a reflex field from \defref{def:reflexfield}. 
\begin{lem}\label{lem:reflex}
    The reflex field associated with~$\Phi_{1}$ and~$\Phi_{3}$ is~$\QQ(\sqrt{-21})$, and the reflex field associated with~$\Phi_{2}$ and~$\Phi_{4}$ is~$\QQ(\sqrt{-3})$.
\end{lem}
\begin{pf}
Note that~$\zeta_{3}=\frac{-1+\sqrt{-3}}{2}$ and that~$\zeta_{3}^{2}=\frac{-1-\sqrt{-3}}{2}$.
Hence, we see~$\zeta_{3}^{\varphi_{1}}=\zeta_{3}^{\varphi_{3}}=\zeta_{3}^2$ and~$\zeta_{3}^{\varphi_{2}}=\zeta_{3}^{\varphi_{4}}=\zeta_{3}$. 
As~$L=\QQ(\pi)$ is of degree~$4$ over~$\QQ$, every element~$x\in L$ is of the form~$a+b\pi+c\pi^2+d\pi^3$, where~$a,b,c,d$ are element in~$\QQ$. 
When the CM type is~$\Phi_{1}$, the summation~$x^{\varphi_{1}}+x^{\varphi_{2}}$ equals~$-7+a\sqrt{-21}$.
Since~$x\in L$ is any element, we see that the reflex field associated with~$\Phi_{1}$ is equal to~$\QQ(\sqrt{-21})$. We use the same method to compute the reflex field related to other CM types~$\Phi_{i}$ for $i=2,3,4$. This completes the proof. \qed
\end{pf}
\begin{lem}
Let~$\pi=\sqrt{7}\zeta_{3}$, and let~$L:=\QQ(\pi)$ be a CM field with maximal totally real subfield~$L_{0}$.
Let~$v$ be the place in~$L_{0}$ lying above~$7$.
Let~$w$ be the place above~$v$ in~$L$, induced by~$\varphi_{1}$ with respect to the embedding~$i:\overline{\QQ}\embed\overline{\QQ}_{7}$.
Let~$\overline{w}$ be another place in~$L$ lying above~$v$.
Then,~$\varphi_{1},\varphi_{3}$ induce~$w$, while~$\varphi_{2},\varphi_{4}$ induce~$\overline{w}$. 
\end{lem}
\begin{pf}
Let~$c$ be the non-trivial element in~$\Aut(L/\QQ)$ of order~$2$ as in \defref{def:CMtype}. 
Then, the element~$c$ maps~$\sqrt{7}$ to~$-\sqrt{7}$ and~$\sqrt{-3}\mapsto\sqrt{-3}$. 
Note that~$\varphi_{3} = \overline{\varphi}_{2}$ and that~$\varphi_{4} = \overline{\varphi}_{1} = \text{id}$, where~$\overline{\varphi}_{2}=\varphi_{2} \circ c$ (resp.~$\overline{\varphi}_{1}= \varphi_{1}\circ c$) is a composition of~$\varphi_{2}$ (resp.~$\varphi_{1}$) and~$c$. 
Now, we let~$w$ be the place in~$\Sigma_{L,7}$ induced by~$\varphi_{1}$.
That is, we have the following commutative diagram \[\begin{tikzcd}
    L \arrow[rr, hook, "i \circ \varphi_{1}"]\arrow[rd,"\textrm{incl}"'] & &\overline{\QQ}_{7}\\
    &L_{w}, \arrow[ru,"\textrm{id}"']&
\end{tikzcd}\] where the lower map~$\textrm{incl}$ is given by the canonical inclusion, and the map~$\textrm{id}$ is given by the identity map in~$\Hom_{\QQ_{7}}(L_{w},\overline{\QQ}_{7})$. 

Let~$v$ be the place in~$L_{0}$ lying below~$w$, which is ramified over~$7$, and~$w$ splits over $v$, say~$v=w\overline{w}$ in~$L$. 
Hence, we have~$L_{w}=L_{\overline{w}}=cL_{w}$,~$L_{w} = L_{0,v} = \QQ_{7}(\sqrt{7})$, and~$\sqrt{-3} \in \QQ_{7}$ by Hensel's lemma since~$T^2+3$ mod $7$ has roots~$\pm 2$ mod~$7$. 
Let~$\Hom_{\QQ_{7}}(L_{w},\overline{\QQ}_{7})=:\{\textrm{id},\rho\}$, where~$\rho$ is non-trivial of order~$2$.
That is,~$\rho(\sqrt{7})=-\sqrt{7}$, and~$\rho$ is the identity map on~$\QQ_{7}$.
Then, we have the following commutative diagram \[\begin{tikzcd}
    L \arrow[rr, hook, "i \circ \varphi_{2}\circ c =i\circ\varphi_{3}"]\arrow[rd,"\textrm{incl}"'] & &\overline{\QQ}_{7}\\
    &L_{w}. \arrow[ru,"\rho"']&
\end{tikzcd}\]
Note that~$\Hom_{c\QQ_{7}}(L_{\overline{w}},\overline{\QQ}_{7})=\{\psi\circ c\,|\,\psi\in \Hom_{\QQ_{7}}(L_{w},\overline{\QQ}_{7})\}$, and that~$i\circ \varphi_{2}$ induces an element in~$\Hom_{c\QQ_{7}}(L_{\overline{w}},\overline{\QQ}_{7})$ because it sends~$\sqrt{-3}$ to~$-\sqrt{-3}$.
With slight ambiguity, we denote~$\textrm{id},\rho$ again the two elements in~$\Hom_{c\QQ_{7}}(L_{\overline{w}},\overline{\QQ}_{7})$, where~$\rho$ is non-trivial of order~$2$. 
That is, we have the following commutative diagram \[\begin{tikzcd}
    L \arrow[rr, hook, "i \circ \varphi_{2}"]\arrow[rd,"\textrm{incl}"'] & &\overline{\QQ}_{7}\\
    &L_{\overline{w}}. \arrow[ru,"\rho"']&
\end{tikzcd}\] 
Similarly, we have the following commutative diagram \[\begin{tikzcd}
    L \arrow[rr, hook, "i \circ \varphi_{4}"]\arrow[rd,"\textrm{incl}"'] & &\overline{\QQ}_{7}\\
    &L_{\overline{w}}. \arrow[ru,"\textrm{id}"']&
\end{tikzcd}\]
So,~$\varphi_{1},\varphi_{3}$ induce~$w$, while~$\varphi_{2},\varphi_{4}$ induce~$\overline{w}$.
This completes the proof.\qed
\end{pf}



We summarize the important information about~$L$ in \autoref{ex1}.
\begin{table}[H]
    \centering
    \begin{tabular}{|c |c| c|c|c|c|}\hline
        \rm{CM type} & 
         $w$ & $\overline{w}$ & \rm{slope} & $(|\Phi_{\bullet, w}|,|\Phi_{\bullet,\overline{w}}|)$ & \rm{reflex field}\\\hline
        $\Phi_{1} =\left\{ \varphi_{1},\varphi_{2}\right\}$& $\varphi_{1}$ & $\varphi_{2}$& $(\frac{1}{2},\frac{1}{2})$& $(1,1)$ &$\QQ(\sqrt{-21})$\\\hline
        $\Phi_{2} =\left\{ \varphi_{2},\varphi_{4}\right\}$& $\emptyset$ & $\varphi_{2},\varphi_4{}$& $(0,1)$&$(0,2)$&$\QQ(\sqrt{-3})$\\\hline
        $\Phi_{3} =\left\{ \varphi_{3},\varphi_{4}\right\}$& $\varphi_{3}$ & $\varphi_{4}$& $(\frac{1}{2},\frac{1}{2})$&$(1,1)$&$\QQ(\sqrt{-21})$\\\hline
        $\Phi_{4} =\left\{\varphi_{1},\varphi_{3} \right\}$ &$\varphi_{1},\varphi_{3}$ & $\emptyset$& $(1,0)$&$(2,0)$&$\QQ(\sqrt{-3})$\\\hline
    \end{tabular}
    \caption{All CM types of~$L$ and their corresponding induced places, slopes,~$(|\Phi_{\bullet, w}|,|\Phi_{\bullet,\overline{w}}|)$, and reflex fields.} 
    \label{ex1}
\end{table}

We exclude certain~$\overline{\QQ}_{7}$-valued CM types of~$L$, which are not realized by~$A_{\pi}/\FF_{7}$, see \defref{def:ST}, and show that~$A_{\pi}/\FF_{7}$ satisfies RRC with the remaining CM types of~$L$ in \propref{all type} and \propref{prop:RRC}.
\begin{prop}\label{all type}
	Let~$\pi:=\sqrt{7}\zeta_{3}$, and let~$A_{\pi}$ be a simple superspecial abelian surface over~$\FF_{7}$ in the isogeny class determined by~$\pi$.
	Then, the abelian surface~$A_{\pi}/\FF_{7}$ does not satisfy RRC with~$\overline{\QQ}_{7}$-valued CM types~$(L,\Phi_{2}),(L,\Phi_{4})$. 
\end{prop}
\begin{pf}
Note that~$A_{\pi}/\FF_{7}$ is supersingular. 
Suppose that the slope of~$A_{\pi}$ can be obtained by~$\Phi_{2}$ or~$\Phi_{4}$ via the Shimura-Taniyama formula \eqref{ST}. 
Then, \[\frac{|\Phi_{i,w}|}{[L_{w}:\QQ_{7}]}\in\{0,1\},\] for~$i=2,4$ by \autoref{ex1}, which implies that~$A_{\pi}$ is ordinary, which is a contradiction. 
This completes the proof.
\qed
\end{pf}

We have the following proposition.
\begin{prop}\label{prop:RRC}
	Any superspecial abelian surface over~$\FF_{7}$ in the isogeny class determined by~$\pi:=\sqrt{7}\zeta_{3}$ satisfies RRC with~$\overline{\QQ}_{7}$-valued CM types~$\Phi_{1},\Phi_{3}$ of~$L$. 
\end{prop}
\begin{pf}
Let~$A_{\pi}/\FF_{7}$ be a simple superspecial abelian surface in the isogeny class determined by~$\pi$.
By \propref{all type}, we consider~$\Phi_{1},\Phi_{3}$.
The Shimura-Taniyama formula \eqref{ST} of the CM type~$(L,\Phi_{i})$ tells \[\frac{|\Phi_{i,w}|}{[L_{w}:\QQ_{7}]}=\frac{|\Phi_{i,\overline{w}}|}{[L_{\overline{w}}:\QQ_{7}]}=\frac{1}{2},\] from \autoref{ex1} for~$i=1,3$, which coincides with the slope of~$A_{\pi}$ since~$A_{\pi}$ is superspecial.

On the other hand, by \lemref{lem:reflex}, we know that the reflex field of~$\Phi_{1}$ and~$\Phi_{3}$ is~$\QQ(\sqrt{-21})$ (regard~$\Phi_{1},\Phi_{3}$ as CM types of~$L$), whose residue field at~$7$ is~$\FF_{7}$ because~$7$ is ramified in~$L_{0}$.
Since there is exactly one place in~$L_{0}$ lying above~$7$, the reflex field condition of RRC is satisfied. 
Since~$A_{\pi}/\FF_{7}$ is chosen arbitrarily, the argument holds for any abelian surface in this isogeny class.
Hence, the proof is complete.\qed
\end{pf}

\subsection{Proofs for the CM liftability}\label{sec:p7}
Let~$A_{\pi}/\FF_{7}$ be a simple superspecial abelian surface over~$\FF_{7}$ in the isogeny class determined by~$\pi=\sqrt{7}\zeta_{3}$. 
In this subsection, 
we prove that any simple superspecial abelian surface over~$\FF_{7}$, in the isogeny class, has a CM lifting (CML). 
We give two proofs.
One of them is an application of many results, while another seems to apply to more general cases, see \remref{rem:moregeneral}.  

We first recall the following theorem from Tate. 
\begin{thm}\cite[Main theorem]{T66}\label{thm:T66}
	Let~$\FF_{q}$ be a finite field. 
	Let~$A,B$ be abelian varieties over~$\FF_{q}$.
	Let~$M_{A},M_{B}$ be corresponding (contravariant) Dieudonn\'e modules of~$A$ and~$B$.
	Then,\[\Hom_{\FF_{q}}(A,B) \otimes\ZZ_{p}\cong \Hom_{W(\FF_{q})[{\rm{F,V}]}}(M_{B},M_{A}),\] where~$W(\FF_{q})$ is a Witt ring of~$\FF_{q}$-vectors. 
\end{thm}
\begin{thm}\label{thm:ApisCML}
	Any simple superspecial abelian surface over~$\FF_{7}$, in the isogeny class determined by~$\pi=\sqrt{7}\zeta_{3}$, has a CM lifting (CML). 
\end{thm}
\begin{proof}[First proof.]
Let~$A_{\pi}/\FF_{7}$ be a simple superspecial abelian surface in the isogeny class determined by~$\pi=\sqrt{7}\zeta_{3}$.
Let~$\textrm{Frob}_{A_{\pi},7}$ be the Frobenius endomorphism of~$A_{\pi}$ over~$\FF_{7}$.
Let~$L:=\QQ(\pi)$ be a CM field with maximal totally real subfield~$L_{0}$.
Then,~$L=\QQ(\pi)=\QQ(\textrm{Frob}_{A_{\pi},7})$ is the endomorphism algebra of~$A_{\pi}$ over~$\FF_{7}$ as~$A_{\pi}/\FF_{7}$ is assumed to be simple. 
Let~${\iota}_{}$ be a CM structure of~$A_{\pi}/\FF_{7}$.
Let~$v$ be the place in~$L_{0}$ lying above~$7$.
Let~$w,\overline{w}$ be the places in~$L$ lying above~$v$.

By \propref{prop:RRC}, we know that~$A_{\pi}/\FF_{7}$ satisfies RRC with respect to a fixed embedding~$\overline{\QQ}\embed\overline{\QQ}_{p}$.
Then, \cite[Theorem~2.5.3]{CCO14} says the following: There exists an abelian surface~$(A,\iota_{A,\Phi})$ with CM structure~$\iota_{A,\Phi}$ over a number field~$K$, whose~$\overline{\QQ}_{7}$-valued CM type is~$\Phi$, such that~$A$ has good reduction~$A_{0}/\FF_{7}$ at the place in~$\Sigma_{K,7}$, which is related to the embedding~$\overline{\QQ}\embed\overline{\QQ}_{7}$ in RRC.
In particular,~$A$ has good reduction~$A_{0}/\FF_{7}$ at the place in~$\Sigma_{K,7}$, which is related to the embedding~$\overline{\QQ}\embed\overline{\QQ}_{7}$ in RRC. 
And,~$A_{0}$ is~$L$-linearly $\FF_{7}$-isogenous to~$A_{\pi}$ with respect to the CM structure of~$A_{0}$, which is reduced from the one of its CM lifting, and the CM structure~$\tilde{\iota}_{}$ of~$A_{\pi}$.
Let~$\varphi:A_{0}\to A_{\pi}$ be the~$L$-linear~$\FF_{7}$-isogeny.
Let~$A_{0}[7^{\infty}]$ be the~$7$-divisible group over~$\FF_{7}$ of~$A_{0}$. 
Moreover,~$A_{0}[7^{\infty}]$ is~$O_{L}$-linear, see \cite[p.~133]{CCO14}.
This implies\[O_{L}\otimes\ZZ_{7}\cong \End_{\FF_{7}}(A_{0})\otimes\ZZ_{7}\] because~$\End_{\FF_{7}}(A_{0})= O_{L}$, see \autoref{table:zeta3}.
The~$7$-divisible groups~$A_{0}[7^{\infty}]$ and~$A_{\pi}[7^{\infty}]$ are~$O_{L}$-linear CM~$7$-divisible groups by \corref{cor:maxp}.
Furthermore, we see that~$A_{0}[7^{\infty}]$ and~$A_{\pi}[7^{\infty}]$ are~$O_{L}$-linearly~$\FF_{7}$-isogenous by the isogeny~$\varphi[7^{\infty}]$ induced by~$\varphi$. 

Let~$M_{0}, M_{\pi}$ be the (contravariant) Dieudonn\'e modules over~$\ZZ_{7}$ of~$A_{0},A_{\pi}$ respectively.
By \propref{lem:sameLietype}, we see that~$A_{0}$ and~$A_{\pi}$ have the same Lie type.
Therefore, by \cite[Proposition~4.2.9]{CCO14}, we know that~$M_{0}$ and~$M_{\pi}$ are~$L\otimes\QQ_{p}$-linearly isomorphic over~$\ZZ_{7}$.
So,~$\Hom_{\FF_{7}}(A_{\pi},A_{0})\otimes\ZZ_{7}\cong \Hom_{\ZZ_{7}[\textrm{F,V}]}(M_{0},M_{\pi})$, by \thmref{thm:T66}, contains an isomorphism~$\rho$.
With the same argument as in the proof of \cite[Theorem~7]{MW71}, there exists a sequence of elements of~$\Hom_{\FF_{7}}(A_{\pi},A_{0})$ approximating~$\rho$. 
In particular, we obtain an element~$\psi\in \Hom_{\FF_{7}}(A_{\pi},A_{0})$, which induces an isomorphism~$\psi': M_{0}\to M_{\pi}$. 
Note that~$\psi$ is an isogeny because if~$\ker \psi$ contained any positive dimension abelian variety, then~$\psi$ will not induce an isomorphism between~$M_{0}$ and~$M_{\pi}$. 
This implies that~$\psi$ is separable because the~$p$-primary part of~$\ker \psi$ corresponds to the cokernel of~$\psi'$, see \cite[p.~525]{W69}. 
It follows from \cite[Proposition~2.22]{BKM} that~$A_{\pi}$ admits a CM lifting.
The proof is complete. \qed 
\end{proof}

\begin{rmk}\label{rmk:RRC7}
We may alternatively apply \cite[Theorem~5.1]{W69} to obtain an~$\FF_{7}$-isogeny between~$A_{0}$ and~$A_{\pi}$. 
Note that~$\End_{\FF_{7}}(A_{\pi})$ and~$\End_{\FF_{7}}(A_{0})$ are equal to~$O_{L}$ from \autoref{table:zeta3}.
This implies that there exists a separable~$\FF_{7}$-isogeny between~$A_{\pi}$ and~$A_{0}$ by \cite[Theorem~5.1]{W69}.
So,~$A_{\pi}/\FF_{7}$ has a CML after applying \cite[Proposition~2.22]{BKM}, because~$A_{0}$ does. 

However, we remark that the endomorphism ring of a simple superspecial abelian surface over~$\FF_{p}$ does not equal~$O_{L}$ for every Weil-$p$ number of our concern.
The proof given in this remark therefore does not generalize, while the above proof does.  
\end{rmk}

Keep notations as in the proof of \thmref{thm:ApisCML}.
We will use Grothendieck-Messing theory, see \cite[Chapter V, Theorem~1.6]{Messing72}, \cite[Main theorem]{N81} and \cite[Proposition~2.3]{OS20}, to give the second proof of \thmref{thm:ApisCML} in the following.
\begin{rmk}\label{rem:moregeneral}
Note that we use the Lie type in the first proof of \thmref{thm:ApisCML}. 
	We expect the method in the second proof to be applied to more general setups of CM liftability than the Lie type because we may not have an embedding from~$O_{L}$ into the endomorphism ring locally at~$p$ even after base change to~$\overline{\FF}_{p}$ when we go beyond the case of superspecial abelian surfaces (for example, supersingular non-superspecial abelian surfaces). 
\end{rmk}

We recall Grothendieck-Messing theory in the following.
\begin{thm}\label{thm:GM}
	Let~$S$ denote the ring of integers, which is a slightly ramified extension of~$\ZZ_{p}$.
That is, the ramification index of the field of fraction~$\textrm{frac}(S)$ of~$S$ over~$\QQ_{7}$ is at most~$p$, see \cite[Definition~2.2]{OS20}.
Let~$\varpi$ be a uniformizer of~$S$.
The deformations of~$A_{\pi}[p^{\infty}]_{w}$ to~$S$ are in bijection with filtrations~${\rm{Fil}}_{w} \subset M_{\pi,w}\otimes_{\ZZ_{7}}S$, such that
\begin{enumerate}
	\item The quotient~$\left(M_{\pi,w}\otimes_{\ZZ_{7}}S\right)/{\rm{Fil}}_{w}$ is torsion-free, 
	\item The filtration~${\rm{Fil}}_{w}$ mod~$\varpi$$={\rm{Fil}}_{w}\otimes_{S}(S/\varpi S)$ equals the kernel of the Frobenius~$\textrm{F}$ on~$M_{\pi,w}\otimes_{\ZZ_{7}} S$ mod~$\varpi$ $=M_{\pi,w}\otimes_{\ZZ_{7}} S\otimes_{S}(S/\varpi S)$, and
	\item The filtration~${\rm{Fil}}_{w}$ is an~$O_{L_{w}}$-module. 
\end{enumerate}
\end{thm}
\begin{rmk}
	Note that the ring~$S$ in the \emph{classical} Grothendieck-Messing theory is with~$p$ locally nilpotent on it, see \cite[Chapter V, Theorem~(1.6)]{Messing72}.
	However, the ring in \thmref{thm:GM} has no local nilpotents. 
	We explain in the remark why the theory is still valid when~$S$ is a slightly ramified extension of~$\ZZ_{p}$, see \cite[Example~3.4]{Marrama22} also.
	
	Let~$S$ be a slightly ramified extension of~$\ZZ_{p}$. 
	Let~$G$ be a~$p$-divisible group over~$S/\varpi$, where~$\varpi$ is a uniformizer of~$S$. 
	We want to lift~$G$ to over~$S$. 
	We first consider the quotient~$S/\varpi$, where~$p$ is a nilpotent. 
	So,~$p$ is locally nilpotent on~$S/\varpi S$.
	Let~$I:=\varpi S$ be an ideal in~$S$.  
	Take~$S_{n}:=S/I^{n+1}$, and then we obtain a~$p$-divisible group~$G_{n}:=G\otimes S_{n}$ over~$S_{n}$ for any~$n\in\ZZ_{\geq 0}$. 
	We may apply \cite[Chapter V, Theorem~(1.6)]{Messing72} on~$G_{n}$ over~$S_{n}$ for any~$n\in\ZZ_{\geq 0}$. 
	
	Since we assume that~$S$ is slightly ramified, it has the so-called \emph{PD-structure~$\gamma_{n}$} of~$(S_{n}, I^n/I^{n+1})$ for every~$n\in\ZZ_{\geq 0}$, see \cite[Chapter III]{Messing72}.
	So,~$G_{n}$ satisfies the Grothendieck-Messing theory over~$S_{n}$ for any~$n\in\ZZ_{\geq 1}$. 
	That is,~$G_{n}$ is lifted to~$\mathcal{G}_{n}$ over~$S_{n-1}$ for any~$n\in\ZZ_{\geq 1}$.  
	Therefore, we obtain a~$p$-divisible group~$\mathcal{G}:=\varprojlim \mathcal{G}_{n}$ over~$S$, which is a lifting of~$G$ by \cite[Chapter II, Lemma~(4.16)]{Messing72}. 
	
	Moreover, let~$L$ be a CM field, such that~$[L:\QQ]$ equals the height of~$G$. 
	Assume that~$O_{L}\otimes\ZZ_{p}\embed\End_{S/\varpi}(G)$, and the~$O_{L}\otimes\ZZ_{p}$-action preserves the filtration in the Grothendieck-Messing theory of each~$G_{n}$ for~$n\in\ZZ_{\geq 1}$. 
	Then~$\mathcal{G}$ is an~$O_{L}$-linear CM~$p$-divisible group. 
	
	Therefore, if we find a filtration~$\textrm{Fil}_{w}$ as in \thmref{thm:GM}, then it will satisfy the previous argument, such that there is an~$O_{L}$-linear CM~$p$-divisible group~$A_{\pi}[p^{\infty}]_{w}$, which lifts~$A_{\pi}[p^{\infty}]_{w}$ to over~$S$.   
\end{rmk}

\begin{pf}\label{rmk:GM7}[Second proof of \thmref{thm:ApisCML}.]
Note that we have a decomposition~$A_{\pi}[7^\infty]=A_{\pi}[7^\infty]_{w}\times A_{\pi}[7^\infty]_{\overline{w}}$ over~$\FF_{7}$, where~$A_{\pi}[7^\infty]_{w}$ (resp.~$A_{\pi}[7^\infty]_{\overline{w}}$) is a~$7$-divisible group with (contravariant) Dieudonn\'e module~$M_{\pi,w}$ (resp.~$M_{\pi,\overline{w}}$). 
Same for~$M_{0}$.
In the following, we will show that~$A_{\pi}[7^{\infty}]_{w}$ and~$A_{\pi}[7^{\infty}]_{\overline{w}}$ have a CM-liftings, respectively.
Note that the Lie algebra of an abelian scheme is isomorphic to the one of its~$p$-divisible group. 
This will then imply the CM liftability of~$A_{\pi}/\FF_{p}$ by Serre-Tate deformation theory on abelian schemes, see \thmref{thm:SerreTate} and \thmref{thm:algebrai}.

Let~$\textrm{F}$ be the Frobenius operator on~$M_{\pi}$. 
Note that~$\textrm{F}^2=7\zeta_{}^{2}$, where~$\zeta_{}\in O_{L_{w}}^{\times}$ is a third root of unity. 
Let~$e_{1}\notin \textrm{F}M_{\pi,w}$, and let~$\{e_{1},e_{2}:=\textrm{F}e_{1}\}$ be~$\ZZ_{7}$-generators of~$M_{\pi,w}$.
Then,~$\textrm{F}e_{2}=7\zeta^2 e_{1}$, and hence~$\textrm{F}$ corresponds to the matrix~$\begin{pmatrix}
	0 & 7\zeta^2 \\
	1 & 0
\end{pmatrix}$ over~$O_{L_{w}}$, which has the characteristic polynomial~$T^2-7\zeta_{}^2\in\ZZ_{7}[T]$.
Let~$S:=\ZZ_{7}[\sqrt{7}]$ with uniformizer~$\varpi$, 
and let~$\textrm{Fil}_{w}$ be the~$S$-submodule of~$M_{\pi,w}\otimes_{\ZZ_{7}}S$, which is generated by~$e_{2}+\sqrt{7}\zeta_{}e_{1}$.
Note that~$\ker \textrm{F}$ modulo~$7\ZZ_{7}$ is~$\FF_{7}$-isomorphic to~$\ker \textrm{F}\otimes_{\ZZ_{7}}S = \ker (\textrm{F}\otimes\textrm{id}_{S})$ modulo~$\sqrt{7}\zeta_{}S$, and hence~$\ker (\textrm{F}\otimes\textrm{id}_{S})$ is generated by~$e_{2}$ modulo~$\sqrt{7}\zeta_{}S$.
On the other hand, from the choice of~$\textrm{Fil}_{w}$, we see that~$\textrm{Fil}_{w}$ modulo~$\sqrt{7}\zeta_{} S$ is generated by~$e_{2}$.
So,\[\textrm{Fil}_{w}\otimes_{S}(S/\varpi S) = \ker\textrm{F}\otimes_{\ZZ_{7}}(S/\varpi S),\] which is the condition~$(2)$ for \thmref{thm:GM}. 

Note that~$O_{L_{w}}=\ZZ_{7}[\pi_{w}]=\ZZ_{7}[\sqrt{7}] =S$, where~$\pi_{w}$ is a uniformizer of~$O_{L_{w}}$, which is equal to~$\sqrt{7}$, up to a unit~$\zeta$ in~$\ZZ_{7}$.
Moreover,~$\pi_{w}$ acts as~$\textrm{F}$ on~$M_{\pi,w}$, and~$\textrm{Fil}_{w}$ is stable under~$\textrm{F}$ by construction.
This implies that~$\textrm{Fil}_{w}$ is stable under the~$O_{L_{w}}$-action, which is the condition~$(3)$ for \thmref{thm:GM}.  

Lastly, we want to check if~$(M_{\pi,w}\otimes_{\ZZ_{7}}S)/\textrm{Fil}_{w}$ is~$S$-torsion-free.
It is sufficient to consider a pure tensor.
Let~$0\neq m\in M_{\pi,w}\otimes_{\ZZ_{7}}S$ such that~$sm \in\textrm{Fil}_{w}$ for some~$s\in S-\{0\}$.
We can write~$m=a_{1}e_{1}+a_{2}e_{2}$ for some~$a_{1},a_{2}\in S$.
Then,~$sm = sa_{1}e_{1}+sa_{2}e_{2} = a_{3}e_{2}+a_{3}\sqrt{7}\zeta_{}e_{1}$ for some~$a_{3}$ because~$sm\in\textrm{Fil}_{w}$ by assumption.
Since~$\{e_{1},e_{2}\}$ are~$\ZZ_{7}$-generators for~$M_{\pi,w}$ (and hence~$S$-generators for~$M_{\pi,w}\otimes S$), and~$S$ is an integral domain, we see that~$a_{1}=a_{2}\sqrt{7}\zeta_{}$.
So,~$m=a_{2}(\sqrt{7}\zeta_{}e_{1}+e_{2})$ is in~$\textrm{Fil}_{w}$.
This is the condition~$(1)$ for \thmref{thm:GM}.

We may use the same argument for~$M_{\pi,\overline{w}}$, and hence we obtain an~$S$-module~$\textrm{Fil}_{\overline{w}}$, which satisfies the conditions in Grothendieck-Messing theory.
Note that~$A_{\pi}[7^{\infty}]=A_{\pi}[7^{\infty}]_{w}\times A_{\pi}[7^{\infty}]_{\overline{w}}$ because~$M_{\pi}=M_{\pi,w}\oplus M_{\pi,\overline{w}}$ as an~$O_{L}\otimes_{\ZZ}\ZZ_{7}$-module.
In particular, with~$\textrm{Fil}_{w}$ and~$\textrm{Fil}_{\overline{w}}$, we know that~$A_{\pi}[7^{\infty}]_{w}$ and~$A_{\pi}[7^{\infty}]_{\overline{w}}$ have an~$O_{L_{w}}$-lifting and an~$O_{L_{\overline{w}}}$-lifting over~$S=\ZZ_{7}[\sqrt{7}]$, respectively, by Grothendieck-Messing theory.
Therefore, we may lift~$A_{\pi}[7^{\infty}]$ to~$S$, and the lifting~$\mathcal{A}[p^{\infty}]$ over~$S$ has CM by~$O_{L}\otimes_{\ZZ}\ZZ_{7}$.
Hence, we obtain a formal CM lifting~$\mathcal{A}$ of~$A_{\pi}/\FF_{7}$ over~$S$, which has smCM by~$L$, follows \thmref{thm:SerreTate}. 
In particular,~$\mathcal{A}[p^{\infty}]$ is the~$p$-divisible group of~$\mathcal{A}$. 
Note that~$\textrm{Lie}(\mathcal{A}[p^{\infty}])$ is isomorphic to~$\textrm{Lie}(\mathcal{A})$, and that the~$L$-linear action on~$\textrm{Lie}(\mathcal{A}[p^{\infty}])$ is given by a~$p$-adic~$\overline{\QQ}_{p}$-valued CM type of~$L\otimes\QQ_{p}$ by \cite[Lemma~3.7.1.3]{CCO14}.
By \eqref{eq:padicCM}, we see that the~$L$-action on~$\textrm{Lie}(\mathcal{A})$ is given by a CM type of~$L$.   
So, the formal CM lifting~$\mathcal{A}$ is algebraizable by \thmref{thm:algebrai}. 



The alternative proof is now complete. \qed
\end{pf}

\section{Simple superspecial abelian surfaces over prime fields.}\label{chap:ss}
Recall Weil-$p$ numbers of our concern~$\pi$ from the beginning of \secref{ex:ss ab surf}.
In this section, we prove the CM liftability of every superspecial abelian surface over~$\FF_{p}$ in the isogeny class determined by any Weil-$p$ number of our concern for any prime~$p>0$.

In some cases, we first show that an superspecial abelian surface~$A_{0}/\FF_{p}$ has a CM lifting after base change to~$\FF_{p^2}$.
Then, with \thmref{lem:changeofcoe}, we may deduce that it has a CM lifting without a base-change.

However, we treat the case~$\pi=\sqrt{3}\zeta_{3}$ separately in \propref{rem:p=3} because the argument is slightly different from other Weil-$p$ number of our concern. 


\subsection{Non-real Weil-$p$ numbers.}\label{sec:nonreal}

We consider a simple superspecial abelian surface~$A_{\pi}$ over~$\FF_{p}$, where~$\pi$ is any Weil-$p$ number of our concern.

We first check if~$A_{\pi}/\FF_{p}$ satisfies RRC. 

\begin{lem}\label{cor:STholdpi}
Let~$\pi$ be a Weil-$p$ number of our concern.
Let~$L:=\QQ(\pi)$ be a CM field with maximal totally real subfield~$L_{0}$, and let~$A_{\pi}/\FF_{p}$ be a simple supersingular abelian surface in the isogeny class determined by~$\pi$, which is obtained by the Honda-Tate theorem.
Then, there is a~$\overline{\QQ}_{p}$-valued CM type~$\Phi$ of~$L$, such that~$(L,\Phi)$ satisfies the Shimura-Taniyama formula for any embedding~$\overline{\QQ}\embed\overline{\QQ}_{p}$. 
\end{lem}
\begin{pf}
We use the notation as in \defref{def:emb}.
Let~$v$ be the place in~$L_{0}$ lying above~$p$.
We have the following three cases.
	\begin{enumerate} 
		\item Suppose that~$L/L_{0}$ is inert or ramified at~$v$.
		In this case, let~$w$ be the place in~$L$ lying above~$p$.
		Then,~$\Phi=\Phi_{w}$ for any CM type~$\Phi$ and any embedding~$\overline{\QQ}\embed\overline{\QQ}_{p}$.
		In particular, we see that~$|\Phi_{w}|/[L_{w}:\QQ_{p}]=1/2$.
		\item Suppose that~$L/L_{0}$ splits completely at~$v$.
		In this case, let~$w,\overline{w}$ be the two distinct places in~$L$ lying above~$p$.
		There exists a~$\overline{\QQ}_{p}$-valued CM type~$\Phi$ of~$L$ such that~$|\Phi_{w}|=|\Phi_{\overline{w}}|=1$.
		Indeed, let~$\varphi_{1},\dots,\varphi_{4}$ be distinct elements in~$\Hom(L,\overline{\QQ}_{p})$ for any embeddings~$L\embed\overline{\QQ}$ and~$\overline{\QQ}\embed\overline{\QQ}_{p}$.
		Let~$c\in\Aut(L/{\QQ})$ as in \defref{def:CMtype}.
		Without loss of generality, we may assume~$\varphi_{1}=\varphi_{4}\circ c$ and~$\varphi_{2}=\varphi_{3}\circ c$.
		We may assume that~$\varphi_{1}$ induces~$w$ and~$\varphi_{2}$ induces~$\overline{w}$.
		By assumption on~$\varphi_{3},\varphi_{4}$, we see that~$\varphi_{3}$ induces~$w$, and~$\varphi_{4}$ induces~$\overline{w}$ because~$\overline{w}=cw$.
	Take~$\Phi:=\{\varphi_{1},\varphi_{2}\}$ or~$\{\varphi_{3},\varphi_{4}\}$.
	With the~$\overline{\QQ}_{p}$-valued CM type~$\Phi$, we have~$|\Phi_{w}|/[L_{w}:\QQ_{p}]=|\Phi_{\overline{w}}|/[L_{\overline{w}}:\QQ_{p}]=1/2$.
	This completes the proof.\qed
	\end{enumerate}
\end{pf}

After \lemref{cor:STholdpi}, for any Weil-$p$ number of our concern, any simple superspecial abelian surface over~$\FF_{p}$, in the isogeny class determined by~$\pi$, satisfies the Shimura-Taniyama formula (for any embedding~$\overline{\QQ}\embed\overline{\QQ}_{p}$) for some CM type of~$L$ 
because the choice of such~$A_{\pi}/\FF_{p}$ is arbitrary in the proof of \lemref{cor:STholdpi}. 

To see whether~$A_{\pi}/\FF_{p}$ satisfies RRC, it remains to check whether the reflex field related to such a CM type satisfies the reflex field condition.
We shall see that not all CM types satisfy the reflex field condition below.
\begin{lem}\label{lem:noRRC}
Let~$\pi$ be a Weil-$p$ number of our concern.
	Let~$A_{\pi}/\FF_{p}$ be a simple supersingular abelian surface in the isogeny class determined by~$\pi$.
	Then,~$A_{\pi}/\FF_{p}$ satisfies RRC with some CM types of~$L$. 
\end{lem}
\begin{pf}
Let~$L:=\QQ(\pi)$, and let~$A_{\pi}/\FF_{p}$ be a simple superspecial abelian surface in the isogeny class determined by~$\pi$.
Note that the reflex field condition holds automatically when~$v$ is ramified or splits completely in~$L$, since the inertia degree at~$w$ over~$p$ is~$1$ in these cases.
In the following, we compute the reflex fields (as a subfield of~$\overline{\QQ}$) attached to all CM types of~$L$ when~$v$ is inert in~$L$.
In this case, the Shimura-Taniyama formula holds for any CM type of~$L$ and for any embedding~$\overline{\QQ}\embed\overline{\QQ}_{p}$ as we have seen in the proof of \lemref{cor:STholdpi}.
And hence,~$A_{\pi}/\FF_{p}$ satisfies RRC if and only if it satisfies the reflex field condition.
\begin{enumerate}
	\item Let~$\pi=\sqrt{p}\zeta_{8}$ and~$p\equiv 3$ mod~$4$.
	In this case,~$L=\QQ(\sqrt{2p},\sqrt{-1})$.
	We compute the reflex field for each CM type. 
	Let \begin{align*}
    \varphi_{1}: \begin{cases}
        \sqrt{2p} \mapsto \sqrt{2p}\\
        \sqrt{-1} \mapsto -\sqrt{-1},
    \end{cases} & \varphi_{2}: \begin{cases}
        \sqrt{2p} \mapsto -\sqrt{2p}\\
        \sqrt{-1}\mapsto \sqrt{-1},
    \end{cases} \end{align*}\begin{align*}
    \varphi_{3}: \begin{cases}
        \sqrt{2p} \mapsto -\sqrt{2p}\\
        \sqrt{-1}\mapsto -\sqrt{-1}, \end{cases}& \varphi_{4}: \begin{cases}
        \sqrt{2p} \mapsto \sqrt{2p}\\
        \sqrt{-1}\mapsto \sqrt{-1}
    \end{cases}
\end{align*}be the elements in~$\Hom(L,\overline{\QQ}_{})$. 
After computation following \defref{def:reflexfield}, we have the following table.

\begin{table}[H]
	\centering
    \begin{tabular}{|c|c|}\hline
    	\rm{CM type}~$\Phi$ & \rm{reflex field} \\ \hline
    	$\{\varphi_{1},\varphi_{2}\}$ & $\QQ(\sqrt{-2p})$  \\ \hline
    	$\{\varphi_{1},\varphi_{3}\}$ & $\QQ(\sqrt{-1})$ \\ \hline
    	$\{\varphi_{2},\varphi_{4}\}$ & $\QQ(\sqrt{-1})$  \\ \hline
    	$\{\varphi_{3},\varphi_{4}\}$ & $\QQ(\sqrt{-2p})$  \\\hline
    \end{tabular}
\caption{Every CM type of~$L$ and its associated reflex field when~$\pi=\sqrt{p}\zeta_{8}$ and~$p\equiv 3$ mod~$4$.}
\end{table}
By assumption on~$p$, we see that the place in~$\QQ(\sqrt{-1})$ is inert over~$p$.
In particular, when we regard~$\QQ(\sqrt{-1})$ as a subfield of~$\overline{\QQ}_{p}$ for any embedding~$\overline{\QQ}\embed\overline{\QQ}_{p}$, the induced residue field is~$\FF_{p^2}$.
That is, when~$\Phi=\{\varphi_{1},\varphi_{3}\}$ or~$\{\varphi_{2},\varphi_{4}\}$, we see that~$A_{\pi}/\FF_{p}$ does not satisfy RRC.
On the other hand, when~$\Phi=\{\varphi_{1},\varphi_{2}\}$ or~$\{\varphi_{3},\varphi_{4}\}$, then~$A_{\pi}/\FF_{p}$ satisfies RRC by the similar argument.

\item Let~$\pi=\sqrt{p}\zeta_{3}$ and~$p\equiv 5,11$ mod~$12$.
	In this case,~$L=\QQ(\sqrt{p},\sqrt{-3})$.
	We compute the reflex field for each CM type by \defref{def:reflexfield}.
	Let \begin{align*}
    \varphi_{1}: \begin{cases}
        \sqrt{p} \mapsto \sqrt{p}\\
        \sqrt{-3} \mapsto -\sqrt{-3},
    \end{cases} & \varphi_{2}: \begin{cases}
        \sqrt{p} \mapsto -\sqrt{p}\\
        \sqrt{-3}\mapsto \sqrt{-3},
    \end{cases} \end{align*}\begin{align*}
    \varphi_{3}: \begin{cases}
        \sqrt{p} \mapsto -\sqrt{p}\\
        \sqrt{-3}\mapsto -\sqrt{-3}, \end{cases}& \varphi_{4}: \begin{cases}
        \sqrt{p} \mapsto \sqrt{p}\\
        \sqrt{-3}\mapsto \sqrt{-3}
    \end{cases}
\end{align*}be the elements in~$\Hom(L,\overline{\QQ}_{})$.
Again, following \defref{def:reflexfield}, we obtain the following table. 
\begin{table}[H]
	\centering
    \begin{tabular}{|c|c|}\hline
    	\rm{CM type}~$\Phi$ & \rm{reflex field} \\ \hline
    	$\{\varphi_{1},\varphi_{2}\}$ & $\QQ(\sqrt{-3p})$  \\ \hline
    	$\{\varphi_{2},\varphi_{4}\}$ & $\QQ(\sqrt{-3})$ \\ \hline
    	$\{\varphi_{3},\varphi_{4}\}$ & $\QQ(\sqrt{-3p})$  \\ \hline
    	$\{\varphi_{1},\varphi_{3}\}$ & $\QQ(\sqrt{-3})$  \\\hline
    \end{tabular}
\caption{Every CM type of~$L$ and its associated reflex field when~$\pi=\sqrt{p}\zeta_{3}$ and~$p\equiv 5,11$ mod~$12$.}\label{table:pzeta3}
\end{table}
Suppose~$p\geq 5$.
The residue field of the reflex field at~$p$ is~$\FF_{p^2}$ when~$\Phi=\{\varphi_{2},\varphi_{4}\}$ or~$\{\varphi_{1},\varphi_{3}\}$ for any embedding~$\overline{\QQ}\embed\overline{\QQ}_{p}$. 
In this case,~$A_{\pi}/\FF_{p}$ does not satisfy RRC.
Otherwise,~$A_{\pi}/\FF_{p}$ satisfies RRC by the similar argument.

Suppose~$p=3$.
Then, we have the Legendre symbol~$\left(\frac{-1}{3}\right)=-1$, and hence the residue field of the reflex field at~$3$ is~$\FF_{9}$ when~$\Phi=\{\varphi_{1},\varphi_{2}\}$ or~$\{\varphi_{3},\varphi_{4}\}$ for any embedding~$\overline{\QQ}\embed\overline{\QQ}_{3}$.
In this case,~$A_{\pi}/\FF_{3}$ does not satisfy RRC.
Otherwise,~$A_{\pi}/\FF_{3}$ satisfies RRC by the similar argument.

	\item Let~$\pi=\sqrt{p}\zeta_{12}$ and~$p\equiv 5,11$ mod~$12$.
	In this case,~$L=\QQ(\sqrt{-p},\sqrt{-3})$.
	We compute the reflex field for each CM type by \defref{def:reflexfield}.
	Let \begin{align*}
    \varphi_{1}: \begin{cases}
        \sqrt{-p} \mapsto -\sqrt{-p}\\
        \sqrt{-3} \mapsto \sqrt{-3},
    \end{cases} & \varphi_{2}: \begin{cases}
        \sqrt{-p} \mapsto \sqrt{-p}\\
        \sqrt{-3}\mapsto -\sqrt{-3},
    \end{cases} \end{align*}\begin{align*}
    \varphi_{3}: \begin{cases}
        \sqrt{-p} \mapsto -\sqrt{-p}\\
        \sqrt{-3}\mapsto -\sqrt{-3}, \end{cases}& \varphi_{4}: \begin{cases}
        \sqrt{-p} \mapsto \sqrt{-p}\\
        \sqrt{-3}\mapsto \sqrt{-3}
    \end{cases}
\end{align*}be the elements in~$\Hom(L,\overline{\QQ}_{})$.
After computation following \defref{def:reflexfield}, we obtain the following table.
\begin{table}[H]
	\centering
    \begin{tabular}{|c|c|}\hline
    	\rm{CM type}~$\Phi$ & \rm{reflex field}  \\ \hline
    	$\{\varphi_{1},\varphi_{4}\}$ & $\QQ(\sqrt{-3})$  \\ \hline
    	$\{\varphi_{1},\varphi_{3}\}$ & $\QQ(\sqrt{-p})$ \\ \hline
    	$\{\varphi_{2},\varphi_{4}\}$ & $\QQ(\sqrt{-p})$  \\ \hline
    	$\{\varphi_{2},\varphi_{3}\}$ & $\QQ(\sqrt{-3})$  \\\hline
    \end{tabular}
\caption{Every CM type of~$L$ and its associated reflex field when~$\pi=\sqrt{p}\zeta_{12}$ and~$p\equiv 5,11$ mod~$12$.}
\end{table}
The residue field of the reflex field at~$p$ is~$\FF_{p^2}$ when~$\Phi=\{\varphi_{1},\varphi_{4}\}$ or~$\{\varphi_{2},\varphi_{3}\}$ for any embedding~$\overline{\QQ}\embed\overline{\QQ}_{p}$. 
In this case,~$A_{\pi}/\FF_{p}$ does not satisfy RRC.
Otherwise,~$A_{\pi}/\FF_{p}$ satisfies RRC by the similar argument.
\end{enumerate}
This completes the proof.\qed
\end{pf}

Now, we state and prove the main theorem. 
\begin{thm}\label{cor:surfacesCML}
	Let~$\pi$ be a Weil-$p$ number of our concern except~$\sqrt{{3}}\zeta_{3}$.
	Let~$A_{\pi}/\FF_{p}$ be a simple superspecial abelian surface in the isogeny class determined by~$\pi$.
	Then,~$A_{\pi}/\FF_{p}$ has a CML.
\end{thm}
\begin{pf}
Let~$L:=\QQ(\pi)$ be a CM field with maximal totally real subfield~$L_{0}$.
Let~$A_{\pi}/\FF_{p}$ be a simple superspecial abelian surface over~$\FF_{p}$, in the isogeny class determined by~$\pi$, which admits smCM by~$L$.
Note that when~$\pi=\sqrt{p}\zeta_{8}$ (resp.~$\sqrt{p}\zeta_{3}$,~$\sqrt{p}\zeta_{12}$ or~$\sqrt{2}\zeta_{3},\sqrt{2}\zeta_{12},\pm\sqrt{2}\zeta_{24}$), $p\equiv 3$ mod~$4$ (resp.~$p\equiv 5,11$ mod~$12$ or~$p=2$), the finite field~$\FF_{p^2}$ (resp.~$\FF_{p}$) is the smallest ground field, where we can consider Lie types of simple superspecial abelian surfaces in the isogeny class determined by~$\pi$, see \remref{rmk:lietypeuse}. 


	By \lemref{lem:noRRC}, there exists at least one CM type of~$L$ such that~$A_{\pi}/\FF_{p}$ satisfies RRC.
	Therefore, there exists a simple abelian variety~$A_{0}$ over~$\FF_{p}$, admitting smCM by~$L$ and admitting a CM lifting defined over a~$p$-adic field~$\widetilde{E}$, such that~$A_{0}$ is~$L$-linearly isogenous to~$A_{\pi}$ over~$\FF_{p}$ by \cite[Theorem~2.5.3]{CCO14}. 
	In particular, the associated~$p$-divisible group~$A_{0}[p^{\infty}]$ of~$A_{0}$ is~$O_{L}$-linear.
	
	Note that the~$p$-divisible group~$A_{\pi}[p^{\infty}]$, which is associated with~$A_{\pi}$, is~$O_{L}$-linear.
	Then,~$A_{\pi}[p^{\infty}]$ and~$A_{0}[p^{\infty}]$ over~$\FF_{p}$ are~$O_{L}\otimes\ZZ_{p}$-linearly~$\FF_{p}$-isogenous.
	The Lie types in the isogeny class determined by~$\pi$, after base change at most to~${\FF}_{p^{2}}$, are the same by \propref{lem:sameLietype}.
	\begin{enumerate}	
		\item Let~$\pi=\sqrt{p}\zeta_{8}$ if~$p\equiv 3$ mod~$4$, and~$\pi=\sqrt{p}\zeta_{3}$ or~$\sqrt{p}\zeta_{12}$ if~$p\equiv 5, 11$ mod~$12$ be a Weil-$p$ number.
	In this case,~$L=\QQ(\pi)/L_{0}$ is inert.
	Let~$\varphi:A_{0}\to A_{\pi}$ be the~$L$-linear~$\FF_{p}$-isogeny obtained by \cite[Theorem~2.5.3]{CCO14}.
	Let~$\tilde{\iota}_{A_{0}}$ be a CM structure of~$A_{0}$, which is reduced from the one of the CM lifting of~$A_{0}$.
	Then, the~$L$-linearity of~$\varphi$ is with respect to CM structures~$\tilde{\iota}$ of~$A_{\pi}$ and~$\tilde{\iota}_{A_{0}}$ of~$A_{0}$.
	We base change~$\varphi$ to~$\FF_{p^2}$, and denote it by~$\varphi\otimes\textrm{id}_{\FF_{p^2}}:A_{0}\otimes\FF_{p^2}\to A_{\pi}\otimes\FF_{p^2}$.
	Moreover, we have an induced morphism~$(\varphi\otimes\textrm{id}_{\FF_{p^2}})[p^{\infty}]$ between the corresponding~$O_{L}$-linear CM~$p$-divisible groups over~$\FF_{p^2}$. 
Replace~$\FF_{7}$ by~$\FF_{p^2}$ in the first proof of \thmref{thm:ApisCML}, we see that~$A_{\pi}\otimes\FF_{p^2}$ has a CM lifting over~$\widetilde{E}$. 
Furthermore,~$A_{\pi}/\FF_{p}$ has a CM lifting by \thmref{lem:changeofcoe}.
	\item Let~$\pi$ be a Weil-$p$ number, such that~$L=\QQ(\pi)/L_{0}$ splits.
	In this case, the proof is the same as the first proof of \thmref{thm:ApisCML} after replacing~$\FF_{7}$ by~$\FF_{p}$.
	\end{enumerate}
The proof is complete because~$A_{\pi}/\FF_{p}$ is arbitrarily chosen. 
	 \qed
\end{pf}

When~$p=3$, the Weil-$3$ number~$\sqrt{3}\zeta_{3}$ is of our concern, which we have excluded in the previous theorem because the endomorphism ring~$\End_{\FF_{3}}(A_{\pi})$ is too small.
Also, it is not enough to do the base-change to~$\FF_{p^2}$ as we will see later.
We consider it now in the following proposition.
\begin{prop}\label{rem:p=3}
	Let~$\pi=\sqrt{3}\zeta_{3}$ be the Weil-$3$ number of our concern.
	Any simple superspecial abelian surface~$A_{\pi}$ over~$\FF_{3}$ in the isogeny class determined by~$\pi$ has a CM lifting.
\end{prop}
\begin{pf}
Let~$\pi:=\sqrt{3}\zeta_{3}$, and let~$L:=\QQ(\pi)=\QQ(\zeta_{12})$ be a CM field with ring of integers~$O_{L}:=\ZZ[\zeta_{12}]$ and with maximal totally real subfield~$L_{0}:=\QQ(\sqrt{3})$.
	Let~$A_{\pi}/\FF_{3}$ be a simple superspecial abelian surface in the isogeny class determined by~$\pi$.
	The Frobenius endomorphism~$\textrm{Frob}_{A_{\pi},3}$ of~$A_{\pi}$ over~$\FF_{3}$ is one of elements in~$\{\pm\sqrt{3}\zeta_{3}, \pm\sqrt{3}\zeta_{3}^2\}$ because they satisfy the irreducible polynomial~$T^4+3T^2+9$.

	By \cite[Equation~(5.9)]{XYY16}, we know there are exactly two (up to~$\FF_{3}$-isomorphism) simple superspecial abelian surfaces~$A_{\pi,1},A_{\pi,2}$ over~$\FF_{3}$ in the isogeny class determined by~$\pi=\sqrt{3}\zeta_{3}$.
	Without loss of generality, we may assume that~$A_{\pi,2}$ is principal, and~$A_{\pi,1}$ is not because the class number of~$O_{L}$ is~$1$. 
	Then, we know that the endomorphism ring of~$A_{\pi,1}$ over~$\FF_{3}$ is isomorphic to~$R_{\textrm{sp}}:=\ZZ[\pi, 3/\pi^{}, \pi^2/3]$ by \cite[Theorem~6.1]{W69}.
	Since~$A_{\pi,1}/\FF_{3}$ satisfies RRC by \lemref{lem:noRRC}, we know that~$A_{\pi,1}$ is~$L$-linearly~$\FF_{3}$-isogenous to an abelian surface~$\widehat{A}$ over~$\FF_{3}$, admitting a CM lifting over a~$p$-adic field~$\widetilde{E}$ by \cite[Theorem~2.5.3]{CCO14}. 
	Moreover, the associated~$3$-divisible group~$\widehat{A}[3^{\infty}]$ over~$\FF_{3}$ is~$O_{L}$-linear.
	
	We claim that~$\widehat{A}$ is principal, and hence is~$\FF_{3}$-isomorphic to~$A_{\pi,2}$.
	Let~$R_{\textrm{ss}}:=\ZZ[\pi,3\pi^{-1}]$ be a subring of~$O_{L}$, whose index~$[O_{L}:R_{\textrm{ss}}]$ is~$9$.
	Then, the inclusions~$R_{\textrm{ss}}\subseteq \End_{\FF_{3}}(\widehat{A})\subseteq O_{L}$ imply that~$O_{L}=\End_{\FF_{3}}(\widehat{A})$ because~$\widehat{A}[3^{\infty}]$ over~$\FF_{3}$ is~$O_{L}$-linear.
	As there is the unique (up to isomorphism) principal abelian surface over~$\FF_{3}$, which is~$A_{\pi,2}/\FF_{3}$ by assumption, we see that~$A_{\pi,2}/\FF_{3}$ has a CML over~$\widetilde{E}$ by \cite[Proposition~2.22]{BKM}.
	
	
	Now, we want to show that~$A_{\pi,1}/\FF_{3}$ has a CM lifting.
	Let~$M_{\pi,1}$ (resp.~$M_{\pi,2}$) 
	be the (contravariant) Dieudonn\'e module over~$\ZZ_{3}$ of~$A_{\pi,1}$ (resp.~$A_{\pi,2}$).
	Since~$A_{\pi,1}/\FF_{3}$ has a Weil-$3$ number~$\sqrt{3}\zeta_{3}$, we see that~$A_{\pi,1}\otimes\FF_{3^6}$ is~$\FF_{3^6}$-isogenous to~$E^2$, where~$E$ is a supersingular elliptic curve over~$\FF_{3^6}$ such that~$\End_{\FF_{3^6}}(E)=\End_{\overline{\FF}_{3}}(E\otimes_{\FF_{3^6}}\overline{\FF}_{3})$ by \cite[Theorem~2(d5)]{T66}. 
	Since~$\End_{\FF_{3^6}}^{0}(A_{\pi,1}\otimes\FF_{3^6})\cong \End_{\FF_{3^6}}^{0}(E^2)$, we have\[\End_{\FF_{3^6}}(A_{\pi,1}\otimes\FF_{3^6})=\End_{\overline{\FF}_{3}}(A_{\pi,1}\otimes\overline{\FF}_{3}).\]
	Similarly,~$\End_{\FF_{3^6}}(A_{\pi,2}\otimes\FF_{3^6})=\End_{\overline{\FF}_{3}}(A_{\pi,2}\otimes\overline{\FF}_{3})$.
	Since~$A_{\pi,1}$ and~$A_{\pi,2}$ are superspecial, they are isomorphic after base change to~$\overline{\FF}_{3}$.
	Therefore,\[\End_{\FF_{3^6}}(A_{\pi,1}\otimes\FF_{3^6}) \cong \End_{\FF_{3^6}}(A_{\pi,2}\otimes\FF_{3^6}).\]
	In particular, we see~$O_{L}\otimes\ZZ_{p}\embed \End_{W(\FF_{3^6})[\textrm{F,V}]}(M_{\pi,1}\otimes W(\FF_{3^6}))$, which implies~$M_{\pi,1}\otimes W(\FF_{3^6}) \cong O_{L}\otimes W(\FF_{3^6})$ by \cite[Proposition~1.4.3.9]{CCO14}.
	Hence,~$M_{\pi,1}\otimes W(\FF_{3^6})$ is isomorphic to~$M_{\pi,2}\otimes W(\FF_{3^6})$.
	Moreover, this implies that there exists a separable\emph{~$\FF_{3^6}$-isogeny}~$\psi: A_{\pi,2}\otimes\FF_{3^6} \to A_{\pi,1}\otimes\FF_{3^6}$, which induced an isomorphism~$M_{\pi,1}\otimes W(\FF_{3^6})\cong M_{\pi,2}\otimes W(\FF_{3^6})$, by \cite[Theorem~1(b)]{T66} (or \cite[Theorem~7]{MW71}) and \cite[p.~525]{W69}.
	Note that~$A_{\pi,2}/\FF_{3}$ has a CML over~$\widetilde{E}$, and that there exists a finite extension~$\widetilde{F}$ of~$\widetilde{E}$ with residue field~$\FF_{3^6}$, see \cite[Chapter III, Theorem~2]{S13}.
	So,~$A_{\pi,2}\otimes\FF_{3^6}$ has a CML over~$\widetilde{F}$.
	This implies that~$A_{\pi,1}\otimes\FF_{3^6}$ has a CML over~$\widetilde{F}$ by \cite[Proposition~2.22]{BKM}.
	The proof is complete after applying \thmref{lem:changeofcoe}.\qed
	
\end{pf}
\subsection{Real Weil-$p$ numbers.} \label{sec:real}
Let~$\pi=\sqrt{p}$, and let~$A_{\pi}/\FF_{p}$ be a simple abelian surface over~$\FF_{p}$, in the isogeny class determined by~$\pi$, obtained by the Honda-Tate theorem.
Then,~$A_{\pi}/\FF_{p}$ is superspecial, see \cite[p.~528]{W69}.

As~$\End_{\FF_{p}}^{0}(A_{\pi})$ is a quaternion algebra over~$\QQ(\sqrt{p})$ by \textit{loc. cit.}, there might be more than one CM field, which can be embedded into it.
In this subsection, we find a CM field~$L$, with which~$A_{\pi}/\FF_{p}$ has a CML (after base change at most to~$\FF_{p^2}$) in \thmref{thm:liftingreal}.
In addition, we choose~$A_{\pi}/\FF_{p}$ arbitrarily, so statements in this subsection hold for any simple superspecial abelian surface over~$\FF_{p}$ in the isogeny class determined by~$\pi=\sqrt{p}$.

The existence of a CM field~$L$, such that~$A_{\pi}/\FF_{p}$ satisfies RRC with its CM type, is proved in the following lemma.
\begin{lem}\label{lem:realRRC}
Let~$\pi:=\sqrt{p}$ for any prime~$p>0$.
Let~$L:=\QQ(\sqrt{p}\zeta_{3})$ be a CM field with maximal totally real subfield~$L_{0}:=\QQ(\pi)$.
Let~$v$ be the place in~$L_{0}$ lying above~$p$, and let~$w$ be a place in~$L$ lying above~$v$. 
Then, a simple superspecial abelian surface~$A_{\pi}/\FF_{p}$, in the isogeny class determined by~$\pi$, satisfies RRC with~$(L,\Phi)$ for some CM type~$\Phi$ of~$L$.
\end{lem}
\begin{pf}
Let~$\pi':=\sqrt{p}\zeta_{3}$, and let~$A_{\pi'}/\FF_{p}$ be a simple superspecial abelian surface, in the isogeny class determined by~$\pi'$. 
Then, we have~$\End_{\FF_{p}}^{0}(A_{\pi'})=\QQ(\pi')$.
By \cite[Theorem~6.1]{W69}, we may assume that~$O_{L}$ is the endomorphism ring of~$A_{\pi'}/\FF_{p}$.
That is,\[O_{L}=\End_{\FF_{p}}(A_{\pi'}).\]
Let~$A_{\pi}/\FF_{p}$ be a simple superspecial abelian surface, in the isogeny class determined by~$\pi$.
Since~$\pi^3={\pi'}^3$, we see that~$A_{\pi}\otimes\FF_{p^3}$ is~$\FF_{p^3}$-isogenous to~$A_{\pi'}\otimes\FF_{p^3}$ by \cite[Theorem~1]{T66}.

Let~$\textrm{Frob}_{A_{\pi},p}$ (resp.~$\textrm{Frob}_{A_{\pi'},p}$) be the Frobenius endomorphism of~$A_{\pi}$ (resp.~$A_{\pi'}$) over~$\FF_{p}$.
Then, there is a~$\QQ$-algebra isomorphism~$\mathfrak{I}:\End_{\FF_{p^3}}^{0}(A_{\pi'}\otimes\FF_{p^3})\cong \End_{\FF_{p^3}}^{0}(A_{\pi}\otimes\FF_{p^3})$, which sends~$\textrm{Frob}_{A_{\pi'},p}^3$ to~$\textrm{Frob}_{A_{\pi},p}^3$.
Since~$\textrm{Frob}_{A_{\pi},p}^3=\pm \textrm{Frob}_{A_{\pi'},p}^3=\pm p\sqrt{p}$, we may assume that~$\mathfrak{I}$ is~$L_{0}$-linear.
As~$\QQ(\pi)=\QQ(\pi^3)$, we have~$\End_{\FF_{p}}^{0}(A_{\pi})=\End_{\FF_{p^3}}^{0}(A_{\pi}\otimes\FF_{p^3})$ by \cite[p.~528]{W69}.
Therefore, we obtain an~$L_{0}$-linear embedding\begin{equation}\label{eq:CMreal}
	{\iota}: L:=\QQ(\pi')=\End_{\FF_{p}}^{0}(A_{\pi'})\xhookrightarrow[]{\textrm{incl}} \End_{\FF_{p^3}}^{0}(A_{\pi'}\otimes\FF_{p^3})\overset{\mathfrak{I}}{\cong}\End_{\FF_{p}}^{0}(A_{\pi}),
\end{equation} 
which is given by the canonical inclusion.

We claim that~$L$ is a CM field, such that~$A_{\pi}/\FF_{p}$ satisfies RRC with some of its CM types.
Fix a suitable embedding~$\overline{\QQ}\embed\overline{\QQ}_{p}$ such that\[\sqrt{p}\mapsto \sqrt{p}, \sqrt{-3}\mapsto -\sqrt{-3}.\]
Let
\begin{align*}
    \varphi_{1}: \begin{cases}
        \sqrt{p} \mapsto \sqrt{p}\\
        \sqrt{-3} \mapsto -\sqrt{-3},
    \end{cases} & \varphi_{2}: \begin{cases}
        \sqrt{p} \mapsto -\sqrt{p}\\
        \sqrt{-3}\mapsto \sqrt{-3},
    \end{cases} \end{align*}\begin{align*}
    \varphi_{3}: \begin{cases}
        \sqrt{p} \mapsto -\sqrt{p}\\
        \sqrt{-3}\mapsto -\sqrt{-3}, \end{cases}& \varphi_{4}: \begin{cases}
        \sqrt{p} \mapsto \sqrt{p}\\
        \sqrt{-3}\mapsto \sqrt{-3}.
    \end{cases}
\end{align*} 
Without loss of generality, we may assume that~$\varphi_{1}$ induces~$w$, and hence~$\varphi_{2}$ induces~$\overline{w}$.
We have the following two cases.
\begin{enumerate}
	\item When~$v$ splits in~$L$, we have the table below, obtained by a similar computation for \autoref{ex1} in \secref{ex:1}.
	\begin{table}[H]
    \centering
    \begin{tabular}{|c |c| c|c|c|c|}\hline
        \rm{CM type} & 
         $w$ & $\overline{w}$ & \rm{slope} & \rm{reflex field} & $(|\Phi_{\bullet, w}|,|\Phi_{\bullet,\overline{w}}|)$ \\ \hline
        $\Phi_{1} =\left\{ \varphi_{1},\varphi_{2}\right\}$& $\varphi_{1}$ & $\varphi_{2}$& $(\frac{1}{2},\frac{1}{2})$ &$\QQ(\sqrt{-3p})$&$(1,1)$\\ \hline
        $\Phi_{2} =\left\{ \varphi_{2},\varphi_{4}\right\}$& $\emptyset$ & $\varphi_{2},\varphi_4{}$& $(0,1)$&$\QQ(\sqrt{-3})$&$(0,2)$\\ \hline
        $\Phi_{3} =\left\{ \varphi_{3},\varphi_{4}\right\}$& $\varphi_{3}$ & $\varphi_{4}$& $(\frac{1}{2},\frac{1}{2})$&$\QQ(\sqrt{-3p})$&$(1,1)$ \\ \hline
        $\Phi_{4} =\left\{\varphi_{1},\varphi_{3} \right\}$ &$\varphi_{1},\varphi_{3}$ & $\emptyset$& $(1,0)$&$\QQ(-3)$&$(2,0)$\\\hline
    \end{tabular}
    \caption{All CM types of~$L$ and their corresponding induced places, slopes, and reflex fields,~$(|\Phi_{\bullet, w}|,|\Phi_{\bullet,\overline{w}}|)$ when~$v$ splits in~$L$.} \label{table:realp}
\end{table}
We immediately see that~$A_{\pi}/\FF_{p}$ does not satisfy the Shimura-Taniyama formula with~$(L,\Phi_{2}), (L,\Phi_{4})$, but~$A_{\pi}/\FF_{p}$ satisfies it with~$(L,\Phi_{1}), (L,\Phi_{3})$ because~$A_{\pi}$ is superspecial.
Next, we want to check whether the reflex field condition holds for~$(L,\Phi_{1})$ or~$(L,\Phi_{3})$.
When~$p\geq 5$ or~$p=2$, we see that~$p$ divides the absolute discriminant of~$\QQ(\sqrt{-3p})$.
So,~$A_{\pi}/\FF_{p}$ satisfies RRC in this case.
When~$p=3$, we see that~$v$ is inert in~$L$, and hence we leave it to the following case. 

\item When~$v$ is inert in~$L$, recall the notations~$\Phi_{w}^{1},\Phi_{w}^{2}$ in \exref{ex:glt}.  We have the following table.
\begin{table}[H]
    \centering
    \begin{tabular}{|c |c| c|c|c|c|}\hline
       \rm{CM type} & 
         $\Phi_{\bullet,w}^{1}$ & $\Phi_{\bullet, w}^{2}$ & \rm{slope} & \rm{reflex field} & $(|\Phi_{\bullet, w}^1|,|\Phi_{\bullet,{w}}^2|)$\\ \hline
        $\Phi_{1} =\left\{ \varphi_{1},\varphi_{2}\right\}$& $\varphi_{1}$ & $\varphi_{2}$& $(\frac{1}{2},\frac{1}{2})$ &$\QQ(\sqrt{-3p})$& $(1,1)$\\ \hline
        $\Phi_{2} =\left\{ \varphi_{2},\varphi_{4}\right\}$& $\emptyset$ & $\varphi_{2}, \varphi_{4}$& $(\frac{1}{2},\frac{1}{2})$&$\QQ(\sqrt{-3})$& $(0,2)$\\ \hline
        $\Phi_{3} =\left\{ \varphi_{3},\varphi_{4}\right\}$& $\varphi_{3}$ & $\varphi_{4}$& $(\frac{1}{2},\frac{1}{2})$&$\QQ(\sqrt{-3p})$ &$(1,1)$\\ \hline
        $\Phi_{4} =\left\{\varphi_{1},\varphi_{3} \right\}$ &$\varphi_{1},\varphi_{3}$ & $\emptyset$& $(\frac{1}{2},\frac{1}{2})$&$\QQ(-3)$&$(2,0)$ \\\hline
    \end{tabular}
    \caption{All CM types of~$L$ and their corresponding induced places, slopes, and reflex fields,~$(|\Phi_{\bullet, w}^1|,|\Phi_{\bullet,{w}}^2|)$ when~$v$ is inert in~$L$.} \label{table:realp1}
\end{table}
We see that~$A_{\pi}/\FF_{p}$ satisfies the Shimura-Taniyama formula with~$(L,\Phi)$ for any~$\Phi\in\{\Phi_{i}\}_{i=1,\dots,4}$.
Next, we want to check whether the reflex field condition is satisfied.
When~$p\neq 3$, we see that~$A_{\pi}/\FF_{p}$ satisfies RRC with the CM types~$(L,\Phi_{1}), (L,\Phi_{3})$.
When~$p=3$, we see that it satisfies RRC with the CM types~$(L,\Phi_{2}),(L,\Phi_{4})$.\qed
\end{enumerate}
\end{pf}


\begin{rmk}
The reason we choose such~$\pi'=\sqrt{p}\zeta_{3}$ is that it is a Weil-$p$ number of our concern for any prime~$p>0$.
In addition, let~$A$ be a~$\ZZ$-order in~$\QQ(\sqrt{p})$, and let~$B$ be an~$A$-order in~$\QQ(\sqrt{p}\zeta_{3})$ such that~$B\cap \QQ(\sqrt{p})=A$
as in \cite[p.~671]{XYY19}.
When~$p\equiv 1$ mod~$4$, all such orders~$B$ with~$[B^{\times}:A^{\times}]>1$ lie in either~$\QQ(\sqrt{p},\sqrt{-1})$ or~$\QQ(\sqrt{p},\sqrt{-3})=\QQ(\sqrt{p}\zeta_{3})$, see \textit{loc. cit}.
\end{rmk}
We recall the following definition, which will be used to compute the Lie type of~$A_{\pi}\otimes\overline{\FF}_{p}$ for any~$p>0$. 
\begin{defn}\cite[Section~3.2]{XYY19}\label{def:optimal}
Let~$p>0$ be any prime.
Let~$L$ be a totally imaginary quadratic extension of~$L_{0}:=\QQ(\sqrt{p})$. 
Suppose that~$L$ is~$L_{0}$-linearly embedded into~$\End_{\FF_{p}}^{0}(A_{\pi})$.
Let~$A$ be a~$\ZZ$-order in~$L_{0}$.
Let~$B$ (resp.~$O$) be an~$A$-order in~$L$ (resp.~$\End_{\FF_{p}}^{0}(A_{\pi})$).
Suppose~$B\cap L_{0}=A$. 
	We denote\[\textrm{Emb}(B,O):=\{\varphi\in\Hom_{L_{0}}(L,\End_{\FF_{p}}^{0}(A_{\pi}))\,|\,\varphi(L)\cap O =\varphi(B)\}.\] 
An element in~$\textrm{Emb}(B,O)$ is called an \emph{optimal embedding from~$B$ to~$O$}.\end{defn}
\begin{ex}
	Let~$p>0$ be a prime such that~$p\equiv 3$ mod~$4$.
	Let~$\pi:=\sqrt{p}\zeta_{3}$, and let~$L:=\QQ(\pi)$ be a CM field with maximal totally real subfield~$L_{0}:=\QQ(\sqrt{p})$.
	Let~$A_{\pi}/\FF_{p}$ be a simple superspecial abelian surface in the isogeny class determined by~$\pi$.
	Let~$\textrm{Frob}_{A_{\pi},p}$ be the Frobenius endomorphism of~$A_{\pi}/\FF_{p}$.
	Note that~$\End_{\FF_{p}}(A_{\pi})$ is a~$\ZZ[\sqrt{p}]$-order because it contains a subring~$R_{\textrm{ss}}:=\ZZ[\textrm{Frob}_{A_{\pi},p},p\textrm{Frob}_{A_{\pi},p}^{-1}]$, as~$A_{\pi}$ is superspecial and hence supersingular.
	Note that~$A_{\pi}/\FF_{p}$ is simple and admits smCM by CM field~$L$.
	So,~$\End_{\FF_{p}}^{0}(A_{\pi})=L$, which implies that~$\Hom_{L_{0}}(L,\End_{\FF_{q}}^{0}(A_{0}))$ equals~$\Gal(L/L_{0})$.
	Take~$A:=\ZZ[\sqrt{p}]$,~$B:=O_{L}$ the ring of integers of~$L$, and~$O:=\End_{\FF_{p}}(A_{\pi})$.
	Then, the condition for~$\varphi\in\Hom_{L_{0}}(L,\End_{\FF_{p}}^{0}(A_{\pi}))$ to be optimal is~$\End_{\FF_{p}}(A_{\pi})=\varphi(O_{L})$.
	That is,~$A_{\pi}/\FF_{p}$ is principal.
	In particular, we cannot apply optimal embeddings for computing the Lie type of~$A_{\pi}$ when it comes to the Weil-$3$ number~$\sqrt{3}\zeta_{3}$, see \propref{rem:p=3}.
\end{ex}



We have the following lemma, which enables us to compute the Lie type of~$A_{\pi}\otimes\overline{\FF}_{p}$.
The orders, mentioned in the lemma, are regarded as~$\ZZ[\sqrt{p}]$-orders.
\begin{lem}\label{lem:optimal}
Let~$\pi:=\sqrt{p}$, and let~$L:=\QQ(\sqrt{p}\zeta_{3})=\QQ(\sqrt{p},\sqrt{-3})$ with maximal totally real subfield~$L_{0}:=\QQ(\pi)$ for any prime~$p>0$.
Let~$A_{\pi}/\FF_{p}$ be any simple superspecial abelian surface over~$\FF_{p}$, in the isogeny class determined by~$\pi$.
Then, there exists a CM structure of~$A_{\pi}/\FF_{p}$, which induces an embedding~$O_{L}\otimes\ZZ_{p}\embed\End_{\FF_{p}}(A_{\pi})\otimes\ZZ_{p}$.
\end{lem}
\begin{pf}
As~$A_{\pi}/\FF_{p}$ is arbitrarily chosen, the arguments below hold for any simple superspecial abelian surface over~$\FF_{p}$ in the same isogeny class.
Note that~$\End_{\FF_{p}}^{0}(A_{\pi})$ is an~$L_{0}$-algebra.
We may assume that the embedding\[{\iota}: L\embed\End_{\FF_{p}}^{0}(A_{\pi})\] is~$L_{0}$-linear as in \eqref{eq:CMreal}. 

Take $A:=\ZZ[\sqrt{p}]$ in \defref{def:optimal}.
	Let~$\pi':=\sqrt{p}\zeta_{3}$, and let~$A_{\pi'}$ be the simple superspecial abelian surface over~$\FF_{p}$, whose endomorphism ring equals~$O_{L}$, as in the proof of \lemref{lem:realRRC}.
	That is, we have \begin{equation}\label{eq:sp}
		O_{L}=\End_{\FF_{p}}(A_{\pi'}).
	\end{equation}
	Let~$\mathfrak{I}$ be the isomorphism as in the proof of \lemref{lem:realRRC}. 
	Note that\[\mathfrak{I}:\End_{\FF_{p^3}}^{0}(A_{\pi'}\otimes\FF_{p^3}) \cong \End_{\FF_{p}}^{0}(A_{\pi})\] is an isomorphism of quaternion algebras over~$L_{0}$, which are ramified only at two infinite primes.
	Let~$v$ be the place in~$L_{0}$ lying above~$p$.
	Since~$v$ is ramified in~$L_{0}$, the completion~$L_{0,v}$ of~$L_{0}$ at~$v$ is a field.
	Therefore,~$\End_{\FF_{p}}^{0}(A_{\pi})\otimes_{L_{0}} L_{0,v}$ is isomorphic to the~$2\times 2$ matrix algebra~$\textrm{M}_{2}(L_{0,v})$ over~$L_{0,v}$ by an isomorphism~$\mathcal{I}$.
	That is,\[\mathcal{I}:\End_{\FF_{p}}^{0}(A_{\pi})\otimes_{L_{0}}L_{0,v}\cong \textrm{M}_{2}(L_{0,v}).\]
	Note that~$\textrm{M}_{2}(L_{0,v})$ has the unique maximal order~$\textrm{M}_{2}(O_{L_{0,v}})$, up to isomorphism.
	Moreover, we have an isomorphism~$\mathrm{I}:\End_{\FF_{p}}^{0}(A_{\pi})\otimes_{L_{0}} L_{0,v}\cong\End_{\FF_{p}}^{0}(A_{\pi})\otimes_{L_{0}} L_{0}\otimes_{\QQ}\QQ_{p}$ because~$L_{0}$ is ramified at~$v$.
	Therefore, \[\mathrm{I}: \End_{\FF_{p}}^{0}(A_{\pi})\otimes_{L_{0}} L_{0,v} \cong \End_{\FF_{p}}^{0}(A_{\pi})\otimes_{\QQ} \QQ_{p}.\]

	Note that~$\End_{\FF_{p}}(A_{\pi})\otimes_{\ZZ}\ZZ_{p}$ (resp.~$\End_{\FF_{p^3}}(A_{\pi'}\otimes\FF_{p^3})\otimes_{\ZZ}\ZZ_{p}$) is maximal in~$\End_{\FF_{p}}^{0}(A_{\pi})\otimes_{\QQ}\QQ_{p}$ (resp.~$\End_{\FF_{p^3}}(A_{\pi'}\otimes\FF_{p^3})\otimes_{\ZZ}\QQ_{p}$) by \cite[Theorem~6.2]{W69}. 
	Therefore, we have~$\mathcal{I}\circ \mathrm{I}^{-1}(\End_{\FF_{p}}(A_{\pi})\otimes\ZZ_{p})=u_{1}\textrm{M}_{2}(O_{L_{0,v}})u_{1}^{-1}$ and~$\mathcal{I}\circ\mathrm{I}^{-1}\circ(\mathfrak{I}\otimes\textrm{id}_{\QQ_{p}})(\End_{\FF_{p^3}}(A_{\pi'}\otimes\FF_{p^3})\otimes\ZZ_{p})=u_{2}\textrm{M}_{2}(O_{L_{0,v}})u_{2}^{-1}$, where~$u_{1},u_{2}\in\textrm{GL}_{2}(O_{L_{0,v}})$, by \cite[Theorem~(17.3)]{R75}.
	Recall the construction of~${\iota}$ from \eqref{eq:CMreal}.
	Then,~${\iota}$ induces an~$L_{0,v}$-linear embedding \begin{equation}\label{eq:idQp}
		{\iota}\otimes\textrm{id}_{\QQ_{p}}: L\otimes\QQ_{p}\xrightarrow{\textrm{incl}\otimes\textrm{id}_{\QQ_{p}}}\End_{\FF_{p^3}}^{0}(A_{\pi'}\otimes\FF_{p^3})\otimes\QQ_{p}\overset{\mathfrak{I}\otimes\textrm{id}_{\QQ_{p}}}{\cong}\End_{\FF_{p}}^{0}(A_{\pi})\otimes\QQ_{p},\end{equation}where~$R_{\textrm{sp}}\otimes\ZZ_{p}$ is mapped to~$\End_{\FF_{p^3}}(A_{\pi'})\otimes\ZZ_{p}$ by~$\textrm{incl}\otimes\textrm{id}_{\QQ_{p}}$ because of \eqref{eq:sp}.
	We consider the following map\begin{equation}\label{eq:matrix}
		\End_{\FF_{p}}^{0}(A_{\pi})\otimes_{\QQ}\QQ_{p} \xrightarrow{\mathcal{I}\circ\mathrm{I}^{-1}} \textrm{M}_{2}(L_{0,v}) \xrightarrow{\textrm{c}_{12}} \textrm{M}_{2}(L_{0,v}) \xrightarrow{\mathrm{I}\circ \mathcal{I}^{-1}} \End_{\FF_{p}}^{0}(A_{\pi})\otimes_{\QQ}\QQ_{p},\end{equation} where~$\textrm{c}_{12}$ is the conjugation by~$u_{1}u_{2}^{-1}$.
	Then, the composition of \eqref{eq:idQp} and \eqref{eq:matrix}, which we denote again by~${\iota}$, is a CM structure of~$A_{\pi}/\FF_{p}$ such that~${\iota}(O_{L}\otimes\ZZ_{p})=\End_{\FF_{p}}(A_{\pi})\otimes_{\ZZ}\ZZ_{p}$.
	This completes the proof.\qed	
\end{pf}
\begin{rmk}\label{rmk:sameLietypereal}
	With the similar argument, \propref{lem:sameLietype} still holds in this case. 
\end{rmk}

After \lemref{lem:optimal}, we can consider the Lie type of~$A_{\pi}/\FF_{p}$.
Now, we state the analogous result to the case where the Weil-$p$ number is~$\sqrt{p}$ for any prime~$p>0$.
\begin{thm}\label{thm:liftingreal}
Let~$L:=\QQ(\sqrt{p}\zeta_{3})$.
Then, any simple superspecial abelian surface over~$\FF_{p}$ in the isogeny class determined by~$\pi:=\sqrt{p}$ has a CML.
\end{thm}
\begin{pf}
Let~$\pi:=\sqrt{p}$, and let~$A_{\pi}/\FF_{p}$ be a simple superspecial abelian surface in the isogeny class determined by~$\pi$.
Let~$L:=\QQ(\sqrt{p}\zeta_{3})$.
Let~$\tilde{\iota}:L\embed\End_{\FF_{p}}^{0}(A_{\pi})$ be a CM structure of~$A_{\pi}/\FF_{p}$ as in 
\lemref{lem:optimal}. 
Then, we see that~$A_{\pi}/\FF_{p}$ admits smCM by CM field~$L$, and it satisfies RRC with some CM type of~$L$ by \lemref{lem:realRRC}. 
Moreover, there exists a simple superspecial abelian surface~$\widehat{A}/\FF_{p}$, admitting smCM by~$L$ and having a CM lifting, and which is~$L$-linearly~$\FF_{p}$-isogenous to~$A_{\pi}$ by \cite[Theorem~2.5.3]{CCO14}.   

We want to prove the following three cases:
\begin{enumerate}
	\item When~$p\not\equiv 5,11$ mod~$12$ and~$p\neq 2$, any simple superspecial abelian surface over~$\FF_{p}$ has a CML, admitting smCM by~$L$;
	\item When~$p\equiv 5,11$ mod~$12$ or~$p=3$, any simple superspecial abelian surface over~$\FF_{p}$ has a CML, admitting smCM by~$L$; 
	\item When~$p=2$, any simple superspecial abelian surface over~$\FF_{2}$ has a CML, admitting smCM by~$L$.
	\end{enumerate}

Let~$L_{0}$ be the maximal totally real subfield of~$L$.
Note that the Lie type of any simple superspecial abelian surface over~$\FF_{p}$ in the isogeny class determined by~$\pi$ is the same, see \remref{rmk:sameLietypereal}.
\begin{enumerate}
	\item In this case,~$L/L_{0}$ splits completely.
So,~$A_{\pi}/\FF_{p}$ has a CML, of which the proof is similar to the one of \thmref{thm:ApisCML}.
	\item In this case,~$L/L_{0}$ is inert.
The remainder of the proof is similar the proof of \thmref{cor:surfacesCML}.
\item 
	We know that~$A_{\pi}/\FF_{2}$ satisfies RRC, which implies that there exists an abelian surface~$\widehat{A}_{}$ over~$\FF_{2}$, which has a CML with CM structure~$\iota_{}$, such that~$\widehat{A}_{}$ is~$L$-linearly~$\FF_{2}$-isogenous to~$A_{\pi}$ by \cite[Theorem~2.5.3]{CCO14}.
	
	By \cite[Theorem~6.2]{W69}, any endomorphism ring in the isogeny class determined by~$\pi:=\sqrt{2}$ is maximal in~$\End_{\FF_{2}}^{0}(A_{\pi})$ as~$p=2$, and the number of~$\FF_{2}$-isomorphism classes in the isogeny class equals the class number of~$\End_{\FF_{2}}^{0}(A_{\pi})$.
	So,~$\widehat{A}_{}$ is~$\FF_{2}$-isomorphic to~$A_{\pi}$ by an isomorphism~$f$ because the class number of~$\End_{\FF_{2}}^{0}(A_{\pi})$ is~$1$ by \cite[Section~6.2.7]{XYY19}.
	Therefore, we see that~$A_{\pi}/\FF_{2}$ has a CML by \cite[Proposition~2.22]{BKM}.

\end{enumerate}
	This completes the proof.\qed
\end{pf}

\section*{Acknowledgments}
The author is supported by the consortium, Rational points: new dimensions, through the grant OCENW.XL21.XL21.011.
The author thanks Valentijn Karemaker for sharing this beautiful topic, weekly inspirational discussions, and helpful comments on earlier drafts.
The author thanks Chia-Fu Yu for his valuable advices and comments on the draft.
The author thanks Ching-Li Chai for his generous suggestions.
The author thanks Thomas Agugliaro, Jonas Bergstr\"om and Steven Groen for the careful readings and useful suggestions.

\bibliography{cite.bib} \bibliographystyle{alpha}
\end{document}